\documentclass[11pt]{amsart}
\usepackage{latexsym,amsfonts,amssymb,amsthm,amsmath,mathrsfs,color,amscd,graphicx,fullpage,parskip,hyperref,bm,bbm,dsfont}
\usepackage[T1]{fontenc}

\usepackage{comment}
\includecomment{uncomment}
\usepackage{commath,mathtools,setspace}
\usepackage{paralist}
\usepackage{extarrows}

\excludecomment{comment} 

\newcommand{\s}[1]{{\mathcal #1}}
\newcommand{\sr}[1]{{\mathscr #1}}
\newcommand{\bb}[1]{{\mathbb #1}}
\newcommand{\vd}[2]{\dfrac{\delta #1}{\delta #2}}
\newcommand{\floor}[1]{\left\lfloor #1 \right\rfloor}

\newcommand{\dmpd}[3]{\dfrac{\partial^2 #1}{\partial #2 \partial #3}}

\newcommand{\dmder}[2]{\dfrac{\delta #1}{\delta #2}}
\newcommand{\ip}[2]{\left\langle #1,#2 \right\rangle}

\newcommand{\firststep}{\setcounter{step}{1}\textbf{Step \arabic{step}:} }
\newcommand{\nextstep}{\stepcounter{step}\textbf{Step \arabic{step}:} }
\newcommand{\R}{\mathbb{R}}

\DeclareMathOperator{\argmax}{argmax}
\DeclareMathOperator{\Lip}{Lip}

\newtheorem{theorem}{Theorem} 
\newtheorem{corollary}[theorem]{Corollary}

\newtheorem{lemma}[theorem]{Lemma}
\newtheorem{proposition}[theorem]{Proposition}

\newtheorem{definition}[theorem]{Definition}
\newtheorem{example}[theorem]{Example}
\newtheorem{remark}[theorem]{Remark}
\newtheorem{assumption}[theorem]{Assumption}

\numberwithin{equation}{section}
\numberwithin{theorem}{section}

\newcounter{step}
\mathtoolsset{showonlyrefs=true}

\begin{document}

	\title[MFG Control with Absorption]
	{Master Equation for Cournot Mean Field Games of Control with Absorption}
	
	\author{P. Jameson Graber}
	\thanks{Jameson Graber gratefully acknowledges support from the National Science Foundation through NSF CAREER Award 2045027 and NSF Grant DMS-1905449.}
	\address{J. Graber: Baylor University, Department of Mathematics;\\
		Sid Richardson Building\\
		1410 S.~4th Street\\
		Waco, TX 76706\\
		Tel.: +1-254-710- \\
		Fax: +1-254-710-3569 
	}
	\email{Jameson\_Graber@baylor.edu}
	
	\author{Ronnie Sircar}
	\thanks{Ronnie Sircar gratefully acknowledges support from the National Science Foundation through award DMS-1736409.}
	\address{R. Sircar: Princeton University, ORFE Department\\
	Sherrerd Hall 208\\
	Princeton NJ 08544.}
	\email{sircar@princeton.edu}
	
	\date{\today}   
	
	\begin{abstract}
		We establish the existence and uniqueness of a solution to the master equation for a mean field game of controls with absorption.
		The mean field game arises as a continuum limit of a dynamic game of exhaustible resources modeling Cournot competition between producers.
		The proof relies on an analysis of a forward-backward system of nonlocal Hamilton-Jacobi/Fokker-Planck equations with Dirichlet boundary conditions.
		In particular, we establish new a priori estimates to prove that solutions are differentiable with respect to the initial measure.
	\end{abstract}
	
	\keywords{mean field games, exhaustible resources}
	
	\maketitle
	
	\tableofcontents
	
	\section{Introduction} \label{sec:intro}
	
	In \cite{harris2010games}, the authors introduced 
	a dynamic game of exhaustible resource production modeling Cournot competition between producers of a good in finite supply, for instance oil, whose Markov perfect (Nash) equilibrium was characterized there by a system of coupled nonlinear PDEs. This built on the influential continuous-time study of the monopoly (single-player) version of the problem by Hotelling from 1931 \cite{hotelling1931economics}. By Cournot competition, we mean that the decision or control variable of the players is their quantity (or rate) of production, the market price or prices of the goods being determined by a decreasing function of the aggregate (or average) production. 
	
	When the goods each player produces are homogeneous, there is a single price $p$ of the good which depends, in the Cournot framework, on the average $\frac1N\sum_j^Nq_j$, where $q_j\geq0$ is player $j$'s quantity, and there are $N<\infty$ players. When the goods are substitutable, for instance oil of different grades from different suppliers, or consumer goods such as televisions, 	
	a typical model has that the price $p_i$ that producer $i$ receives for its good depends in a decreasing manner on $q_i+\frac{\epsilon}{N-1}\sum_{j \neq i} q_j$. That is, its price is influenced by the average of the other players' quantities (thereby viewing them as exchangeable), where $\epsilon>0$ measures the degree of interaction. A dynamic exhaustible resources problem in this case is analyzed in \cite{ledvina2012oligopoly}. 
	
Mean field games, in which there is a continuum of players, have been much-studied in the past 15 years. We refer, for instance, to \cite{bensoussan2013mean} and \cite{carmona2017probabilistic}, for surveys from PDE and probabilistic perspectives respectively.  In the context of the Cournot model, the homogeneous goods case leads to a continuum approximation model whose optimal strategies are of (unrealistic) bang-bang type: the players either produce nothing or as quickly as possible. The substitutable goods case has a more reasonable mean field game model, as studied in \cite{chan2015bertrand} and \cite{chan2017fracking}. As mean field games of control, and because the state variable is absorbed at zero (exhaustion of the resource), they differ from the vast majority of problems studied in the literature where interaction is through the mean of the state variable, which lives on the full space. Rigorous existence results are thus more recent and under various restrictions, for instance \cite{graber2018existence,graber2018variational,graber2020commodities,graber2021nonlocal}.
We refer the reader to \cite{cardaliaguet2017mfgcontrols,gomes2014extended,gomes2016extended,kobeissi2021classical} for benchmark results on mean field games of controls.

There has been much recent interest in describing mean field games through a Master Equation \cite{cardaliaguet2019master,bensoussan2015master,carmona2017probabilisticII}.
The study of such equations now has a large body of literature, going back to such works as \cite{gangbo2015existence,chassagneux2019numerical}.
Again the existing results in the literature concern mean field interaction through the state.
See the recent results found in~\cite{gangbo2020global,gangbo2021mean,mayorga2019short,mou2019wellposedness}.
As for boundary conditions, most references contain results only for master equations on the whole space or with periodic boundary conditions.
See, however, the recent work by Ricciardi for Neumann boundary conditions \cite{ricciardi2021master}. 
Here we introduce and analyze the Master Equation of Cournot mean field games of control with absorption. 
Our main result is the existence and uniqueness of a classical solution.

Once one has a unique classical solution to the master equation, a natural application is to the convergence problem for $N$-player games corresponding to a mean field game.
Using the arguments of \cite[Chapter 6]{cardaliaguet2019master}, one can hope to obtain estimates that prove the closed-loop Nash equilibrium strategies for $N$-player games converge to the mean field equilibrium strategy.
In our case, the infinite time horizon, the dependence of the dynamics on the distribution of controls, and especially the absorbing boundary conditions add technical obstacles to a straightforward application of the arguments found in \cite{cardaliaguet2019master}.
We leave this application to future research.

In the rest of this section, we introduce the main notation needed and give our main results. In Section \ref{sec:description}, we give the precise description of the Cournot model as a mean field game and write the corresponding master equation.
In Section \ref{sec:metric}, we define a metric on the space of measures and introduce a notion of derivative for functions defined on this space.
In Section \ref{sec:main result}, we give the definition of a solution to the master equation and present Theorem \ref{thm:main result}, which gives precise conditions under which a unique solution exists.
Finally, in Section \ref{sec:structure} we present the outline of the rest of the paper, which is devoted to the proof of Theorem \ref{thm:main result}.

	
	\subsection{Description of the model} \label{sec:description}
	Let 
	$P:\intco{0,\infty} \to \bb{R}$ be a given price function, satisfying the following:
	\begin{assumption}
		\label{as:P}
		$P$ is continuous on $\intco{0,\infty}$ with $P(0) > 0$.
		For some $n \geq 4$,
		$P$ is $n$ times continuously differentiable on $(0,\infty)$, $P^{(n)}$ is locally Lipschitz, and $P' < 0$.
		In addition, $\limsup_{q \to 0+} P'(q)$ is strictly less than zero (it could be $-\infty$), and there exists a finite \emph{saturation point} $\eta > 0$ such that $P(\eta) = 0$.
	\end{assumption}
	The \emph{profit function} $\pi : \intco{0,\infty}^4 \to \bb{R}$ for an individual producer is given by
	\begin{equation} \label{eq:revenue}
		\pi(\epsilon,q,Q,a) = \begin{cases}
			q\del{P(\epsilon Q+q)-a} &\text{if}~q > 0,\\
			0 &\text{if}~q = 0.
		\end{cases}
	\end{equation}
	Here $q$ is the rate of production chosen by the producer, $Q$ is the market's aggregate rate of production, $a$ is the marginal cost of production, and $\epsilon \geq 0$ is a fixed parameter that determines the substitutability of goods.
	
	It will be convenient to define the \emph{relative prudence}
	\begin{equation}
		\rho(Q) := -\frac{QP''(Q)}{P'(Q)}
	\end{equation}
	Notice that by Assumption \ref{as:P}, $\rho$ is continuously differentiable on $(0,\infty)$.
	If, for example, we take $P'(q) = -q^{-\rho}$ for some fixed $\rho \in\R$
	(cf.~\cite{harris2010games}), then $\rho(Q) = \rho$ (constant relative prudence).
	\begin{assumption}[Relative prudence]
		\label{as:prudence}
		We assume
		\begin{equation}
			\bar \rho := \sup_{Q \in (0,\infty)} \rho(Q) < \frac{2+\epsilon}{1+\epsilon} \leq 2.
		\end{equation}
	\end{assumption}
Assumptions \ref{as:P} and \ref{as:prudence} guarantee a Hamiltonian of the following continuous time game is well-defined.
	
	In the finite $N$-player differential game introduced in \cite{harris2010games}, each player $i$ has remaining stock (or reserves) $x_i(t)$ at time $t\geq0$
	and we denote by $\bar q_i(t)\geq0$ their chosen rate of production, so $x_i(t)$ satisfies the stochastic differential equation 
	\begin{equation}
		\dif{x}_i(t) = \left(-\bar q_i(t) \dif{t}+ \sigma \dif W_i(t)\right)\bb{I}_{\cbr{x_i(t) > 0}},
	\end{equation}
	where each $W_i(t)$ is an independent standard Brownian motion representing, for instance, uncertainty in the extraction process.
	The producers start with initial ($t=0$) reserves ${\bf x}\in\bb{R}^N_+$
and each
	maximizes expected discounted lifetime profit. The value function $u_i:\bb{R}^N_+\to\bb{R}$ of player $i$ is given by 
	\begin{equation}
		u_i({\bf x})= \sup_{\bar q_i} \bb{E}\int_0^{\tau_i} e^{-rt}\pi\del{\epsilon,\bar q_i(t),\bar Q_{-i}(t), 0}
		\dif t, \label{uidef}
	\end{equation}
	where 
	$\tau_i$ is the first time $x_i$ hits (and is absorbed at) zero, $r\geq0$ is the common discount rate on future profits, $\bar Q_{-i}(t)$ is the mean production rate of the other producers:
	$$ \bar Q_{-i}(t) = \frac{1}{N-1}\sum_{j\neq i}\bar q_j(t), $$
	 and we assume for simplicity that marginal costs of production are zero.

The Hamilton-Jacobi-Bellman equation corresponding to each player's optimal control problem in \eqref{uidef} is as follows.
	Define 
	\begin{equation*}
		H(\epsilon,Q,a) := \sup_{q \geq 0} \pi(\epsilon,q,Q,a)
		\quad
		\mbox{from which it follows}
		\quad 
		\argmax_{q \geq 0} \pi(\epsilon,q,Q,a) = -\pd{H}{a}(\epsilon,Q,a).
	\end{equation*}
In a Markov perfect (Nash) equilibrium of the $N$-player differential game
the associated system of Hamilton-Jacobi-Bellman (HJB) partial differential equations (PDEs) for the value functions is
	\begin{equation} \label{eq:hjb}
	H\del{\epsilon,\bar Q_{-i}^*({\bf x}),\dpd{u_i}{x_i}} + \sum_{j \neq i} \dpd{H}{a}\del{\epsilon,\bar Q_{-j}^*({\bf x}),\dpd{u_j}{x_j}}\dpd{u_i}{x_j} - ru_i  + \frac{\sigma^2}{2}\sum_{j=1}^N \dpd[2]{u_i}{x_j} = 0,
	\end{equation}
	coupled with 
	\begin{equation} \label{eq:total quantity}
	\bar Q_{-i}^*({\bf x}) = -\frac{1}{N-1}\sum_{j \neq i}\dpd{H}{a}\del{\epsilon,\bar Q_{-j}^*({\bf x}),\dpd{u_j}{x_j}}.
	\end{equation}
	See \cite[Equation (3.4)]{harris2010games}; here we have additional diffusion terms due to the Brownian noise in the dynamics.
	
The mean field game (MFG) version of this problem, corresponding to a continuum of players 
with density of initial reserves $m$ 
was introduced in \cite{chan2015bertrand} and further studied in \cite{chan2017fracking}, where it is characterized by two PDEs and a fixed point condition (which are given here in Section \ref{sec:structure}). An explicit solution of the deterministic MFG ($\sigma=0$) when the price function $P$ is linear is given in \cite{PUR}.

We next introduce the master equation formulation of this MFG.
	 
	\subsection{Master Equation Heuristics}
		Let $m$ be a measure representing the initial distribution of stock over all producers.
	Let $U(x,m)$ be the maximum discounted lifetime profit for an individual producer that starts with a stock of $x$.
	If we assume that $U$ is smooth with respect to both variables (see Definition \ref{def:derivative} below for derivatives in the space of measures), then $U$ will satisfy
	\begin{multline} \label{eq:master equation}
		H\del{\epsilon, Q^*,\dpd{U}{x}(m,x)} + \int_{\s{D}} \dpd{H}{a}\del{\epsilon, Q^*,\dpd{U}{x}(m,y)}\dpd{}{y}\vd{U}{m}(m,x,y)\dif m(y) - rU(m,x)\\
		 +\dfrac{\sigma^2}{2}\del{\dpd[2]{U}{x}(m,x) + \int_{\s{D}} \dpd[2]{}{y}\vd{U}{m}(m,x,y)
			\dif m(y)} = 0,
	\end{multline}
	where $Q^*$ is defined as the unique fixed solution of the equation
	\begin{equation} \label{eq:clearing condition}
		Q^* = -\int_{\s{D}}\dpd{H}{a}\del{\epsilon,Q^*,\dpd{U}{x}(m,y)} \dif m(y).
	\end{equation}
	Equation \eqref{eq:master equation} is called the \emph{master equation}.
	
	\begin{comment}
	 and 
	$a_1,\ldots,a_N \geq 0$ given constant (marginal) costs.
	The profit function is
	
	where $Q_{-i} := \frac{1}{N-1}\sum_{j \neq i} q_j$ and $0 < \epsilon \leq N-1$.
	(If $\epsilon = N-1$ then $P$ depends precisely on the total demand; this is the homogeneous goods case.)
	
	A \emph{Nash equilibrium} is a vector ${\bf q^*} = (q_1^*,\ldots,q_N^*) \in [0,\infty)^N$ such that, for all $i$,
	\begin{equation}
	\pi(\epsilon,q_i^*,Q_{-i}^*,a_i) = \max_{q_i \in [0,\infty)} \pi(\epsilon,q_i,Q_{-i}^*,a_i).
	\end{equation}
	
	Conditions are given in \cite[Proposition 2.5]{harris2010games} 
	and \cite{ledvina2012oligopoly} that guarantee a unique Nash equilibrium for every ${\bf a} = (a_1,\ldots,a_N) \in [0,\infty)^N$.
	In this case we denote the equilibrium by ${\bf q^*}({\bf a})$, and we define $G_i({\bf a}) = \pi(\epsilon,q_i^*,Q_{-i}^*,a_i)$.
	\end{comment}
	
	Formally, the master equation can be derived from the system of Hamilton-Jacobi-Bellman (HJB) equations \eqref{eq:hjb} for the $N$-player game.
	Letting $N \to \infty$, we formally interpret each sum as an integral with respect to the population distribution.
	See \cite{cardaliaguet2019master,carmona2017probabilisticII} for a detailed interpretation of the master equation.

	\subsection{Metric and derivative on a space of measures}
	\label{sec:metric}
	
	Before we can state our main result, we will need to define a notion of derivative with respect to a measure.
	Let $\s{M} = \s{M}(\s{D})$ be the space of all finite signed Radon measures $\mu$ on $\s{D}$.
	We denote by $\s{M}_{+}$ the subset of $\s{M}$ consisting only of positive measures.
	The topology on $\s{M}$ is that of narrow convergence.
	We say that a sequence $\{\mu_n\}$ in $\s{M}$ converges narrowly if for every bounded continuous function $\phi$ on $\s{D}$, we have
	\begin{equation*}
		\int_{\s{D}} \phi(x)\dif \mu_n(x) \to \int_{\s{D}} \phi(x)\dif \mu(x).
	\end{equation*}

	We now introduce the derivative on $\s{M}(\s{D})$.
	\begin{definition}[Differentiability with respect to measures]
		\label{def:derivative}
		Let $\sr{M}$ be any dense subset of $\s{M}_+$.
		Given a function $F:\sr{M} \to \bb{R}$,
		we say that $F$ is continuously differentiable if there exists a continuous function $f:\sr{M} \times \s{D} \to \bb{R}$, satisfying
		\begin{equation*}
			\abs{f(m,x)} \leq C(m) \ \forall x \in \s{D}
		\end{equation*}
		for some constant $C(m)$, such that
		\begin{equation} \label{eq:derivative}
			\lim_{t \to 0+} \frac{F\del{m + t(\hat m - m)}- F(m)}{t} = \int_{\s{D}} f(m,x)\dif\,(\hat m - m)(x) \quad\forall m,\hat m \in \sr{M}.
		\end{equation}
		The function $f(m,x)$ is unique, and we denote it $f(m,x) = \dmder{F}{m}(m,x)$.
	\end{definition}
	Definition \ref{def:derivative} is essentially the classical G\^ateaux derivative, though we only take $m,\hat m$ from the convex subset $\sr{M}$ of the vector space $\s{M}$.
	Uniqueness follows from the fact that the measure $\hat m -m$ in \eqref{eq:derivative} can be taken to be an essentially arbitrary positive measure (by density of $\sr{M}$ in $\s{M}_+$); contrast with the situation in which $m,\hat m$ must be probability measures (cf.~\cite{cardaliaguet2019master}).
	
	\subsection{Statement of the main result}
	\label{sec:main result}
	
	To state our main result, we will first define a set of measures on which the master equation \eqref{eq:master equation} is supposed to hold.
	Fix $\alpha \in (0,1)$ and let $\s{M}^{2+\alpha}$ denote the set of all positive measures $m$ on $\s{D} = (0,\infty)$ satisfying the condition
	\begin{equation*}
		\int_{\s{D}} x^{-2-\alpha}\dif m(x) < \infty.
	\end{equation*}
	\begin{definition}
		We say that a function $U:\s{D} \times \s{M}^{2+\alpha} \to \bb{R}$ is a (classical) solution of the master equation \eqref{eq:master equation}-\eqref{eq:clearing condition} with absorbing boundary conditions provided it satisfies the following:
		\begin{enumerate}
			\item $U(0,m) = 0$ for every $m \in \s{M}^{2+\alpha}$;
			\item $U$ and $\vd{U}{m}$ are twice continuously differentiable with respect to $x$;
			\item for every $m \in \s{M}^{2+\alpha}$ and $x > 0$,
			Equation \eqref{eq:master equation} is satisfied.
		\end{enumerate}
	\end{definition}
	The Dirichlet boundary condition $U(0,m) = 0$ is an absorbing type boundary condition, representing the fact that players exit the game as they run out of resources (cf.~\cite{chan2015bertrand,hambly2017stochastic}).
	Theorem \ref{thm:main result} is the first result, as far as we know, on the Master Equation with boundary conditions of this type.
	
	Our main result in this paper is as follows.
	\begin{theorem} \label{thm:main result}
		Under Assumptions \ref{as:P} and \ref{as:prudence}, the following assertions hold.
		\begin{enumerate}
			\item There exist constants $r^* > 0$ (large) and $\epsilon^* > 0$ (small) such that whenever $r \geq r^*$ and $0 < \epsilon \leq \epsilon^*$, the master equation \eqref{eq:master equation} has a solution, which is unique under the condition \eqref{eq:dUdm estimate} (cf.~Section \ref{sec:holder estimates}).
			\item If $P$ is linear, and in particular if (without loss of generality) $P(q) = 1-q$, then there exists a constant $r^*$ such that for every $r \geq r^*$ and $\epsilon < 2$, the master equation \eqref{eq:master equation} has a solution, which is unique under the condition \eqref{eq:dUdm estimate} (cf.~Section \ref{sec:holder estimates}).
		\end{enumerate}
	\end{theorem}
	\begin{remark}
		The precise conditions on $r^*$ and $\epsilon^*$ in Theorem \ref{thm:main result} are contained in Assumptions \ref{as:r big, epsilon small} and \ref{as:r big, P linear}.
		Although these two conditions are essentially in dichotomy, nevertheless in this paper we make an attempt to utilize as much as possible a unified method of proof for both cases.
		See Remark \ref{rem:why ep small} for more details.
	\end{remark}
	
	\subsection{Structure of the proof} 
	\label{sec:structure}
	
	In a generalized sense, we use the method of characteristics to solve the master equation \eqref{eq:master equation}-\eqref{eq:clearing condition}.
	Consider the HJB/Fokker-Planck system
	\begin{equation} \label{eq:mfg infty}
		\begin{cases}
			(i) & \dpd{u}{t} + \dfrac{\sigma^2}{2}\dpd[2]{u}{x} + H\del{\epsilon,Q^*(t),\dpd{u}{x}} - ru = 0,\\
			(ii) & \dpd{m}{t} - \dfrac{\sigma^2}{2}\dpd[2]{m}{x} + \dpd{}{x}\del{\dpd{H}{a}\del{\epsilon,Q^*(t),\dpd{u}{x}}m} = 0,\\
			(iii) & Q^*(t) = -\displaystyle \int_{\s{D}} \dpd{H}{a}\del{\epsilon,Q^*(t),\dpd{u}{x}}\dif m(t),\\
			(iv) & m|_{x=0} = u|_{x=0} = 0, \
			m|_{t=0} = m_0 \in \s{M}_+(\s{D})
		\end{cases}
	\end{equation}
	where $\s{D} := (0,\infty)$.
	We can think of System \eqref{eq:mfg infty} as the characteristics of Equation \eqref{eq:master equation}.
	Indeed, suppose $U$ is a smooth solution to \eqref{eq:master equation} and $(u,m)$ is a smooth solution to \eqref{eq:mfg infty}.
	Then formally the two are related by the equation $u(x,t) = U(x,m(t))$, and in particular $U(x,m_0) = u(x,0)$.
	In the proof of our main result, our strategy will be to \emph{define} a function $U$ in this way, then prove that it satisfies \eqref{eq:master equation}.
	To do this, we follow these steps:
	\begin{enumerate}
		\item Prove that \eqref{eq:mfg infty} has a unique solution $(u,m)$ for any $m_0 \in \s{M}^\alpha$.
		Define $U(x,m_0) = u(x,0)$.
		\item Prove that $U$ is differentiable with respect to the measure variable $m_0$:
		\begin{enumerate}
			\item Formally differentiate \eqref{eq:mfg infty} with respect to the measure variable to obtain a linearized system.
			\item Prove that the linearized system has a unique solution. 
			\item Prove that the unique solution thus obtained is indeed the derivative of $U$ with respect to the measure.
		\end{enumerate}
		\item Use the smoothness of $U$ to establish that System \eqref{eq:master equation}-\eqref{eq:clearing condition} is satisfied.
	\end{enumerate}
	The remainder of this paper is structured as follows.
	In Section \ref{sec:prelim} we establish notation and define function spaces as needed.
	In Section \ref{sec:fokkerplanck} we study the Fokker-Planck equation with absorbing boundary conditions and establish some results that allow us to prove existence of solutions to System \eqref{eq:mfg infty}; they may also have independent interest.
	In Section \ref{sec:fwdbckwd} we present existence, uniqueness, and regularity results on System \eqref{eq:mfg infty}.
	Section \ref{sec:sensitivity} is the core this paper, in which we derive all of the a priori estimates on linearized systems that will allow us to prove differentiability of the master field $U(x,m)$.
	Here the reader will find some parallels with a recent work by Graber and Laurel that also deals with linearized systems in order to analyze sensitivity of solutions to the parameter $\epsilon$ \cite{graber2022parameter}.
	In the present work, the analysis is considerably more sophisticated because we are taking derivatives with respect to a \emph{measure} and not a scalar parameter; this requires estimates on a linearized system in appropriate norms, in particular dual spaces that introduce a great deal of technicalities.
	The main result is proved in Section \ref{sec:solution master}, essentially as a corollary of Section \ref{sec:sensitivity}.
	Proofs of some technical results are left in the appendix.
	
	\section{Preliminaries}
	\label{sec:prelim}
	
	\subsection{Function spaces}
	Let $\s{D} = (0,\infty)$.
	For $n \in \bb{N}$, we denote by $\s{C}^n = \s{C}^n(\overline{\s{D}})$ the space of all $n$ times continuously differentiable functions on $\overline{\s{D}}$ such that the norm
	\begin{equation*}
	\enVert{f}_{\s{C}^n(\overline{\s{D}})} = \sum_{k=0}^n \sup_{x \in \overline{\s{D}}} \abs{\dod[k]{f}{x}(x)}
	\end{equation*}
	is finite; $\s{C}^n(\overline{\s{D}})$ is a Banach space endowed with this norm.
	In particular, $\s{C}^0(\overline{\s{D}})$ is simply the space of all continuous functions, endowed with the supremum norm.
	We denote by $\s{C}^n_c = \s{C}_c^n(\s{D})$ the space of all $n$ times continuously differentiable functions which have compact support contained in $\s{D}$; this is a subspace of $\s{C}^n(\overline{\s{D}})$, and $\s{C}^n_0(\s{D})$ denotes its closure.
	We also denote $\s{C}_c^\infty(\s{D}) = \cap_{n=1}^\infty \s{C}_c^n(\s{D})$.
	
	For any $\alpha \in (0,1)$, define the H\"older seminorm
	\begin{equation*}
	[f]_\alpha := \sup_{x,y \in \overline{\s{D}}, x \neq y} \frac{\abs{f(x)-f(y)}}{\abs{x-y}^\alpha}.
	\end{equation*}
	Define $\s{C}^{n+\alpha} = \s{C}^{n+\alpha}(\overline{\s{D}})$ to be the space of all $n$ times continuously differentiable functions $f$ whose $n$th derivative is H\"older continuous, such that the norm
	\begin{equation*}
	\enVert{f}_{\s{C}^{n+\alpha}(\overline{\s{D}})} = \enVert{f}_{\s{C}^n(\overline{\s{D}})} + \intcc{\dod[n]{f}{x}}_\alpha
	\end{equation*}
	is finite.
	In particular, when $n = 0$ the space $\s{C}^\alpha(\overline{\s{D}})$ is simply the space of all $\alpha$-H\"older continuous functions with standard norm.
	We define $\s{C}_\diamond^\alpha = \s{C}_\diamond^\alpha(\s{D})$ to be the space of all $f \in \s{C}^\alpha(\overline{\s{D}})$ such that $f(0) = 0$.

	When $\alpha = 1$, the quantity $[f]_\alpha$ defined above is referred to as the Lipschitz constant of $f$, denoted $\Lip(f)$ instead of $[f]_1$.
	We define $\Lip(\overline{\s{D}})$ to be the space of all Lipschitz continuous functions on $\overline{\s{D}}$, with norm
	\begin{equation*}
		\enVert{f}_{\Lip(\overline{\s{D}})} = \enVert{f}_{\s{C}^0} + \Lip(f),
	\end{equation*}
	and the subspace $\Lip_\diamond(\s{D})$ the set of all $f \in \Lip(\overline{\s{D}})$ such that $f(0) = 0$.
	
	We now define H\"older spaces of functions on space-time.
	Let $I = [0,T]$ or $I = \intco{0,\infty}$.
	For any number $\beta \geq 0$ we define the space $\s{C}^{\beta,0}(\overline{\s{D}} \times I)$ to be the set of all functions $u:\overline{\s{D}} \times I \to \bb{R}$ such that the following norm is finite:
	\begin{equation*}
		\enVert{u}_{\s{C}^{\beta,0}} = \enVert{u}_{\s{C}^{\beta,0}(\overline{\s{D}} \times I)} := \sup_{t \in I} \enVert{u(\cdot,t)}_{\s{C}^\beta(\overline{\s{D}})}.
	\end{equation*}
	For any $\alpha \in \intoo{0,1}$ define
	\begin{equation*}
	[u]_{\alpha,\alpha/2} := \sup_{x,y \in \overline{\s{D}}, t,s \in I, x \neq y, t \neq s} \frac{\abs{u(x,t)-u(y,s)}}{\abs{x-y}^\alpha + \abs{t-s}^{\alpha/2}}.
	\end{equation*}
	We denote by $\s{C}^{\alpha,\alpha/2}(\overline{\s{D}} \times I)$ the subspace of $\s{C}^{0,0}(\overline{\s{D}} \times I)$ such that the norm
	\begin{equation*}
	\enVert{u}_{\s{C}^{\alpha,\alpha/2}(\overline{\s{D}} \times I)} := \enVert{u}_{\s{C}^{0,0}(\overline{\s{D}} \times I)} + [u]_{\alpha,\alpha/2}
	\end{equation*}
	is finite.
	The space $\s{C}^{2,1}(\overline{\s{D}} \times I)$ consists of functions such that
	\begin{equation*}
	\enVert{u}_{\s{C}^{2,1}(\overline{\s{D}} \times I)} := \enVert{u}_{\s{C}^{0,0}(\overline{\s{D}} \times I)} + \enVert{\dpd{u}{x}}_{\s{C}^{0,0}(\overline{\s{D}} \times I)} + \enVert{\dpd[2]{u}{x}}_{\s{C}^{0,0}(\overline{\s{D}} \times I)} + \enVert{\dpd{u}{t}}_{\s{C}^{0,0}(\overline{\s{D}} \times I)}
	\end{equation*}
	is finite, and the subspace $\s{C}^{2+\alpha,1+\alpha/2}(\overline{\s{D}} \times I)$ such that
	\begin{equation*}
	\enVert{u}_{\s{C}^{2+\alpha,1+\alpha/2}(\overline{\s{D}} \times I)} := \enVert{u}_{\s{C}^{2,1}(\overline{\s{D}} \times I)} + \intcc{\dpd[2]{u}{x}}_{\alpha,\alpha/2} + \intcc{\dpd{u}{t}}_{\alpha,\alpha/2}
	\end{equation*}
	is finite. Cf.~\cite[Section 1.1]{ladyzhenskaia1968linear}.
	Note that there exist constants $C_\alpha$ such that
	\begin{equation*}
	 \enVert{u}_{\s{C}^{2+\alpha,1+\alpha/2}(\overline{\s{D}} \times I)} \leq C_\alpha\del{\enVert{u}_{\s{C}^{0,0}(\overline{\s{D}} \times I)} + \intcc{\dpd[2]{u}{x}}_{\alpha,\alpha/2} + \intcc{\dpd{u}{t}}_{\alpha,\alpha/2}}.
	\end{equation*}
	We define the Lebesgue spaces $L^p$ in the usual way, and we write the norms $\enVert{f}_p = \enVert{f}_{L^p}$ interchangeably.

	\subsection{Norms on the space of measures}
	We define the total variation norm $\enVert{\mu}_{TV} = \abs{\mu}(\s{D})$, which can also be expressed as
	\begin{equation*}
		\enVert{\mu}_{TV} = \sup\cbr{\int_{\s{D}} \phi(x)\dif\mu(x) : \phi \in \s{C}^0(\s{D}), \ \enVert{\phi}_{\s{C}^0} \leq 1}.
	\end{equation*}
	Under this norm, $\s{M}$ becomes a Banach space.
	On the other hand, it is not necessary to converge in this norm in order to converge narrowly.
	For this it suffices to consider $\s{M}$ as a subspace of the dual of $\s{C}^\alpha_{\diamond}$, with norm
	\begin{equation*}
		\enVert{\mu}_{\del{\s{C}^\alpha_{\diamond}}^*} = \sup\cbr{\int_{\s{D}} \phi(x)\dif\mu(x) : \phi \in \s{C}^\alpha_{\diamond}(\s{D}), \ \enVert{\phi}_{\s{C}^\alpha} \leq 1}.
	\end{equation*}
	We may also replace $\s{C}^\alpha_{\diamond}$ with $\Lip_{\diamond}$.
	\begin{lemma} \label{eq:Ca* implies narrow}
		Let $\{\mu_n\}$ be a sequence in $\s{M}$.
		If $\enVert{\mu_n}_{TV}$ is bounded, if $\enVert{\mu_n - \mu}_{\del{\s{C}^\alpha_{\diamond}}^*} \to 0$, and if $\mu_n(\s{D}) \to \mu(\s{D})$, then $\mu_n$ converges narrowly to $\mu$.
	\end{lemma}
	
	\begin{proof}
		Let $\phi$ be a bounded, continuous function on $\s{D}$, and let $\varepsilon > 0$.
		Choose $\psi \in \s{C}^\alpha_{\diamond}$ such that $\enVert{\phi - \phi(0) - \psi}_{\s{C}^0} < \varepsilon$.
		Then
		\begin{equation*}
			\abs{\int_{\s{D}} \phi \dif\,(\mu_n - \mu)} \leq \varepsilon \del{\enVert{\mu_n}_{TV} + \enVert{\mu}_{TV}} + \abs{\phi(0)}\abs{\mu_n(\s{D}) - \mu(\s{D})} + \abs{\int_{\s{D}} \psi \dif\,(\mu_n - \mu)}.
		\end{equation*}
		Using the fact that $\enVert{\mu_n}_{TV}$ is bounded, we let $n \to \infty$ and then $\varepsilon \to 0$ to conclude.
	\end{proof}

	\subsection{Remark on constants}
	
	Throughout this manuscript, $C$ will denote a generic positive constant, whose precise value may change from line to line.
	When $C$ depends on the data from the problem, will attempt to specify all the parameters on which $C$ depends.
	In particular, we may write $C(a_1,\ldots,a_n)$ to denote a positive number which depends on given parameters $a_1,\ldots,a_n$.
	When no parameters are specified, this means $C$ depends only on the number of steps in the proof (and is generally an increasing function thereof).

	\section{Fokker-Planck equation with absorbing boundary conditions}
	\label{sec:fokkerplanck}
	
	Recall $\s{D} := (0,\infty)$.
	In this section we study weak solutions to a Fokker-Planck equation with Dirichlet boundary conditions:
	\begin{equation} \label{eq:fp}
		\begin{cases}
			\dpd{m}{t} - \dfrac{\sigma^2}{2}\dpd[2]{m}{x} - \dpd{}{x}\del{bm} = 0,\\
			m|_{x=0} = 0, \
			m|_{t=0} = m_0
		\end{cases}
	\end{equation}
	for a given velocity function $b = b(x,t)$.
	We want an interpretation of \eqref{eq:fp} that makes sense for any $m_0 \in \s{M}(\s{D})$.	
	Thus we say that $m \in \s{C}^0\del{[0,T];\s{M}(\s{D})}$ is a \emph{weak solution} of \eqref{eq:fp} provided that, for all $\phi \in \s{C}_c^\infty(\s{D} \times \intco{0,T})$, we have
	\begin{equation} \label{eq:fp weak}
		\int_0^T \int_{\s{D}} \del{-\dpd{\phi}{t} - \dfrac{\sigma^2}{2}\dpd[2]{\phi}{x} + b\dpd{\phi}{x}}m(\dif x, t)\dif t = \int_{\s{D}} \phi(x,0)m_0(\dif x).
	\end{equation}

	Our main existence/uniqueness result is contained in the following lemma.
	Its proof is fairly standard and is found in Appendix \ref{ap:fokkerplanck proofs}.
	
	\begin{lemma} \label{lem:fp}
		Let $b$ be a bounded continuous function on $\s{D} \times [0,T]$, and let $m_0 \in \s{M}_{1,+}(\s{D})$.
		Then there exists a unique weak solution $m$ of \eqref{eq:fp}.
		It satisfies
		\begin{equation} \label{eq:TV}
			\enVert{m(t)}_{TV} \leq \enVert{m_0}_{TV} \quad \forall t \geq 0.
		\end{equation}
		It is also H\"older continuous with respect to the $\s{C}^\alpha_{\diamond}(\s{D})^*$ and $\Lip_{\diamond}(\s{D})^*$ metrics, and in particular
		\begin{equation} \label{eq:m holder in time}
			\begin{aligned}
				\enVert{m(t)}_{\s{C}^\alpha_{\diamond}(\s{D})^*} &\leq \enVert{m_0}_{TV}\del{\int_{\s{D}} x^\alpha m_0(\dif x) + 2\del{\enVert{b}_\infty^\alpha + \sigma^\alpha} \max \cbr{t^{\alpha},t^{\alpha/2}}}\\
				\enVert{m(t_1)-m(t_2)}_{\s{C}^\alpha_{\diamond}(\s{D})^*} &\leq 2\enVert{m_0}_{TV}\del{\enVert{b}_\infty^\alpha + \sigma^\alpha}\abs{t_1-t_2}^{\alpha/2} \ \forall t_1,t_2 \geq 0 \ \text{s.t.} \ \abs{t_1-t_2} \leq 1,\\
			\end{aligned}
		\end{equation}
		where for $\alpha = 1$ we replace $\s{C}^1_{\diamond}(\s{D})^*$ with $\Lip_{\diamond}(\s{D})^*$.
		Its total mass function $\eta(t)$ is continuous and decreasing on $[0,T]$.
	\end{lemma}
	
	Lemma \ref{lem:fp} has the following straightforward corollary, whose proof we omit.
	\begin{corollary} \label{cor:fp signed}
		Let $b$ be a bounded continuous function on $\s{D} \times [0,T]$, let $m_0 \in \s{M}_{1}(\s{D})$, and let $m_0^+$ and $m_0^-$ denote the positive and negative parts, respectively, of $m_0$.
		Then there exists a unique weak solution $m$ of \eqref{eq:fp}, whose positive part $m^+$ is precisely the solution of \eqref{eq:fp} with $m_0$ replaced by $m_0^+$, and whose negative part $m^-$ is the solution of \eqref{eq:fp} with $m_0$ replaced by $m_0^-$.
		The estimates \eqref{eq:m holder in time} still hold, with $m_0$ replaced by $\abs{m_0}$.
	\end{corollary}

	\subsection{The mass function}
	\label{sec:mass function}
	
	Let $m$ be a weak solution to \eqref{eq:fp}.
	We define the \emph{total mass function} $\eta:[0,T] \to \bb{R}$ by
	\begin{equation}
		\eta(t) := \int_{\s{D}} m(\dif x,t).
	\end{equation}
	Notice that $\eta$ is in general not constant.
	Since the equations in System \eqref{eq:mfg infty} depend on $\eta$, we are motivated to study the regularity of $\eta$ as a function of time, and in particular we would like to know when it is H\"older continuous in order to establish the existence of classical solutions to the system.
	Note that it is insufficient to know how regular it is only for $t$ away from zero, because the behavior of the population mass as $t \to 0$ influences the regularity of solutions to the backward-in-time Hamilton-Jacobi equation.
	
	As a first step, we analyze the case where $b = 0$, so that \eqref{eq:fp} reduces to the heat equation with absorbing boundary conditions.
	Our goal is to determine whether the heat semigroup itself produces a H\"older continuous flow of total population mass.
	Recall that the heat kernel is given by
	\begin{equation} \label{eq:heat kernel}
		S(x,t) = (2\sigma^2\pi t)^{-1/2}\exp\cbr{-\frac{x^2}{2\sigma^2t}}
	\end{equation}
	and that the solution of the heat equation with absorbing boundary condition at $x=0$
	\begin{equation} \label{eq:heat}
		\dpd{m}{t} = \frac{\sigma^2}{2}\dpd[2]{m}{x}, \quad m|_{t=0} = m_0, \quad m|_{x=0} = 0
	\end{equation}
	is given by
	\begin{equation} \label{eq:heat solution}
		m(x,t) = \int_{\s{D}} \del{S(x-y,t) - S(x+y,t)}m_0(\dif y).
	\end{equation}
	For a measure $m_0 \in \s{M}(\s{D})$ the corresponding mass function generated by the heat equation is
	\begin{equation}
		\label{eq:eta function}
		\eta^h[m_0](t) := \int_{\s{D}}\int_{\s{D}} \del{S(x-y,t) - S(x+y,t)}m_0(\dif y)\dif x.
	\end{equation}
	By Fubini's theorem, one can reverse the order of integration in \eqref{eq:eta function} and then write $\eta^h[m](t)$ explicitly in terms of the cdf of $m$:
	\begin{equation*}
			\eta^h[m_0](t) = \frac{2}{\sqrt{2\sigma^2\pi}}\int_0^\infty \exp\cbr{-\frac{x^2}{2\sigma^2}} m_0\del{\intoo{t^{1/2}x,\infty}}  \dif x.
	\end{equation*}
	
	To the question, ``Is $\eta^h[m_0](\cdot)$ H\"older continuous on $[0,T]$ for every measure $m_0 \in \s{M}(\s{D})$?" the answer is a straightforward ``no," as the following example shows.
	\begin{example} \label{counterexample}
		Define $m$ as a density
		\begin{equation*}
			m(x) = \frac{1}{x(\ln x)^2}\bb{I}_{(0,e^{-1})}(x).
		\end{equation*}
		Note that $m$ is a probability density on $\s{D}$ with cdf
		\begin{equation*}
			F(x) = \int_0^x m(s)\dif s = -\frac{1}{\ln x}\bb{I}_{(0,e^{-1})}(x) + \bb{I}_{\intco{e^{-1},\infty}}(x).
		\end{equation*}
		Assume that $\eta^h[m](\cdot)$ is $\alpha$-H\"older continuous on $[0,T]$ for some $\alpha \in (0,1)$.
		Then there exists a constant $C$ such that
		\begin{equation*}
			1 - \eta^h[m](s) = \frac{2}{\sqrt{2\sigma^2 \pi}}\int_0^\infty  F(\sqrt{s}x)e^{-\frac{x^2}{2\sigma^2}}\dif x \leq Cs^\alpha \ \forall s > 0,
		\end{equation*}
		and so, by Fatou's Lemma,
		\begin{equation*}
			\frac{2}{\sqrt{2\sigma^2 \pi}}\int_0^\infty \liminf_{s \to 0^+}s^{-\alpha} F(\sqrt{s}x)^{-\frac{x^2}{2\sigma^2}}\dif x \leq C.
		\end{equation*}
		But for any $x > 0$, we have
		\begin{equation*}
			\lim_{s \to 0+} s^{-\alpha}F(\sqrt{s}x) = \lim_{s \to 0+} \frac{-1}{s^{\alpha}\ln (\sqrt{s}x)} = +\infty.
		\end{equation*}
		This is a contradiction.
	\end{example}

	For $0 < \alpha < 1$ we define $\sr{M}_\alpha(\s{D})$ to be the space of all $m \in \s{M}(\s{D})$ on $\s{D}$ such that $\eta^h[m] \in \s{C}^\alpha\del{\intco{0,\infty}}$, with norm
	\begin{equation}
	\enVert{m}_{\sr{M}_\alpha} = \enVert{\eta^h[m]}_{\s{C}^\alpha\del{\intco{0,\infty}}}
	+ \enVert{m}_{TV}.
	\end{equation}
	It is straightforward to see that $\sr{M}_\alpha$ is a Banach space.
	The heat equation \eqref{eq:heat} generates a semigroup of contractions on $\sr{M}_\alpha$.
	Indeed, let $m(t)$ denote the (measure-valued) solution at time $t$.
	First we deduce $\enVert{m(t)}_{TV} \leq \enVert{m_0}_{TV}$ by integrating \eqref{eq:heat solution}.
	Moreover, by the semigroup property (i.e.~by uniqueness of solutions to the heat equation) we have $\eta^h[m(t)](s) = \eta^h[m_0](t+s)$, so that
	\begin{equation*}
		\enVert{\eta^h[m(t)]}_{\s{C}^\alpha\del{\intco{0,\infty}}} 
	= \enVert{\eta^h[m_0](t+\cdot)}_{\s{C}^\alpha\del{\intco{0,\infty}}}
	\leq \enVert{\eta^h[m_0]}_{\s{C}^\alpha\del{\intco{0,\infty}}} \quad \forall t \geq 0.
	\end{equation*}

	Example \ref{counterexample} shows that measures which have a steep concentration of mass near 0 will fail to be in $\sr{M}_\alpha$.
	We now show prove that the converse is true, i.e.~an estimate on the concentration of mass near zero will guarantee inclusion in $\sr{M}_\alpha$.
	For any $\alpha > 0$, denote by $\s{M}^\alpha$ the set of all $m \in \s{M}$ satisfying
	\begin{equation} \label{eq:inverse moment}
		\int_{\s{D}} \abs{x}^{-\alpha}\dif\,\abs{m}(x) < \infty.
	\end{equation}
	For instance, $\s{M}^\alpha$ contains all finite measures with support in $\intco{z,\infty}$ for some $z > 0$.
	In particular, $\s{M}^\alpha$ is dense in $\s{M}$.	
	If we endow $\s{M}^\alpha$ with the norm
	\begin{equation}
		\enVert{m}_{\s{M}^\alpha} = \enVert{m}_{TV} + \int_{\s{D}} \abs{x}^{-\alpha}\dif\,\abs{m}(x) = \int_{\s{D}} \del{1 + \abs{x}^{-\alpha}}\dif\,\abs{m}(x),
	\end{equation}
	then it is straightforward to see that $\s{M}^\alpha$ is a Banach space.
	We will also denote $\s{M}^\alpha_+ = \s{M}^\alpha \cap \s{M}_+$, i.e.~the set of all \emph{positive} measures such that \eqref{eq:inverse moment} holds.
	\begin{proposition} \label{pr:inverse moment}
		Let $\alpha \in (0,2)$.
		Then $\s{M}^\alpha \subset \sr{M}_{\alpha/2}$, and there exists a constant $C(\alpha)$ such that
		\begin{equation}
			\enVert{m}_{\sr{M}_{\alpha/2}} \leq C(\alpha)\enVert{m}_{\s{M}^{\alpha}} \quad \forall m \in \s{M}^{\alpha}.
		\end{equation}
		In particular, $\sr{M}_{\alpha/2}$ is dense in $\s{M}$.
	\end{proposition}
	
	\begin{proof}
		We can write
		\begin{equation}
			\eta^h[m](t) = \int_{\s{D}} f(y,t)m_0(\dif y),
		\end{equation}
		where
		\begin{equation}
			f(y,t) = \int_0^\infty \del{S(x-y,t) - S(x+y,t)}\dif x.
		\end{equation}
		We observe that
		\begin{equation}
		\begin{split}
			\dpd{f}{t}(y,t) &= \int_0^\infty \del{\dpd{S}{t}(x-y,t) - \dpd{S}{t}(x+y,t)}\dif x\\
		&= \frac{\sigma^2}{2}\int_0^\infty \del{\dpd[2]{S}{x}(x-y,t) - \dpd[2]{S}{x}(x+y,t)}\dif x\\
		&= -\sigma^2\dpd{S}{x}(y,t)
		 = \frac{y}{\sqrt{2\sigma^2\pi}t^{3/2}}e^{-\frac{y^2}{2\sigma^2 t}}.
		\end{split}
		\end{equation}
		Let $p \geq 1$, $y > 1$.
		By a change of variables $s = y^2/t$, we deduce 
		\begin{equation*}
		\del{\int_0^\infty \abs{f_t(y,s)}^p \dif s}^{1/p} = C(p)y^{-2/p'}, \ p' := p/(p-1).
		\end{equation*}
		Therefore
		\begin{equation*}
		\abs{f(y,t_1)-f(y,t_2)} \leq C(p)y^{-2/p'}\abs{t_1-t_2}^{1/p'}
		\end{equation*}
		We choose $p = 2/(2-\alpha)$, or equivalently $p' = 2/\alpha$.
		Then we have
		\begin{equation*}
		\abs{\eta^h[m](t_1)-\eta^h[m](t_2)} \leq \int_0^\infty \abs{f(y,t_1)-f(y,t_2)}m(y)\dif y
		\leq C(\alpha)\abs{t_1-t_2}^{\alpha/2}\int_0^\infty y^{-\alpha}m(y)\dif y.
		\end{equation*}
		The claim follows.
	\end{proof}

	Recall that the heat semigroup is a semigroup of contractions on $\sr{M}_\alpha$.
	It turns out that the heat semigroup is also bounded on $\s{M}^\alpha$ for arbitrary $\alpha > 0$, as the following lemma implies.
	\begin{lemma}
		\label{lem:Malpha bound}
		Let $m_0$ be a positive measure satisfying \eqref{eq:inverse moment} for some $\alpha > 0$.
		There exists a constant $C(\alpha)$ such that if $m$ is the solution of the heat equation \eqref{eq:heat}, then
		\begin{equation} \label{eq:Malpha bound}
			\int_{\s{D}} \abs{x}^{-\alpha}m(\dif x,t) \leq C(\alpha)\int_{\s{D}} \abs{x}^{-\alpha}m_0(\dif x).
		\end{equation}
	\end{lemma}

	The proof of Lemma \ref{lem:Malpha bound}, which is found in Appendix \ref{ap:fokkerplanck proofs}, relies on the following result, which will be useful for other estimates on parabolic equations.
	\begin{lemma}\label{lem:m_n}
		Let $S(x,t)$ be the heat kernel, defined in \eqref{eq:heat kernel}.
		For all $n = 0,1,2,\ldots$, there exists a (Hermite) polynomial $P_n$ of degree $n$ such that
		\begin{equation} \label{eq:dnSdx}
			\dpd[n]{S}{x}(x,t) = \del{\sigma^2 t}^{-n/2}P_n\del{\frac{\abs{x}}{\sqrt{\sigma^2 t}}}S(x,t).
		\end{equation}
		As a corollary, for all $n = 0,1,2,\ldots,$ and $k = 1,2,3,\ldots$ the constants
		\begin{equation} \label{eq:m_n}
			m_n := \sup_{x,t} \abs{x}^{n+1}\abs{\dpd[n]{S}{x}(x,t)},
			\quad
			m_{n,k} := \sup_{x,t} \abs{x}^{n+1-k}\del{\sigma^2 t}^{k/2}\abs{\dpd[n]{S}{x}(x,t)}
		\end{equation}
		are finite and depend only on $n$ and $k$.
	\end{lemma}
	\begin{proof}
		The proof of \eqref{eq:dnSdx} is elementary using induction.
		The second claim follows from the fact that $\sup_{x \geq 0} x^\alpha e^{-x}$ is finite for any $\alpha \geq 0$.
	\end{proof}

	We conclude this section by generalizing our results to the Fokker-Planck equation for an arbitrary bounded continuous drift term $b(x,t)$.
	The proofs are found in Appendix \ref{ap:fokkerplanck proofs}.

	\begin{lemma} \label{lem:fp eta holder}
		Let $b$ be a bounded continuous function on $\s{D} \times [0,T]$, let $m_0 \in \s{M}_{+}(\s{D}) \cap \sr{M}_{\alpha}(\s{D})$, and let $m$ be the unique weak solution $m$ of \eqref{eq:fp}, given by Lemma \ref{lem:fp}.
		Then the total mass function $\eta(t) := \int_{\s{D}} m(\dif x,t)$ is $\beta$-H\"older continuous for $\beta = \min\{\alpha,1/2\}$, with
		\begin{equation}
		\label{eq:eta holder}
		\enVert{\eta}_{C^{\beta}([0,T])} \leq C(\sigma)\del{\enVert{m_0}_{\sr{M}_\alpha} + \enVert{b}_\infty}.
		\end{equation}
	\end{lemma}

	\begin{lemma}
		\label{lem:Malpha fp bound}
		Let $b$ be a bounded continuous function on $\s{D} \times [0,T]$, let $m_0 \in \s{M}_{+}^\alpha(\s{D})$ for some $\alpha > 0$, and let $m$ be the unique weak solution $m$ of \eqref{eq:fp}, given by Lemma \ref{lem:fp}.
		Then there exists some constants $C(\alpha)$ and $C(\alpha,\sigma)$ such that
		\begin{equation} \label{eq:inverse moment bound fp}
			\int_{\s{D}} \abs{x}^{-\alpha}m(\dif x,t) \leq C(\alpha)e^{C(\alpha,\sigma)\enVert{b}_\infty t}\int_{\s{D}} \abs{x}^{-\alpha}m_0(\dif x).
		\end{equation}
	\end{lemma}

	\section{Forward-backward system}
	\label{sec:fwdbckwd}
	
	In this section we prove existence and uniqueness of solutions to \emph{infinite time horizon forward-backward system} \eqref{eq:mfg infty}.
	Many of the ideas in this section can already be found in \cite{graber2021nonlocal}.
	Our result is novel in that (i) the time horizon is infinite and (ii) the initial measure $m_0$ need not be smooth nor even a density.
	The proof is based on a priori estimates followed by an application of the Leray-Schauder fixed point theorem (see e.g.~\cite[Theorem 11.3]{gilbarg2015elliptic}).
	Most of the proofs in this section involve either standard computations or ideas that can be found in the previous works \cite{graber2021nonlocal,graber2020commodities,graber2018existence},
	and so we relegate them to Appendix \ref{ap:forward-backward proofs}.
	However, in the sequel we will make frequent reference to the \emph{estimates} found in this section.

	\subsection{The Hamiltonian}
	\label{sec:hamiltonian}
	
	In this subsection we deduce a number of structural features of the Hamiltonian, using only Assumptions \ref{as:P} and \ref{as:prudence}.
	The proofs can be found in Appendix \ref{ap:hamiltonian proofs}.
	
	\begin{lemma} \label{lem:opt quant}
		[Unique optimal quantity]
		The function $q^*:\intco{0,\infty}^3 \to \intco{0,\infty}$ given by $q^*(\epsilon,Q,a) = \argmax_{q \geq 0} \pi(\epsilon,q,Q,a)$ is well-defined and locally Lipschitz continuous.
		It is non-increasing in the variable $a$.
		With respect to $\epsilon$ and $Q$, it satisfies
		\begin{equation} \label{eq:qstarQ leq}
		-\epsilon \leq \dpd{q^*}{Q} \leq \epsilon\frac{\bar \rho-1}{2-\bar \rho}, \ -Q \leq \dpd{q^*}{\epsilon} \leq Q\frac{\bar \rho-1}{2-\bar \rho}.
		\end{equation}
		
		Define $H(\epsilon,Q,a) = \pi\del{\epsilon,q^*(Q,a),Q,a} \geq 0$.
		Then $H$ is locally Lipschitz, decreasing in all variables, and convex in $a$; its derivative $\dpd{H}{a} = -q^*$ is also locally Lipschitz.
	\end{lemma}

	\begin{corollary}
		[Smoothness and uniform convexity]
		\label{cor:smoothness}
		Let $\bar\epsilon \geq 0,\bar Q \geq 0,$ and $\bar a > 0$ be constants such that $\bar a < P(\bar \epsilon \bar Q)$.
		Consider the restriction of $H = H(\epsilon,Q,a)$ to the domain $[0,\bar \epsilon] \times [0,\bar Q] \times [0,\bar a]$.
		Then $H$ is $n$ times continuously differentiable with Lipschitz continuous derivatives, where $n$ is the same as in Assumption \ref{as:P}.
		It is also uniformly convex in the $a$ variable, and in particular there exists a constant $C_H = C(\bar\epsilon,\bar Q,\underline{a},\bar a) \geq 1$ such that
		\begin{equation} \label{eq:C_H}
		C_H^{-1} \leq \dpd[2]{H}{a}(\epsilon,Q,a) \leq C_H \quad \forall (\epsilon,Q,a) \in [0,\bar \epsilon] \times [0,\bar Q] \times [\underline{a},\bar a].
		\end{equation}
	\end{corollary}

	\begin{corollary}
		[$Q$ dependence] \label{cor:dHdQ}
		We have the following estimates in the region where $P(\epsilon Q) > a$:
		\begin{equation} \label{eq:dHdQ estimate}
		\abs{\dpd{H}{Q}} \leq \epsilon (P(0)-a),
		\quad
		\abs{\dmpd{H}{Q}{a}} \leq \epsilon\max\cbr{\abs{\frac{\bar \rho - 1}{\bar \rho - 2}},1} =: \bar P \epsilon.
		\end{equation}
	\end{corollary}

	\begin{lemma}[Unique aggregate quantity] \label{lem:market clearing Cournot}
		Let $\epsilon \geq 0, \phi \in L^\infty(\s{D})$ and $m \in \s{M}_+(\s{D})$ with $\int_{\s{D}} \dif m(x) \leq 1$ and $\phi \geq 0$ (a.e.).
		Then there exists a unique $Q^* = Q^*(\epsilon,\phi,m) \geq 0$ such that
		\begin{equation}\label{eq:market clearing Cournot}
		Q^* = \int_{\s{D}} q^*(\epsilon,Q^*,\phi(x))\dif m(x) = -\int_{\s{D}} \dpd{H}{a}\del{\epsilon,Q^*,\phi(x)}\dif m(x).
		\end{equation}
		Moreover, $Q^*$ satisfies the a priori estimate
		\begin{equation} \label{eq:Q upper bound}
		Q^* \leq c(\bar \rho,\epsilon)q^*(0,0,0), \ c(\bar \rho,\epsilon) := \max\cbr{\frac{2-\bar \rho}{2+\epsilon - (1+\epsilon)\bar \rho},1}.
		\end{equation}
		Finally, $Q^*$ is locally Lipschitz in the following sense.
		If $\epsilon_1,\epsilon_2 \in [0,\epsilon]$, $\phi_1,\phi_2$ Lipschitz functions with $\enVert{\phi_i}_\infty \leq M$, and $m_1,m_2 \in \s{M}_{1,+}(\s{D})$ with $\int_{\s{D}} \dif m_i(x) \leq 1$, set $Q^*_i = Q^*(\epsilon_i,\phi_i,m_i)$ to be the solution of \eqref{eq:market clearing Cournot} corresponding to $\epsilon_i,\phi_i,m_i$ for $i = 1,2$. Then there exists a constant $C = C(\epsilon,\bar \rho,M)$ such that
		\begin{equation}
		\label{eq:Q1-Q2}
		\begin{split}
		\abs{Q_1^*-Q_2^*} &\leq C\del{\abs{\epsilon_1-\epsilon_2} + \int_{\s{D}} \abs{\phi_1(x)-\phi_2(x)}\dif m_1(x) + \max_{i=1,2}\enVert{\dod{\phi_i}{x}}_\infty{\bf d}_1(m_1,m_2) + \abs{\int_{\s{D}} \dif (m_1-m_2)(x)}}\\
		&\leq C\del{\abs{\epsilon_1-\epsilon_2} + \enVert{\phi_1-\phi_2}_\infty + \max_{i=1,2}\enVert{\dod{\phi_i}{x}}_\infty{\bf d}_1(m_1,m_2) + \abs{\int_{\s{D}} \dif (m_1-m_2)(x)}}.
		\end{split}
		\end{equation}
	\end{lemma}
	\begin{remark}
		\label{rmk:c(rho,ep) increasing}
		The function $c(\bar \rho,\epsilon)$ in equation \eqref{eq:Q upper bound} is an increasing function of $\epsilon$.
	\end{remark}

\begin{corollary} \label{cor:q* a priori bound}
	Let $\epsilon,\phi,m,$ and $Q^* = Q^*(\epsilon,\phi,m)$ be as in Lemma \ref{lem:market clearing Cournot}.
	Then
	\begin{equation}\label{eq:q* a priori bound}
	q^*\del{\epsilon,Q^*,\phi(x)} \leq c(\bar \rho,\epsilon)q^*\del{0,0,0}
	\end{equation}
	for a.e.~$x \in \s{D}$.
\end{corollary}
	
	\subsection{Finite time horizon problem}

	In this section we fix a final time $T > 0$ and consider the forward-backward system only on this time horizon.
	For technical reasons, we will need to replace the constant $\epsilon$ with a function $\epsilon(t)$ such that $\epsilon(T) = 0$.
	System \eqref{eq:mfg infty} becomes
	\begin{equation} \label{eq:mfg T}
	\begin{cases}
	(i) & \dpd{u}{t} + \dfrac{\sigma^2}{2}\dpd[2]{u}{x} + H\del{\epsilon(t),Q^*(t),\dpd{u}{x}} - ru = 0,\\
	(ii) & \dpd{m}{t} - \dfrac{\sigma^2}{2}\dpd[2]{m}{x} + \dpd{}{x}\del{\dpd{H}{a}\del{\epsilon(t),Q^*(t),\dpd{u}{x}}m} = 0,\\
	(iii) & Q^*(t) = -\displaystyle \int_{\s{D}} \dpd{H}{a}\del{\epsilon(t),Q^*(t),\dpd{u}{x}}\dif m(t),\\
	(iv) & m|_{x=0} = u|_{x=0} = 0, \
	m|_{t=0} = m_0 \in \s{P}(\s{D}), \ u|_{t=T} = u_T \in C^{2+\alpha}.
	\end{cases}
	\end{equation}
	We define $(u,m)$ to be a solution to \eqref{eq:mfg T} provided that $u$ is a smooth function on $\overline{\s{D}} \times [0,T]$ (twice continuously differentiable with respect to $x$, continuously differentiable with respect to $t$), $m \in \s{C}([0,T];\s{P}(\s{D}))$, Equations (i) and (iii) are satisfied pointwise,  the boundary conditions for $u$ in (iv) are satisfied pointwise, and Equation (ii) with the boundary conditions for $m$ from (iv) holds in the sense of distributions (see Section \ref{sec:fokkerplanck}).
	Note that a solution $(u,m)$ must satisfy $\dpd{u}{x} \geq 0$, because the domain of $H$ is $\intco{0,\infty}^3$.
	It is possible to relax this somewhat by extending the domain of $H(\epsilon,Q,a)$ to include all $a > \lim_{q \to \infty} P(q)$, but we need not do so here.
	
	\begin{assumption}
		[Structure of $\epsilon(t)$] \label{as:epsilon(t)}
		We assume $\epsilon$ is a smooth, non-negative, non-increasing function on $\intcc{0,T}$ such that $\epsilon(T) = 0$ and $\enVert{\epsilon'}_\infty \leq 1$.
	\end{assumption}
	
	\begin{assumption}[Structure of $u_T$] \label{as:uT}
		For each $T > 0$, the function $u_T$ is an element of $\s{C}^{2+\alpha}\del{\overline{\s{D}}}$ that satisfies the following conditions:
		\begin{enumerate}
			\item $u_T(0) = 0$;
			\item $\frac{\sigma^2}{2}u_T''(0) + H\del{0,0,u_T'(0)} = 0$;
			\item $0 \leq u_T(x) \leq c_1$ for all $x \in \s{D}$, where $c_1 > 0$ is some constant;
			\item there exists a constant $c_3 > 0$, independent of $T$, such that $0 \leq u_T'(x) \leq c_3$ for all $x \in \s{D}$ and all $T > 0$;
			\item there exists a constant $C_\alpha$, independent of $T$, such that $\enVert{u_T}_{\s{C}^{2+\alpha}\del{\overline{\s{D}}}} \leq C_\alpha$ for all $T > 0$.
		\end{enumerate}
	\end{assumption}
	\begin{remark} \label{rmk:c c'}
		It is always possible to satisfy Assumption \ref{as:uT} for an arbitrary constant $c_3 > 0$.
		Here we give one possible construction.
		Set $h = \frac{2}{\sigma^2}H(0,0,c_3)$, so that condition (2) becomes $u_T''(0) = -h$.
		If $h > 0$,
		then Assumption \ref{as:uT} is satisfied by the function
		\begin{equation}
			u_T(x) = \frac{2(c_3)^2}{3h} + \frac{h^2}{12c_3}\del{x - \frac{2c_3}{h}}_-^3,
		\end{equation}
		where $x_- := \min\cbr{x,0}$.
		In the case where $h = 0$, Assumption \ref{as:uT} is satisfied by the function
		\begin{equation}
			u_T(x) = \begin{cases}
				c_3 x - \frac{x^3}{6} &\text{if}~x \leq (c_3)^{1/2},\\
				\frac{1}{2}(c_3)^{3/2} + \frac{1}{3}\del{x - 2(c_3)^{1/2}}^3, &\text{if}~(c_3)^{1/2} \leq x \leq 2(c_3)^{1/2},\\
				\frac{1}{2}(c_3)^{3/2}, &\text{if}~2(c_3)^{1/2} \leq x.
			\end{cases}
		\end{equation}
		Note also that these examples can be slightly modified to produce globally $\s{C}^\infty$ functions satisfying Assumption \ref{as:uT}.
	\end{remark}

	\subsection{Estimates on the Hamilton-Jacobi equation}
	\label{sec:estimates on hj}
	
	\begin{lemma}
		[A priori estimates for HJ equation]
		\label{lem: a priori HJ}
		Let $Q^*(t)$ be any bounded, non-negative function.
		Let $u$ be a solution of the Hamilton-Jacobi equation
		\begin{equation} \label{eq:HJ}
		\dpd{u}{t} + \dfrac{\sigma^2}{2}\dpd[2]{u}{x} + H\del{\epsilon(t),Q^*(t),\dpd{u}{x}} - ru = 0, \ x \in \s{D}, \ t \in \intco{0,T}
		\end{equation}
		with Dirichlet boundary conditions $u(0,t) = 0$ and final condition $u(x,T) = u_T(x)$, which satisfies Assumption \ref{as:uT}.
		Then for all $x \in \s{D}$ and $t \in [0,T]$, we have
		\begin{equation} \label{eq:u ux bounded}
		0\leq u(x,t) \leq \frac{1}{r}H(0,0,0) + c_1, \quad   0 \leq u_x(x,t) \leq M,
		\end{equation}
		where 
		\begin{equation} \label{eq:max ux}
		M = M(\sigma,r,c_1,c_3) := \begin{cases}
			2\sqrt{\frac{2H(0,0,0)\del{H(0,0,0) + rc_1}}{\sigma^2 r}} &\text{if}~ c_3 \leq \sqrt{\frac{2}{\sigma^2 r}}H(0,0,0)\\
			c_3 + \frac{2}{\sigma^2r c_3}H(0,0,0)^2
			+ \frac{2c_1}{\sigma^2 c_3}H(0,0,0) &\text{if}~c_3 \geq \sqrt{\frac{2}{\sigma^2 r}}H(0,0,0)
		\end{cases}.
		\end{equation}
	\end{lemma}

	\begin{proof}
		See Appendix \ref{ap:hj proofs}.
		Cf.~\cite[Section 4]{graber2021nonlocal}.
	\end{proof}

\subsection{Estimates on the coupling}

\label{sec:estimates on coupling}

\begin{lemma} \label{lem:Q(t) Holder}
	Let $(u,m)$ be a solution of \eqref{eq:mfg T}.
	Then $Q^*$, given by \eqref{eq:mfg T}(iii), satisfies the following bounds:
	\begin{equation}
	\label{eq:Q(t) upper bound}
	0 \leq Q^*(t) \leq c(\bar \rho,\epsilon(0))q^*(0,0,0) = -c(\bar \rho,\epsilon(0))\dpd{H}{a}\del{0,0,0},
	\end{equation}
	where $c(\bar \rho,\epsilon)$ is defined in \eqref{eq:Q upper bound}.
	
	Suppose, moreover, that $m_0 \in \sr{M}_{\alpha/2}$ for some $\alpha \in \intoc{0,1}$.
	Then $Q^*(t)$ is H\"older continuous on $\intcc{0,T}$ with
	\begin{equation}\label{eq:Q(t) Holder} 
		\enVert{Q^*}_{\s{C}^{\alpha/2}} \leq C\del{\enVert{\dpd{u}{x}}_{\s{C}^{\alpha,\alpha/2}} + \enVert{\dpd[2]{u}{x}}_{\infty} + 1}.
	\end{equation}
	for some $C = C\del{\bar \rho,\epsilon(0),\sigma,M,\enVert{m_0}_{\sr{M}_{\alpha/2}}}$, where $M$ is the constant from Lemma \ref{lem: a priori HJ}  that gives an upper bound on $\enVert{\pd{u}{x}}_\infty$.
	\begin{comment}
	As a corollary, for any $\delta > 0$ arbitrary, there exists\\ $C_\delta = C\del{\delta,\bar \rho,\epsilon(0),\sigma,M,\enVert{m_0}_{\sr{M}_{\alpha/2}},\alpha}$ such that
	\begin{equation}\label{eq:Q(t) Holder delta} 
	\enVert{Q^*}_{\s{C}^{\alpha/2}} \leq \delta \enVert{u}_{\s{C}^{2+\alpha,1+\alpha/2}} + C_\delta.
	\end{equation}
	\end{comment}
\end{lemma}

\begin{proof}
	See Appendix \ref{ap:coupling proofs}.
\end{proof}

\subsection{Parabolic estimates}

\label{sec:parabolic estimates}

Before stating our result on the existence of smooth solutions to the system, we present some estimates on solutions to parabolic problems that \emph{do not depend on the time horizon}.
These estimates will be useful in study of the linearized system (Section \ref{sec:sensitivity}).

\begin{lemma} \label{lem:C2alpha estimates}
	Let $T > 0, r > 0$ be given.
	For any $f \in \s{C}^{\alpha,\alpha/2}(\overline{\s{D}} \times \intcc{0,T})$, and $u_0 \in \s{C}^{2+\alpha}(\overline{\s{D}})$,
	there exists a unique solution $u \in \s{C}^{2+\alpha,1+\alpha/2}(\overline{\s{D}} \times \intcc{0,T})$ of
	\begin{equation} \label{eq:parabolic + ru}
		\dpd{u}{t} + ru - \frac{\sigma^2}{2}\dpd[2]{u}{x} = f, \ \forall x \in \s{D}, t > 0; u(0,t) = 0 \ \forall t > 0; \ u(x,0) = u_0(x) \ \forall x \in \s{D}
	\end{equation}
	satisfying
	\begin{equation} \label{eq:global holder estimate}
		\enVert{u}_{\s{C}^{2+\alpha,1+\alpha/2}(\overline{\s{D}} \times \intcc{0,T})} \leq C(\sigma,r,\alpha)\del{\enVert{f}_{\s{C}^{\alpha,\alpha/2}(\overline{\s{D}} \times \intcc{0,T})} + \enVert{u_0}_{\s{C}^{2+\alpha}(\overline{\s{D}})}}.
	\end{equation}
	The constant $C(\sigma,r,\alpha)$ in \eqref{eq:global holder estimate} does \emph{not} depend on $T$.
	More specifically, we can say that if $r \geq 1$,
	\begin{equation} \label{eq:global holder estimate1}
		\enVert{u}_{\s{C}^{2+\alpha,1+\alpha/2}(\overline{\s{D}} \times \intcc{0,T})} \leq C(\sigma,\alpha)\del{\sbr{f}_{\alpha,\alpha/2} + r^{\frac{\alpha}{2}}\enVert{f}_0 + \sbr{u_0}_{2+\alpha}
			+ r^{1 + \frac{\alpha}{2}}\enVert{u_0}_0}.
	\end{equation}
\end{lemma}

\begin{proof}
	The result follows from potential estimates found in \cite[Chapter IV]{ladyzhenskaia1968linear}.
	See Appendix \ref{ap:parabolic proofs}.
\end{proof}

\subsection{Existence of solutions}

\label{sec:existence}

\begin{lemma} \label{lem:uC2alpha}
	Let $m_0 \in \sr{M}_{\alpha/2}$ and $0 < \alpha \leq 1$.
	Then there exists a constant $$C = C(\bar \rho, \epsilon(0),\sigma,M,c_1,\enVert{m_0}_{\sr{M}_{\alpha/2}},\alpha)$$ such that for any solution $(u,m)$ of \eqref{eq:mfg T},	
	\begin{equation} \label{eq:uC2alpha}
	\enVert{u}_{\s{C}^{2+\alpha,1+\alpha/2}} \leq C\del{1 + r^{\frac{\alpha}{2}} + C_\alpha + r^{1 + \frac{\alpha}{2}}c_1},
	\end{equation}
	where $M$ is the constant from Lemma \ref{lem: a priori HJ} and $c_1,c_3,C_\alpha$ are the constants from Assumption \ref{as:uT}.
\end{lemma}

\begin{remark}
	The constant on the right-hand side of \eqref{eq:uC2alpha} does not depend on $T$.
\end{remark}

For the proof of Lemma \ref{lem:uC2alpha}, see Appendix \ref{ap:existence proofs}.

	\begin{theorem}
		[Existence of classical solutions for \eqref{eq:mfg T}] \label{thm:exist mfg T}
		Let $m_0 \in \sr{M}_{\alpha/2}$ and $0 < \alpha \leq 1$.
		Then there exists a solution $(u,m)$ satisfying the finite time horizon problem \eqref{eq:mfg T} and having the following regularity: $u \in \s{C}^{2+\alpha,1+\alpha/2}\del{\overline{\s{D}} \times [0,T]}$, $m \in \s{C}^{1/2}\del{[0,T];\s{M}_{1,+}(\s{D})}$.
		Thus, Equation \eqref{eq:mfg T}(i) is satisfied in a classical sense, while Equation \eqref{eq:mfg T}(ii) is satisfied in the weak sense defined in \eqref{eq:fp weak}, and Equation \eqref{eq:mfg T}(iii) holds pointwise.
	\end{theorem}

	\begin{proof}
		We use the Leray-Schauder fixed point theorem in a more or less standard way, cf.~\cite{graber2021nonlocal,graber2020commodities,graber2018existence}.
		The details are given in Appendix \ref{ap:existence proofs}.
	\end{proof}

\begin{theorem}[Existence of solutions to the infinite horizon problem \eqref{eq:mfg infty}]
	\label{thm:exist mfg infty}
	Let $m_0 \in \sr{M}_{\alpha/2}$ and $0 < \alpha \leq 1$.
	Then there exists a solution $(u,m) \in \s{C}^{2+\alpha,1+\alpha/2}\del{\overline{\s{D}} \times \intco{0,\infty}} \times \in \s{C}^{1/2}\del{\intco{0,\infty};\s{M}_{1,+}(\s{D})}$ solving the infinite time horizon problem \eqref{eq:mfg infty} and satisfying the following estimates:
	\begin{equation} \label{eq:main estimates}
	\begin{split}
	&\enVert{u}_{\s{C}^{2+\alpha,1+\alpha/2}} 
	\leq C(\bar \rho, \epsilon,\sigma,M,\enVert{m_0}_{\sr{M}_{\alpha/2}},\alpha)\del{1 + r^{\frac{\alpha}{2}} + C_\alpha},\\
	&{\bf d}_1\del{m(t_1),m(t_2)} \leq 2\del{M + \sigma}\abs{t_1-t_2}^{1/2} \ \forall \abs{t_1-t_2} \leq 1,\\
	&0\leq u(x,t) \leq \frac{1}{r}H(0,0,0), \quad   0 \leq \dpd{u(x,t)}{x} \leq M \ \forall x \in \overline{\s{D}}, t \geq 0,\\
	&0 \leq Q^*(t) \leq \bar Q,\\
	&0 \leq -\dpd{H}{a}\del{\epsilon,Q^*(t),\dpd{u}{x}(x,t)} \leq \bar Q \quad \forall (x,t) \in \s{D} \times \intco{0,\infty}
	\end{split}
	\end{equation}
	where $M$ and $\bar Q$ are defined by
	\begin{equation} \label{eq:M}
		M := 2\sqrt{\frac{2}{\sigma^2 r}}H(0,0,0), \quad \bar Q := -c\del{\bar \rho,\epsilon}\dpd{H}{a}(0,0,0)
	\end{equation}
\end{theorem}

\begin{proof}
	For each $T > 0$, we will let $\epsilon(t)$ be a function satisfying Assumption \ref{as:epsilon(t)} as well as $\epsilon(0) = \epsilon$, and we let $u_T$ be a function satisfying Assumption \ref{as:uT}.
	By Theorem \ref{thm:exist mfg T} there exists a solution of \eqref{eq:mfg T}, which we denote $(u^T,m^T)$.
	Fix an arbitrary $T_0 > 0$.
	By Lemmas \ref{lem:uC2alpha} and \ref{lem:fp}, $(u^T,m^T)$ is uniformly bounded in $\s{C}^{2+\alpha,1+\alpha/2}\del{\overline{\s{D}} \times \intcc{0,T_0}} \times \s{C}^{1/2}\del{\intcc{0,T_0};\s{M}_{1,+}(\s{D})}$ for all $T \geq T_0$, with norms bounded by a constant that does not depend on $T_0$.
	Thus, by standard diagonalization, we may pass to a subsequence, still denoted $(u^T,m^T)$, that converges  to some fixed $(u,m)$, where the convergence is in $\s{C}^{2,1}\del{\overline{\s{D}} \times \intcc{0,T_0}} \times \s{C}^{0}\del{\intcc{0,T_0};\s{M}_{1,+}(\s{D})}$ for \emph{every} $T_0$.
	By the uniform estimates on $(u^T,m^T)$ it also follows that $(u,m) \in \s{C}^{2+\alpha,1+\alpha/2}\del{\overline{\s{D}} \times \intco{0,\infty}} \times \s{C}^{1/2}\del{\intco{0,\infty};\s{M}_{1,+}(\s{D})}$.
	To see that $(u,m)$ is indeed a solution to \eqref{eq:mfg infty}, it suffices to pass to the limit in the equations satisfied by $(u^T,m^T)$ on arbitrary time horizons.
	Finally, note that the following estimates hold:
	\begin{equation} \label{eq:main estimates1}
		\begin{split}
			&\enVert{u}_{\s{C}^{2+\alpha,1+\alpha/2}} 
			\leq C(\bar \rho, \epsilon(0),\sigma,M,c_1,c_3,\enVert{m_0}_{\sr{M}_{\alpha/2}},\alpha)\del{1 + r^{\frac{\alpha}{2}} + C_\alpha + r^{1 + \frac{\alpha}{2}}c_1},\\
			&{\bf d}_1\del{m(t_1),m(t_2)} \leq 2\del{M + \sigma}\abs{t_1-t_2}^{1/2} \ \forall \abs{t_1-t_2} \leq 1,\\
			&0\leq u(x,t) \leq \frac{1}{r}H(0,0,0) + c_1, \quad   0 \leq \dpd{u(x,t)}{x} \leq M(\sigma,r,c_1,c_3) \ \forall x \in \overline{\s{D}}, t \geq 0,\\
			&0 \leq Q^*(t) \leq -c\del{\bar \rho,\epsilon}\dpd{H}{a}(0,0,0),\\
			&0 \leq -\dpd{H}{a}\del{\epsilon,Q^*(t),\dpd{u}{x}(x,t)} \leq -c\del{\bar \rho,\epsilon}\dpd{H}{a}(0,0,0) \quad \forall (x,t) \in \s{D} \times \intco{0,\infty}
		\end{split}
	\end{equation}
	where $M = M(\sigma,r,c_1,c_3)$ is defined in \eqref{eq:max ux}.
	This follows because they hold  for $(u^T,m^T)$ uniformly in $T$ (Lemmas \ref{lem:fp}, \ref{lem: a priori HJ}, and \ref{lem:Q(t) Holder}, also Corollary \ref{cor:q* a priori bound}).
	Now by Remark \ref{rmk:c c'}, $c_1,c_2$ and $c_3$ can be made arbitrarily close to zero.
	Letting $c_1,c_3 \to 0$ and using the continuity of $H$ and $\pd{H}{a}$, we deduce the estimates \eqref{eq:main estimates}.
\end{proof}

\subsection{Uniqueness and smoothness of the Hamiltonian}
\label{sec:uniqueness}
When the demand schedule is linear, uniqueness of solutions to \eqref{eq:mfg infty} follows with no further conditions on the data, cf.~\cite{graber2020commodities}.
In the case of a general, nonlinear demand schedule satisfying Assumptions \ref{as:P} and \ref{as:prudence}, we can prove uniqueness of solutions for small enough parameter $\epsilon$.
Cf.~\cite{graber2021nonlocal}.
The smallness of $\epsilon$ makes two contributions.
First, it ensures that the Hamiltonian $H$ is a smooth, uniformly convex function on the domain where solutions exist.
Second, it ensures that certain ``energy estimates" \`a la Lasry-Lions (see \cite{lasry07}) hold, which prove uniqueness.
The case where $\epsilon$ is small has independent interest, aside from being a technical condition that yields uniqueness.
(Cf.~Remark \ref{rem:why ep small}.)
\begin{remark}
	\label{rem:why ep small}
	The inspiration for taking $\epsilon > 0$ small is taken from the basic idea that Chan and Sircar use to compute solutions \cite{chan2015bertrand,chan2017fracking}
	Namely, it is natural to try take a formal Taylor expansion of the solution with respect to $\epsilon$ around zero, since at $\epsilon = 0$ the system of equations is decoupled.
	(See \cite{graber2022parameter} for a justification of this technique.)
	Now when $\epsilon > 0$ is small enough, one might think to simplify our approach by devising a contraction mapping argument.
	In the present work, we do not take this approach, but instead seek to unify as much as possible with the case where the demand schedule is linear.
	For in this latter case, it is essentially from the structure of the Hamiltonian that one obtains the ``propagation of monotonicity" (cf.~\cite{gangbo2021mean}) that is needed to prove uniqueness.
	We show that the same is true when $\epsilon$ is small, and we do so by proving the same type of estimates as we do for the linear demand schedule.
	One could, in principle, generalize this idea to other ``smallness" conditions; for example, if the demand schedule is ``close enough to linear" in a suitable sense, then our arguments for uniqueness will go through for a wide range of parameters $\epsilon$.
	In the present work, however, we do not pursue this direction, so as to avoid a multiplication of technicalities.
\end{remark}

In this section we consider both the smoothness of the Hamiltonian and uniqueness of solutions separately.
The former can at first be viewed as a tool for proving the latter, in the case of a nonlinear demand schedule.
However, when we prove the regularity of the master field in Sections \ref{sec:sensitivity} and \ref{sec:solution master}, the smoothness of the Hamiltonian will be required even when the demand schedule is linear.
Therefore we address it in a separate subsection.

\subsubsection{Assumptions ensuring that the Hamiltonian is smooth}

 \label{sec:H is smooth}

The following assumption ensures in general that $H$ can be treated as a smooth, uniformly convex function in System \eqref{eq:mfg infty}.

\begin{assumption}
	\label{as:H must be smooth}
	We assume that $M < P\del{\epsilon \bar Q}$, where $M$ and $\bar Q$ are defined in \eqref{eq:M}.
\end{assumption}
\begin{remark} \label{rmk:H must be smooth}
	[Sufficient conditions to give Assumption \ref{as:H must be smooth}]
	There necessarily exists $r^*$ large enough so that
	\begin{equation*} 
		2\sqrt{\frac{2}{\sigma^2 r}}H(0,0,0) < P(0) \quad \forall r \geq r^*.
	\end{equation*}
	Then, since $\epsilon \mapsto P\del{\epsilon \bar Q}$ is a continuous, decreasing function of $\epsilon$, there exists $\epsilon^* > 0$ such that Assumption \ref{as:H must be smooth} holds for all $0 < \epsilon \leq \epsilon^*$ and all $r \geq r^*$.
\end{remark}

Under Assumption \ref{as:H must be smooth}, it follows from Corollary \ref{cor:smoothness} and the a priori estimates \eqref{eq:main estimates} from Theorem \ref{thm:exist mfg infty} that in System \eqref{eq:mfg infty} (or \eqref{eq:mfg T}, provided $c_2$ from Assumption \ref{as:uT} is chosen small enough), $H$ can be treated as $n$ times continuously differentiable with Lipschitz continuous derivatives, and moreover it is uniformly convex.
In particular, from \eqref{eq:C_H} there exists a constant $C_H  \geq 1$ such that
\begin{equation} \label{eq:C_H smoothness}
	C_H^{-1} \leq \dpd[2]{H}{a}\del{\epsilon,Q(t),\dpd{u}{x}(x,t)} \leq C_H \quad \forall (x,t) \in \s{D} \times (0,\infty)
\end{equation}
whenever $u$ is a solution of \eqref{eq:mfg infty}.

An interesting special case is when the demand schedule is linear; without loss of generality we take $P(q) = 1-q$.
In this case (and in general when $\bar \rho \leq 1$) we have $c(\bar \rho,\epsilon) = 1$, and a simple computation shows $\bar Q = 1/2$ and $M = (2\sigma^2 r)^{-1/2}$.
For any $\epsilon^* < 2$, it is possible to take $r^*$ sufficiently large so that Assumption \ref{as:H must be smooth} holds for any $r \geq r^*$ and any $\epsilon \leq \epsilon^*$.
In this case, the smoothness of $H$ on the domain where solutions lie implies that the solution to \eqref{eq:mfg infty} is the same as the solution to
\begin{equation} \label{eq:mfg lin dem}
	\begin{cases}
		(i) & \dpd{u}{t} + \dfrac{\sigma^2}{2}\dpd[2]{u}{x} + \dfrac{1}{4}\del{1 - \epsilon Q^*(t) - \dpd{u}{x}}^2 - ru = 0,\\
		(ii) & \dpd{m}{t} - \dfrac{\sigma^2}{2}\dpd[2]{m}{x} - \dpd{}{x}\del{\dfrac{1}{2}\del{1 - \epsilon Q^*(t) - \dpd{u}{x}}m} = 0,\\
		(iii) & Q^*(t) = -\displaystyle \int_{\s{D}} \dfrac{1}{2}\del{1 - \epsilon Q^*(t) - \dpd{u}{x}(\cdot,t)}\dif m(t),\\
		(iv) & m|_{x=0} = u|_{x=0} = 0, \
		m|_{t=0} = m_0
	\end{cases}
\end{equation}

\subsubsection{Uniqueness}

\begin{comment}
Finally, we emphasize that our stability result will hold in general for $\epsilon$ small enough.
\begin{assumption}
\label{as:epsilon small}
We have
\begin{equation} \label{eq:epsilon small}
\epsilon \max\cbr{\abs{\frac{\bar \rho - 1}{\bar \rho - 2}},P(0) + 1}  c(\bar \rho,\epsilon)\max\{C_H,\bar Q\}\del{1 + \frac{C_H^2}{\bar Q^2}} \leq \frac{1}{12C_H}.
\end{equation}
We also assume, without loss of generality, that $\epsilon \leq 1$.
\end{assumption}
\begin{remark}
\label{rmk:c(rho,eps) no eps}
As a corollary of our assumption that $\epsilon$ is small, the constant $c(\bar \rho,\epsilon)$, which is increasing in $\epsilon$, can be bounded above by a fixed constant (e.g.~$c(\bar \rho,1)$).
On the other hand, it is bounded below by 1, i.e.~it does not vanish as $\epsilon$ vanishes. For this reason, unless we are attempting to optimize the constants in our estimates, we can treat $c(\bar \rho,\epsilon)$ as if it did not depend on $\epsilon$ at all.
\end{remark}
\end{comment}

\begin{theorem} \label{thm:uniqueness small ep}
	In addition to Assumption \ref{as:H must be smooth}, suppose that 
	\begin{align} \label{eq:r big for uniqueness}
		r &\geq 1000\max\cbr{1 + c(\bar \rho,\epsilon)\bar P \epsilon,1 + c(\bar \rho,\epsilon)\bar Q,\bar Q + \epsilon P(0) + 1}^2 \quad \text{and}\\
		\label{eq:ep small for uniqueness}
		\epsilon &\leq \del{4C_H c(\bar \rho,\epsilon)\del{1 + \bar Q}\del{C_H\del{P(0)+1} + \bar P}}^{-1},
	\end{align}
	where $C_H$ is the constant from \eqref{eq:C_H smoothness}.
	Then there is at most one solution $(u,m,Q^*)$ of \eqref{eq:mfg T}, and likewise there is at most one solution $(u,m,Q^*)$ of \eqref{eq:mfg infty} such that $u$ and $\pd{u}{x}$ are bounded.
\end{theorem}

\begin{proof}
	Suppose that $(u,m,Q^*)$ and $(\hat u,\hat m,\hat Q^*)$ are both solutions of \eqref{eq:mfg T}, or of \eqref{eq:mfg infty} with $u,\pd{u}{x},\hat u,$ and $\pd{\hat u}{x}$ bounded.
	We will employ the results of Sections \ref{sec:Xn and Xnstar} and \ref{sec:energy estimates}, which are proved independently.
	Equation \eqref{eq:r big for uniqueness} (which is surely an overestimate, see Remark \ref{rmk:kappa if n=0}) implies that Assumption \ref{as:r bigger than kappa} holds.
	Then Equation \eqref{eq:ep small for uniqueness} implies that Lemma \ref{lem:energy differences} holds.
	Since the initial conditions are the same, i.e.~$\hat m_0 = m_0$, we have
	\begin{equation*}
		\int_0^T \int_{\s{D}}  e^{-r t} \abs{\dpd{u}{x}-\dpd{\hat u}{x}}^2\del{m(\dif x,t) + \hat m(\dif x,t)} \dif t = 0,
	\end{equation*}
	where $T$ is the (finite or infinite) time horizon.
	It follows that $\pd{u}{x} = \pd{\hat u}{x}$ on the support of $m$ and $\hat m$, and so by Lemma \ref{lem:market clearing Cournot} we deduce that $Q^* = \hat Q^*$.
	Then by standard uniqueness for parabolic equations, it follows that $m = \hat m$; we also get $u = \hat u$ in a straightforward way if $T < \infty$.
	
	For the infinite time horizon case, let $w(x,t) = e^{-rt}\del{u(x,t) - \hat u(x,t)}$ and note that it satisfies
	\begin{equation}
		-\dpd{w}{t} - \frac{\sigma^2}{2} \dpd[2]{w}{x} = e^{-rt}\del{H\del{\epsilon,Q^*(t),\dpd{u}{x}} - H\del{\epsilon,Q^*(t),\dpd{\hat u}{x}}} \leq C\abs{\dpd{w}{x}},
	\end{equation}
	since $\pd{u}{x}$ and $\pd{\hat u}{x}$ are bounded.
	Let $c > 0$.
	Multiply by $(w - c)_+ := \max\cbr{w-c,0}$ and integrate to get
	\begin{multline}
		\int_0^\infty (w - c)_+(x,t)^2 \dif x + \frac{\sigma^2}{2} \int_t^T \int_0^\infty \abs{\dpd{(w - c)_+}{x}}^2 \dif x \dif \tau\\
		\leq \int_0^\infty (w - c)_+(x,T)^2 \dif x + C\int_t^T \int_0^\infty \abs{\dpd{(w - c)_+}{x}}(w - c)_+ \dif x \dif \tau,
	\end{multline}
	from which we deduce
	\begin{equation}
		\int_0^\infty (w - c)_+(x,t)^2 \dif x
		\leq \int_0^\infty (w - c)_+(x,T)^2 \dif x + C\int_t^T \int_0^\infty (w - c)_+^2 \dif x \dif \tau.
	\end{equation}
	By Gronwall's inequality (applied backward in time), we obtain
	\begin{equation}
		\int_0^\infty (w - c)_+(x,t)^2 \dif x
		\leq e^{C(T-t)}\int_0^\infty (w - c)_+(x,T)^2 \dif x.
	\end{equation}
	Since $u,\hat u$ are bounded, taking $T$ large enough we deduce $w(x,T) \leq c$, and thus the right-hand side is zero.
	We deduce that $w \leq c$ everywhere.
	Since $c$ is arbitrary, it follows that $w \leq 0$, i.e.~$u \leq \hat u$.
	By reversing the roles of $u$ and $\hat u$ we see that $u = \hat u$.
\end{proof}

The following result does not require any of the assumptions made in this section, but simply imposes a linear demand schedule.
\begin{theorem} \label{thm:uniqueness lin demand}
	Under the assumption $P(q) = 1-q$ (but no additional assumptions), there is at most one solution to the finite horizon problem \eqref{eq:mfg T}, and likewise at most one solution to the infinite time horizon problem \eqref{eq:mfg infty} such that $u$ and $\pd{u}{x}$ are bounded.
\end{theorem}

\begin{proof}
	Let $(u,m,Q^*)$ and $(\hat u,\hat m,\hat Q^*)$ be two solutions to the PDE system \eqref{eq:mfg infty}, then set $q^* := q^*(\epsilon,Q^*,\pd{u}{x})$ and $\hat q^* = q^*(\epsilon,\hat Q^*,\pd{\hat u}{x})$.
	Following the calculations in \cite{graber2020commodities}, we derive
	\begin{multline} \label{eq:lin energy differences}
		\int_0^T \int_0^\infty e^{-rt}\del{\hat q^* - q^*}^2 (m + \hat m)\dif x \dif t
		+ \epsilon\int_0^T e^{-rt} \del{Q^*(t) - \hat Q^*(t)}^2 \dif t\\
		\leq \int_0^\infty \del{e^{-rT}\del{u - \hat u}(x,T)\del{m - \hat m}(x,T) - \del{u - \hat u}(x,0)\del{m - \hat m}(x,0)}\dif x.
	\end{multline}
	Because the initial/final data are the same, the right-hand side is zero, and we conclude using the same arguments as in the proof of Theorem \ref{thm:uniqueness small ep}.
\end{proof}

\section{A priori estimates on the linearized system}

\label{sec:sensitivity}


In this section our goal is to prove a priori estimates and existence of solutions for a system of the form
\begin{equation} \label{eq:linearized}
	\begin{cases}
		(i) \ \dpd{w}{t} + \dfrac{\sigma^2}{2}\dpd[2]{w}{x} + V_1(x,t)\dpd{w}{x} + V_2(x,t)\s{Q}(t) - rw = f,\\
		(ii) \ \dpd{\mu}{t} - \dfrac{\sigma^2}{2}\dpd[2]{\mu}{x} + \dpd{}{x}\del{V_3(x,t)\mu} + \dpd{}{x}\del{\del{V_4(x,t)\dpd{w}{x} + V_5(x,t)\s{Q}(t)}m + \nu}= 0,\\
		(iii) \ \s{Q}(t) = \displaystyle  \del{1 + \int_{\s{D}} V_5(\cdot,t)\dif m(t)}^{-1}\del{- \int_{\s{D}} \dif \nu(t) - \int_{\s{D}} V_3(\cdot,t)\dif \mu(t) - \int_{\s{D}} V_4(\cdot,t)\dpd{w}{x}(\cdot,t)\dif m(t)},\\
		(iv) \ \mu|_{x=0} = w|_{x=0} = 0, \
		\mu|_{t=0} = \mu_0.
	\end{cases}
\end{equation}
It is useful to study System \eqref{eq:linearized} at a sufficiently high level of abstraction because our estimates will serve three purposes:
\begin{enumerate}
	\item proving that $U$ is Lipschitz with respect to the measure variable,
	\item proving the existence of a candidate for $\vd{U}{m}$, and
	\item proving that the candidate is indeed a derivative in the sense of Definition \ref{def:derivative}.
\end{enumerate}
To see this, let $(u,m,Q^*)$ and $(\hat u,\hat m,\hat Q^*)$ be the solutions of \eqref{eq:mfg infty} corresponding to initial conditions $m_0$ and $\hat m_0$, respectively.
For $s \in [0,1]$ define
\begin{equation*}
	u_s = s \hat u + (1-s)u,
	\quad
	Q^*_s = s\hat{Q}^* + (1-s)Q^*.
\end{equation*}
If $w = \hat u - u, \mu = \hat m - m$, and $\s{Q} = \hat Q^* - Q^*$, then \eqref{eq:linearized} is satisfied with
\begin{equation} \label{eq:V_i for differences}
	\begin{aligned}
		V_1(x,t) &= \int_0^1 \dpd{H}{a}\del{\epsilon,Q_s^*(t),\dpd{u_s}{x}}\dif s,\\
		V_2(x,t) &= \int_0^1 \dpd{H}{Q}\del{\epsilon,Q_s^*(t),\dpd{u_s}{x}}\dif s,\\
		V_3(x,t) &= \dpd{H}{a}\del{\epsilon,\hat{Q}^*(t),\dpd{\hat u}{x}},\\
		V_4(x,t) &= \int_0^1 \dpd[2]{H}{a}\del{\epsilon,Q_s^*(t),\dpd{u_s}{x}}\dif s, \\
		V_5(x,t) &= \int_0^1 \dmpd{H}{Q}{a}\del{\epsilon,Q_s^*(t),\dpd{u_s}{x}}\dif s,
	\end{aligned}
\end{equation}
with $f = 0$ and $\nu = 0$.

Next, we formally take the derivative of System \eqref{eq:mfg infty} with respect to the measure.
The result is System \eqref{eq:linearized} if we define
\begin{equation} \label{eq:V_i for derivative}
	\begin{cases}
		V_1(x,t) = \dpd{H}{a}\del{\epsilon,Q^*(t),\dpd{u}{x}},\\
		V_2(x,t) = \dpd{H}{Q}\del{\epsilon,Q^*(t),\dpd{u}{x}},\\
		V_3(x,t) = \dpd{H}{a}\del{\epsilon,Q^*(t),\dpd{u}{x}},\\
		V_4(x,t) = \dpd[2]{H}{a}\del{\epsilon,Q^*(t),\dpd{u}{x}},\\
		V_5(x,t) = \dmpd{H}{Q}{a}\del{\epsilon,Q^*(t),\dpd{u}{x}},
	\end{cases}
\end{equation}
with $f = 0$ and $\nu = 0$.
If $(w,\mu)$ is the solution to System \eqref{eq:linearized} assuming \eqref{eq:V_i for derivative} and initial conditions $\mu_0 = \delta_y$, then $w(x,0) = \vd{U}{m}(m,x,y)$ is the candidate derivative of the master field $U(m_0,x)$ with respect to $m_0$, where $m_0$ is a given initial condition in System \eqref{eq:mfg infty}.

Finally, let $\tilde w = \hat u - u - w, \tilde \mu = \hat m - m - \mu, \widetilde{\s{Q}} = \hat Q^* - Q^* - \s{Q}$.
Then $(\tilde w,\tilde \mu,\widetilde{\s{Q}})$ satisfies \eqref{eq:linearized} with $V_1,\ldots,V_5$ defined as in \eqref{eq:V_i for derivative} and with
\begin{multline} \label{eq:f nu}
	f(x,t) = -\int_0^1 \del{\dpd{H}{a}\del{\epsilon,Q_s^*(t),\dpd{u_s}{x}} - \dpd{H}{a}\del{\epsilon,Q^*(t),\dpd{u}{x}}}\del{\dpd{\hat u}{x}-\dpd{u}{x}}\dif s\\
	- \int_0^1 \del{\dpd{H}{Q}\del{\epsilon,Q_s^*(t),\dpd{u_s}{x}} - \dpd{H}{Q}\del{\epsilon,Q^*(t),\dpd{u}{x}}}\del{\hat Q^*(t) - Q^*(t)}\dif s,\\
	\nu(t) = \dmpd{H}{Q}{a}\del{\epsilon,Q^*(t),\dpd{u}{x}}(\hat Q - Q)(\hat m - m)
	+ \dpd[2]{H}{a}\del{\epsilon,Q^*(t),\dpd{u}{x}}\del{\dpd{\hat u}{x}-\dpd{u}{x}}(\hat m - m)\\
	+ \hat m \int_0^1 \del{\dmpd{H}{Q}{a}\del{\epsilon,Q_s^*(t),\dpd{u_s}{x}}
		 - \dmpd{H}{Q}{a}\del{\epsilon,Q^*(t),\dpd{u}{x}}}\del{\hat Q^*(t) - Q^*(t)}\dif s\\
	+ \hat m \int_0^1 \del{\dpd[2]{H}{a}\del{\epsilon,Q_s^*(t),\dpd{u_s}{x}} - \dpd[2]{H}{a}\del{\epsilon,Q^*(t),\dpd{u}{x}}}\del{\dpd{\hat u}{x}-\dpd{u}{x}}\dif s.
\end{multline}
Our a priori estimates on $(\tilde w,\tilde \mu,\widetilde{\s{Q}})$ will allow us to conclude that our candidate satisfies the definition of derivative given in Definition \ref{def:derivative}.

Conceptually, the a priori estimates are organized in the following progression.
A crucial point is to obtain \emph{energy estimates}, which are derived by developing $\od{}{t}\ip{w}{\mu}$ using the equations and isolating positive terms.
However, it was already noticed in \cite{graber2022parameter} that the integral terms appearing in system such as \eqref{eq:linearized} interfere with the energy estimates.
Because of this, we first introduce a set of technical estimates on the Fokker-Planck equation, which require substantial preliminary results on parabolic equations.
Once this major step is accomplished, we are then to proceed to the energy estimates, followed by H\"older regularity in time, and concluded by full Schauder type estimates.
Combining the a priori estimates with the Leray-Schauder fixed point theorem, we also deduce an existence result for System \eqref{eq:linearized}.

\subsection{Preliminaries: global-in-time interior estimates}

\label{sec:interior estimates}



In the context of our study of System \eqref{eq:linearized}, the main purpose of this section is to introduce some function spaces which, together with their \emph{dual} spaces, will be useful for technical reasons in the sequel.
There is a more general motivation, however, which is to find higher-order estimates on parabolic equations with Dirichlet boundary conditions, while bypassing the compatibility conditions on the boundary.
So as not to distract the reader from the main purpose of this section, we have moved all the proofs to the appendix.

\subsubsection{Interior estimates on the heat equation}
\label{sec:interior heat}

Define $d(x) := \min\cbr{x,1}$.
Let $n$ be a non-negative integer and let $k \geq 0$.
For a function $\phi:\intco{0,\infty} \to \bb{R}$, we define the seminorm
\begin{equation}
	\intcc{\phi}_{n,k} := \enVert{d^{n+k}\phi^{(n)}}_0 = \sup_{x \geq 0} d(x)^{n+k} \abs{\phi^{(n)}(x)}
\end{equation}
and the norm
\begin{equation}
	\enVert{\phi}_{n,k} := \max_{0 \leq j \leq n}
	\intcc{\phi}_{j,k}.
\end{equation}
When $k = 0$ we will simply write $\intcc{\phi}_{n,0} = \intcc{\phi}_{n}$ and $\enVert{\phi}_{n,0} = \enVert{\phi}_{n}$.
We will define $X_{n,k}$ to be the space of all function $\phi:\intco{0,\infty} \to \bb{R}$ such that $\enVert{\phi}_{n,k}$ is finite, and $X_{n} := X_{n,0}$.

We will also make use of the following norm:
\begin{equation}
	\enVert{\phi}_{n,1}^* := \sup_{0 \leq x \leq 1}\abs{\int_0^x \phi(\xi)\dif \xi} + \enVert{\phi}_{n,1}.
\end{equation}

Consider now the following potentials:
\begin{equation} \label{eq:potentials}
	\begin{split}
		u(x,t) &= \int_{0}^\infty S(x-y,t)u_0(y)\dif y,\\
		v(x,t) &= \int_{0}^t\int_{0}^\infty S(x-y,t-s)f(y,s)\dif y\dif s,\\
		w(x,t) &= -2\int_{0}^t \dpd{S}{x}(x,t-s)\psi(s)\dif s.
	\end{split}
\end{equation}

\begin{proposition} \label{pr:potentials}
	Let $u_0 \in X_n$, $f \in \s{C}\del{[0,T];X_{n-1,1}}$, and $\psi \in \s{C}([0,T])$.
	Then there exists a constant $M_n$, depending only on $n$, such that for $u,v,w$ defined as in \eqref{eq:potentials}, we have
	\begin{equation} \label{eq:potential estimates}
		\begin{split}
			\enVert{u(\cdot,t)}_n
			&\leq M_n \enVert{u_0}_n,\\
			\enVert{v(\cdot,t)}_n
			&\leq M_n\int_0^t (t-s)^{-1/2}\enVert{f(\cdot,s)}_{n-1,1}^*\dif s,\\
			\enVert{w(\cdot,t)}_n &\leq M_n \sup_{0 \leq s \leq t}\abs{\psi(s)}.
		\end{split}
	\end{equation}
\end{proposition}

\begin{proof}
	See Appendix \ref{ap:proofs interior heat}.
\end{proof}

A corollary of Proposition \ref{pr:potentials} is an estimate of solutions to the Dirichlet problem:
\begin{equation}
	\label{eq:dirichlet}
	\dpd{u}{t} = \frac{\sigma^2}{2}\dpd[2]{u}{x} + f(x,t), \quad u(0,t) = \psi(t), \quad u(x,0) = u_0(x).
\end{equation}
\begin{theorem} \label{thm:dirichlet estimate}
	Let $u_0 \in X_n$, $f \in \s{C}\del{[0,T];X_{n-1,1}}$, and $\psi \in \s{C}([0,T])$.
	Let $u$ be the solution of \eqref{eq:dirichlet}.
	Then there exists a constant $M_n$, depending only on $n$, such that
	\begin{equation}
		\label{eq:dirichlet estimate}
		\enVert{u(\cdot,t)}_n \leq M_n \del{\enVert{u_0}_n + t^{1/2}\sup_{0 \leq s \leq t}\enVert{f(\cdot,s)}_{n-1,1}^* + \sup_{0 \leq s \leq t} \abs{\psi(s)}}.
	\end{equation}
\end{theorem}

\begin{proof}
	See Appendix \ref{ap:proofs interior heat}.
\end{proof}

\subsubsection{Application to MFG system}
\label{sec:interior mfg}

Here and in what follows we will let $n$ be a positive integer such that $P$ is $n+2$ times differentiable; by Assumption \ref{as:P} it is possible to take $n = 2$.
Then we deduce that $H$ is $n+1$ times differentiable.
A corollary of the results in Section \ref{sec:interior heat} is the following:

\begin{proposition}
	Let $(u,m)$ be the solution to the mean field games system on a finite or infinite time horizon $T$, i.e.~either of System \eqref{eq:mfg T} or \eqref{eq:mfg infty}.
	Suppose
	\begin{equation} \label{eq:r assm n}
		r > \max\cbr{(2\bar Q M_n)^2,1} \ln(2M_n),
	\end{equation}
	where $\bar Q$ is defined in Equation \eqref{eq:M} and $M_n$ is the constant from Theorem \ref{thm:dirichlet estimate}.
	Then for any $n$ such that $H$ is $n+1$ times differentiable, we have
	\begin{align} \label{eq:dHda_n}
		\sup_{t \geq 0} \enVert{\dpd{H}{a}\del{\epsilon,Q^*(t),\dpd{u}{x}(\cdot,t)}}_n &\leq D_n(r),\\
		\sup_{t \geq 0} \enVert{\dpd{H}{Q}\del{\epsilon,Q^*(t),\dpd{u}{x}(\cdot,t)}}_n &\leq \epsilon \tilde D_n(r), \label{eq:dHdQ_n}
	\end{align}
	where $D(r),\tilde D_n(r) \geq 1$ are constants that decrease as $r$ increases.
\end{proposition}

\begin{proof}
	See Appendix \ref{ap:proofs interior mfg}.
\end{proof}

\begin{remark}[Constants for $n = 0$] \label{rmk:constants for n=0}
	It is worth noting that in the case $n = 0$, the constants used in this Section are already known. 
	In particular, $M_0 = 1$, $D_0(r) = \bar Q$ (see Equation \eqref{eq:q* a priori bound}), and $\tilde D_0(r) = P(0)$ (see Corollary \ref{cor:dHdQ}).
\end{remark}

\subsection{Assumptions on the data}

We will study \eqref{eq:linearized} on a time horizon $T$ which could be finite or infinite.
When $T < \infty$ we take a final condition $w(x,T) = 0$ and assume that $\epsilon = \epsilon(t)$ satisfies Assumption \ref{as:epsilon(t)}.
We will denote $\epsilon(0) = \epsilon$.
If $T = \infty$ then $\epsilon$ is assumed to be constant, and we assume that
\begin{equation}
	\label{eq:w(T) T =infty}
	\lim_{t \to \infty} e^{-\frac{r}{2}t}\enVert{w(\cdot,t)}_n =
	\lim_{t \to \infty} e^{-\frac{r}{2}t}\enVert{\dpd{w}{x}(\cdot,t)}_n = 0,
	\quad
	t \mapsto e^{-\frac{r}{2}t}\enVert{\mu(t)}_{-n} \quad \text{is bounded.}
\end{equation}
In addition, we will state many of the following results in terms of an arbitrary positive integer $n$, which satisfies the restriction that $P$ is $n + 2$ times differentiable and therefore $H$ is $n + 1$ times differentiable.
Assumption \ref{as:H must be smooth} and Equation \eqref{eq:r assm n} will be in force throughout this section.
Hence Proposition \ref{pr:dudxn estimates} and its corollaries \eqref{eq:dHda_n} and \eqref{eq:dHdQ_n} apply.

We now state assumptions on the coefficients $V_1,\ldots,V_5$, which are abstracted from the particular cases \eqref{eq:V_i for derivative} and \eqref{eq:V_i for differences}.
\begin{assumption}
	\label{as:Vi}
	\begin{enumerate}
		\item $\enVert{V_1(\cdot,t)}_n \leq D_n(r)$ for all $t$, where $D_n(r)$ is the same as in Equation \eqref{eq:dHda_n}, and we assume without loss of generality that $D_n(r) \geq 1$;
		\item $\enVert{V_2(\cdot,t)}_n \leq \epsilon\tilde D_n(r)$ for all $t$, where $\tilde D_n(r)$ is the same as in Equation \eqref{eq:dHdQ_n};
		\item $\enVert{V_3(\cdot,t)}_n \leq D_n(r)$ for all $t$;
		\item $C_H^{-1} \leq V_4(x,t) \leq C_H$ for all $(x,t)$;
		\item $V_5(x,t) \in \intcc{\epsilon \frac{1-\bar \rho}{2-\bar \rho},\epsilon}$ for all $(x,t)$, and thus $\enVert{V_5}_0 \leq \epsilon \max\cbr{\abs{\frac{\bar \rho - 1}{\bar \rho - 2}}, 1} =: \bar P \epsilon$.
	\end{enumerate}
\end{assumption}

\begin{lemma}
	Let $V_1,\ldots,V_5$ be given using formula \eqref{eq:V_i for derivative} or \eqref{eq:V_i for differences}.
	Then Assumption \ref{as:Vi} holds.
\end{lemma}

\begin{proof}
	This follows from Corollaries \ref{cor:smoothness}, \ref{cor:dHdQ},  and \ref{cor:q* a priori bound}; Equations \eqref{eq:dHda_n} and \eqref{eq:dHdQ_n}; and the a priori estimates from Theorem \ref{thm:exist mfg infty}.
\end{proof}

\textbf{Notation:}
If $g = g(y,t)$ is a function depending on $t$ and other variables $y$ and $\rho$ is a real number, we will denote by $g_\rho$ the function
\begin{equation}
	\label{eq:g_rho}
	g_\rho(y,t) = e^{-\rho t}g(y,t).
\end{equation}
The \emph{energy} with parameter $\rho$ is denoted
\begin{equation} \label{eq:energy}
	E_\rho(t) = \int_{\s{D}} \abs{\dpd{w_\rho}{x}(\cdot,t)}^2 \dif m(t) = \int_{\s{D}} e^{-2\rho t}\abs{\dpd{w}{x}(\cdot,t)}^2 \dif m(t).
\end{equation}
This quantity will appear often in our estimates, and we will prove a priori bounds on $\int_0^T E_\rho(t)\dif t$ in Section \ref{sec:energy estimates}.

\subsection{Estimates in $X_n$ and $X_n^*$} \label{sec:Xn and Xnstar}

We will denote by $X_n^*$ the dual of the space $X_n$, and by $\enVert{\cdot}_{-n}$ the dual norm
\begin{equation} \label{eq:-n norm}
	\enVert{\mu}_{-n} = \sup_{\enVert{\phi}_n \leq 1} \ip{\phi}{\mu}.
\end{equation}
Note that $\enVert{\mu}_{-0} = \enVert{\mu}_{TV}$ by the Riesz representation theorem:
\begin{equation}
	\enVert{\mu}_{-0} = \sup_{\enVert{\phi}_0 \leq 1} \int_{\s{D}} \phi(x)\dif \mu(x) = \enVert{\mu}_{TV}.
\end{equation}
In this subsection we provide a priori estimates on $\mu(t)$ in $X_n^*$, where $(w,\mu)$ is a solution of the linearized system.
First, we introduce a technical lemma, somewhat reminiscent of Gr\"onwall's inequality.
Cf.~\cite[Lemma 2.1]{graber2022parameter}.
\begin{lemma} \label{lem:int estimate}
	Let $A,B,\delta > 0$ be given constants.
	Suppose $f,g:\intco{0,\infty} \to \intco{0,\infty}$ are functions that satisfy
	\begin{equation} \label{eq:ineq}
		f(t_1) \leq Af(t_0) + \int_{t_0}^{t_1}(t_1-s)^{-1/2}\del{B f(s) + g(s)}\dif s \quad \forall 0 \leq t_0 \leq t_1 \leq t_0 + \delta
	\end{equation}
	Then for any $\lambda > \frac{1}{\delta}\ln(A)$, we have
	\begin{equation} \label{eq:estimate}
		\del{1 - \frac{2\delta^{1/2}B}{1 - Ae^{-\lambda \delta}}}\int_0^T e^{-\lambda t}f(t)\dif t \leq \frac{A}{\lambda - \delta^{-1}\ln(A)}f(0) 
		+
		\frac{2\delta^{1/2}}{1 - Ae^{-\lambda \delta}}\int_0^T e^{-\lambda t}g(t) \dif t.
	\end{equation}
\end{lemma}

\begin{proof}
	See Appendix \ref{ap:int estimate}.
\end{proof}

\begin{lemma} \label{lem:mu-n int}
	Let $(w,\mu)$ be a solution of \eqref{eq:linearized}.
	Fix $\rho \geq \kappa(r)$, where
	\begin{equation} \label{eq:kappa}
		\kappa(r) := 32\del{1 + c(\bar \rho,\epsilon)\bar P \epsilon}^2D_n(r)^2M_n^2\ln(8M_n^2).
	\end{equation}
	Then we have
	\begin{equation} \label{eq:mu-n int}
		\int_0^T \enVert{\mu_\rho(t)}_{-n}^2 \dif t \leq \enVert{\mu_0}_{-n}^2
		+ \int_0^T \del{\enVert{V_4}_0^2 E_\rho(s) + \enVert{\nu_\rho(s)}_{-n}^2}\dif s.
	\end{equation}
\end{lemma}

\begin{proof}
	\firststep
	Fix $t_1 > t_0 \geq 0$ and let $\phi_1 \in X_n$.
	Define $\phi$ to be the solution of the Dirichlet problem
	\begin{equation}
		-\dpd{\phi}{t} = \frac{\sigma^2}{2}\dpd[2]{\phi}{x}, \quad \phi(0,t) = 0, \quad \phi(x,T) = \phi_1(x).
	\end{equation}
	By the reflection principle, a formula for $\phi$ is
	\begin{equation}
		\phi(x,t) = \int_0^\infty S(x-y,t_1-t)\phi_1(y)\dif y.
	\end{equation}
	By applying Theorem \ref{thm:dirichlet estimate}, we get
	\begin{equation} \label{eq:phi-est1}
		\enVert{\phi(\cdot,t)}_n \leq M_n\enVert{\phi_1}_n \quad \forall t \in [0,t_1].
	\end{equation}
	Moreover, by the same argument as in \ref{pr:potentials}, we get
	\begin{equation} \label{eq:phi-est2}
		\enVert{\dpd{\phi}{x}(\cdot,t)}_n \leq M_n\enVert{\phi_1}_n (t_1-t)^{-1/2}.
	\end{equation}
	Now use $\phi$ as a test function in \eqref{eq:linearized}(ii) to get
	\begin{multline}
		\ip{\phi_1}{\mu(t_1)} = \ip{\phi(t_0)}{\mu(t_0)} + \int_{t_0}^{t_1} \ip{\dpd{\phi}{x}(\cdot,t)V_3(\cdot,t)}{\mu(t)}\dif t\\
		- \int_{t_0}^{t_1} \ip{\dpd{\phi}{x}(\cdot,t)}{\del{V_4(\cdot,t)\dpd{w}{x}(\cdot,t) + V_5(\cdot,t)\s{Q}(t)}m(t) + \nu(t)}\dif t.
	\end{multline}
	Applying \eqref{eq:phi-est1} and \eqref{eq:phi-est2} as well as the Cauchy-Schwartz inequality, recalling that $\enVert{m(t)}_{TV} \leq 1$, we get
	\begin{multline} \label{eq:phi1 mu}
		\abs{\ip{\phi_1}{\mu(t_1)}} \leq M_n \enVert{\phi_1}_n \enVert{\mu(t_0)}_{-n}
		+ M_n \enVert{\phi_1}_n \int_{t_0}^{t_1} (t_1-t)^{-1/2}\enVert{V_3(\cdot,t)}_n\enVert{\mu(t)}_{-n}\dif t\\
		+ M_n\enVert{\phi_1}_n\int_{t_0}^{t_1} (t_1-t)^{-1/2}\del{\enVert{V_4}_0 E_0(t)^{1/2} + \enVert{V_5}_0 \abs{\s{Q}(t)} + \enVert{\nu(t)}_{-n}}\dif t.
	\end{multline}
	
	\nextstep
	Next, we need to estimate $\s{Q}(t)$ using \eqref{eq:linearized}(iii).
	We get
	\begin{equation} \label{eq:abs Q}
		\abs{\s{Q}(t)} \leq c(\bar \rho,\epsilon)\del{\enVert{\nu(t)}_{-n} + \enVert{V_3(\cdot,t)}_n \enVert{\mu(t)}_{-n} + \enVert{V_4}_0 E_0(t)^{1/2}}.
	\end{equation}
	Plugging \eqref{eq:abs Q} into \eqref{eq:phi1 mu} and using Assumption \ref{as:Vi}, we deduce
	\begin{multline} \label{eq:phi1 mu 2}
		\abs{\ip{\phi_1}{\mu(t_1)}} \leq M_n \enVert{\phi_1}_n \enVert{\mu(t_0)}_{-n}
		+ \del{1 + c(\bar \rho,\epsilon)\bar P \epsilon}D_n(r)M_n \enVert{\phi_1}_n \int_{t_0}^{t_1} (t_1-t)^{-1/2}\enVert{\mu(t)}_{-n}\dif t\\
		+ \del{1 + c(\bar \rho,\epsilon)\bar P \epsilon}M_n\enVert{\phi_1}_n\int_{t_0}^{t_1} (t_1-t)^{-1/2}\del{\enVert{V_4}_0 E_0(t)^{1/2} + \enVert{\nu(t)}_{-n}}\dif t.
	\end{multline}
	Taking the supremum over all $\phi_1 \in X_n$, we get
	\begin{multline} \label{eq:mu-n 1}
		\enVert{\mu(t_1)}_{-n} \leq M_n  \enVert{\mu(t_0)}_{-n}
		+ \del{1 + c(\bar \rho,\epsilon)\bar P \epsilon}D_n(r)M_n  \int_{t_0}^{t_1} (t_1-t)^{-1/2}\enVert{\mu(t)}_{-n}\dif t\\
		+ \del{1 + c(\bar \rho,\epsilon)\bar P \epsilon}M_n\int_{t_0}^{t_1} (t_1-t)^{-1/2}\del{\enVert{V_4}_0 E_0(t)^{1/2} + \enVert{\nu(t)}_{-n}}\dif t, \quad \forall 0 \leq t_0 < t_1.
	\end{multline}
	\nextstep 
	Square both sides of \eqref{eq:mu-n 1} and use Cauchy-Schwartz to get
	\begin{multline} \label{eq:mu-n 2}
		\enVert{\mu(t_1)}_{-n}^2 \leq 4M_n^2 \enVert{\mu(t_0)}_{-n}^2
		+ \tilde B (t_1-t_0)^{1/2}\int_{t_0}^{t_1} (t_1-t)^{-1/2}\enVert{\mu(t)}_{-n}^2\dif t\\
		+ \tilde B (t_1-t_0)^{1/2}\int_{t_0}^{t_1} (t_1-t)^{-1/2}\del{\enVert{V_4}_0^2 E_0(t) + \enVert{\nu(t)}_{-n}^2}\dif t, \quad \forall 0 \leq t_0 < t_1
	\end{multline}
	where $\tilde B := 8\del{1 + c(\bar \rho,\epsilon)\bar P \epsilon}^2D_n(r)^2M_n^2$.
	Now we will apply Lemma \ref{lem:int estimate} with
	\begin{equation*}
		\begin{split}
			&A = 4M_n^2,
			\ B = \tilde B \delta^{1/2}, \
			\delta = (8\tilde B)^{-1}, \
			f(t) = \enVert{\mu(t)}_n^2,\\\
			&g(t) = B\del{\enVert{V_4}_0^2 E_0(t) + \enVert{\nu(t)}_{-n}^2},\ \text{and} \
			\lambda = 2\rho.
		\end{split}
	\end{equation*}
	Comparing the definition in Equation \eqref{eq:kappa}, we see that
	\begin{equation}
		\lambda \geq 2\kappa(r) = \delta^{-1}\ln(2A) > \delta^{-1}\ln(A)
		\quad \Rightarrow \quad
		1 - Ae^{-\lambda \delta} \leq \frac{1}{2}.
	\end{equation}
	We also have $2\delta^{1/2}B = 2\delta \tilde B \leq 1/4$, and thus \eqref{eq:estimate} implies
	\begin{equation}
		\frac{1}{2}\int_0^T e^{-\lambda t}f(t)\dif t \leq \frac{A}{\delta^{-1}\ln(2)}f(0) 
		+
		4\delta^{1/2}\int_0^T e^{-\lambda t}g(t) \dif t.
	\end{equation}
	By comparing the constants defined above, we deduce
	\begin{equation}
		\int_0^T e^{-\lambda t}f(t)\dif t \leq f(0) 
		+
		B^{-1}\int_0^T e^{-\lambda t}g(t) \dif t,
	\end{equation}
	which implies \eqref{eq:mu-n int}, as desired.
\end{proof}

\begin{corollary}
	\label{cor:int sQ estimate}
	Let $(w,\mu)$ be a solution of \eqref{eq:linearized}, and suppose $\rho \geq \kappa(r)$ with $\kappa(r)$ defined in \eqref{eq:kappa}.
	Then
	\begin{equation} \label{eq:int sQ estimate}
		\del{\int_0^T \abs{\s{Q}_\rho(t)}^2 \dif t}^{1/2} \leq \hat D_n(r)\del{\enVert{\mu_0}_{-n} + \del{\int_0^T \enVert{\nu_\rho(t)}_{-n}^2 \dif t}^{1/2}
			+ \del{\int_0^T E_\rho(t) \dif t}^{1/2}}
	\end{equation}
	where $\hat D_n(r) = c(\bar \rho,\epsilon)\del{1+D_n(r)}$.
\end{corollary}

\begin{proof}
	Multiply \eqref{eq:abs Q} by $e^{-\rho t}$, take the $L^2(0,T)$ norm and then apply Lemma \ref{lem:mu-n int}.
\end{proof}

\begin{lemma} \label{lem:w_rho int}
	Let $(w,\mu)$ be a solution of \eqref{eq:linearized} with time horizon $T$.
	There exists a constant $\kappa_1(r)$, which depends only on $n,\sigma$, and $r$ and is decreasing with respect to $r$, such that if
	\begin{equation}
		\label{eq:rho upper bound}
		\rho \leq r - \kappa_1(r)
	\end{equation}
	and if
	\begin{equation} \label{eq:dwrho goes to 0}
		\enVert{\dpd{w_\rho}{x}(\cdot,t)}_n \to 0 \quad \text{as} \quad t \to T,
	\end{equation}
	then the following a estimate holds:
	\begin{equation} \label{eq:int dw_rdx}
		\int_0^T \enVert{\dpd{w_\rho}{x}(\cdot,t)}_n^2\dif t
		\leq \hat D_n(r)^2\enVert{\mu_0}_{-n}^2 + \hat D_n(r)^2\int_{0}^{T}  \del{E_\rho(t) + \enVert{\nu_\rho(t)}_{-n}^2
			+ \enVert{f_\rho(\cdot,t)}_n^2}\dif t,
	\end{equation}
	where $\hat D_n(r) = c(\bar \rho,\epsilon)\del{1+D_n(r)}$.
\end{lemma}

\begin{proof}
	\firststep
	Fix some $T' < T$, where $T \in \intoc{0,\infty}$ is the time horizon.
	For any function $g = g(y,t)$ depending on $t$ and possibly other variables, let $\tilde g(y,t) = g(y,T'-t)$.
	By reversing time in Equation \eqref{eq:linearized}(i), we see that $\tilde w_r$ satisfies
	\begin{equation}
		\dpd{\tilde w_r}{t} = \frac{\sigma^2}{2}\dpd[2]{\tilde w_r}{x} + \tilde V_1\dpd{\tilde w_r}{x} + \tilde V_2 \widetilde{\s{Q}_r}(t) - \tilde f_r.
	\end{equation}
	Since $\tilde w_r(0,t) = 0$, we have
	\begin{multline} \label{eq:tilde w_r formula}
		\tilde w_r(x,t) = \int_0^\infty G_{1}(x,y,t-t_0)\tilde w_r(y,t_0)\dif y\\
		+ \int_{t_0}^t \int_0^\infty G_{1}(x,y,t-s)\del{\tilde V_1(y,s)\dpd{\tilde w_r}{y}(y,s) + \tilde V_2(y,s) \widetilde{\s{Q}_r}(s) - \tilde f_r(y,s)}\dif y \dif s
		\quad \forall t \geq t_0 \geq 0
	\end{multline}
	where we define
	\begin{equation}
		G_{(-1)^n}(x,y,t) = (-1)^n S(x-y,t) - S(x+y,t).
	\end{equation}
	Using an argument similar to the proof of Theorem \ref{thm:dirichlet estimate}, we deduce
	\begin{multline} \label{eq:tilde dw_rdx 2}
		\enVert{\dpd{\tilde w_r}{x}(\cdot,t)}_n
		\leq A_n \enVert{\dpd{\tilde w_r}{x}(\cdot,t_0)}_n\\
		+ B_{n,\sigma}\del{D_n(r) + \epsilon\tilde D_n(r) + 1} \int_{t_0}^t (t-s)^{-1/2}\del{\enVert{\dpd{\tilde w_r}{x}(\cdot,s)}_n
			+ \abs{\widetilde{\s{Q}_r}(s)} 
			+ \enVert{\tilde f_r(\cdot,s)}_n}\dif s, 
	\end{multline}
	where $A_n$ depends only on the constants $m_1,\ldots,m_n$, $B_{n,\sigma}$ depends only on the constants $m_{1,\sigma},\ldots,m_{n,1}$, and $\tilde D_n(r)$ is the constant from \eqref{eq:dHdQ_n}.
	
	\nextstep
	Square both sides of \eqref{eq:tilde dw_rdx 2} to get
	\begin{multline} \label{eq:tilde dw_rdx 3}
		\enVert{\dpd{\tilde w_r}{x}(\cdot,t)}_n^2
		\leq \tilde A_n \enVert{\dpd{\tilde w_r}{x}(\cdot,t_0)}_n^2\\
		+ \tilde B_n(t-t_0)^{1/2} \int_{t_0}^t (t-s)^{-1/2}\del{\enVert{\dpd{\tilde w_r}{x}(\cdot,s)}_n^2
			+ \abs{\widetilde{\s{Q}_r}(s)} ^2
			+ \enVert{\tilde f_r(\cdot,s)}_n^2}\dif s, 
	\end{multline}
	where
	\begin{equation}
		\tilde A_n := 4A_n^2,
		\quad 
		\tilde B_{n} := 8B_{n,\sigma}^2\del{D_n(r) + \epsilon\tilde D_n(r) + 1}^2.
	\end{equation}
	We will apply Lemma \ref{lem:int estimate} with
	\begin{equation*}
		\delta = (8\tilde B_n)^{-1}, A = \tilde A_n, \ B = \tilde B_n \delta^{1/2}, \ g(t) = B\del{\abs{\widetilde{\s{Q}_r}(s)} ^2
			+ \enVert{\tilde f_r(\cdot,s)}_n^2}.
	\end{equation*}
	\begin{comment}
	For every $\lambda > \delta^{-1}\ln(\tilde A_n)$, we get
	\begin{multline*} 
		\del{1 - \frac{2\delta\tilde B_n}{1 - \tilde A_n e^{-\lambda \delta}}}\int_0^T e^{-\lambda t}\enVert{\dpd{\tilde w_r}{x}(\cdot,t)}_n^2\dif t
		\\
		 \leq \frac{\tilde A_n}{\lambda - \delta^{-1}\ln(\tilde A_n)}\enVert{\dpd{\tilde w_r}{x}(\cdot,0)}_n^2
		+
		\frac{2\delta \tilde B_n}{1 - \tilde A_ne^{-\lambda \delta}}\int_0^T e^{-\lambda t}\del{\abs{\widetilde{\s{Q}_r}(s)} ^2
			+ \enVert{\tilde f_r(\cdot,s)}_n^2} \dif t.
	\end{multline*}
	\end{comment}
	We deduce that for every $\lambda \geq \delta^{-1}\ln(2\tilde A_n)$,
	\begin{equation} 
		\int_0^{T'} e^{-\lambda t}\enVert{\dpd{\tilde w_r}{x}(\cdot,t)}_n^2\dif t
		\leq \frac{\tilde A_n}{4\tilde B_n\ln(2)}\enVert{\dpd{\tilde w_r}{x}(\cdot,0)}_n^2
		+
		\int_0^{T'} e^{-\lambda t}\del{\abs{\widetilde{\s{Q}_r}(t)} ^2
			+ \enVert{\tilde f_r(\cdot,t)}_n^2} \dif t.
	\end{equation}
	Define
	\begin{equation}
		\kappa_1(r) := 4\tilde B_n \ln(2 \tilde A_n) =
		32B_n^2\del{D_n(r) + \epsilon \tilde D_n(r) + 1}^2 \ln(2 \tilde A_n),
	\end{equation}
	which satisfies the hypotheses given in the statement of the lemma.
	Then set $\rho = r - \frac{\lambda}{2}$; we have define $\kappa_1(r)$ so that $\rho \leq r- \kappa_1(r)$ is equivalent to $\lambda \geq  \delta^{-1}\ln(2\tilde A_n)$.
	Now make the substitution $t \mapsto T' - t$, then let $T' \to T$ and use \eqref{eq:dwrho goes to 0} to get
	\begin{equation} 
		\int_0^{T} \enVert{\dpd{ w_\rho}{x}(\cdot,t)}_n^2\dif t
		\leq \int_0^{T} \del{\abs{\s{Q}_\rho(t)} ^2
			+ \enVert{f_\rho(\cdot,t)}_n^2} \dif t.
	\end{equation}
	Finally, we use Corollary \ref{cor:int sQ estimate} to get \eqref{eq:int dw_rdx}.
\end{proof}

We can also estimate $\enVert{\mu_\rho(t)}_{-n}$ pointwise, provided we are willing to include some dependence on $\enVert{\dpd{w_\rho}{x}}_0$, which will be estimated below.
\begin{lemma} \label{lem:mu-n0}
	Let $(w,\mu)$ be a solution of \eqref{eq:linearized}.
	Suppose
	\begin{equation} \label{eq:kappa_0}
		\rho \geq 
		36\del{1 + c(\bar \rho,\epsilon)}^2D_n(r)^2M_n^2 =: \kappa_0(r).
	\end{equation}
	Then
	\begin{equation} \label{eq:mu-n 0}
		\sup_{0 \leq t \leq T}\enVert{\mu_{\rho}(t)}_{-n} \leq 2M_n \enVert{\mu_0}_{-n}
		\\
		+ \sup_{0 \leq \tau \leq T}\enVert{\nu_{\rho}(\tau)}_{-n}
		+ C_n \enVert{\dpd{w_\rho}{x}}_0^{1/2}\del{\int_{0}^{T}E_{\rho}(s)\dif s}^{1/4}.
	\end{equation}
	where
	\begin{equation}
		C_n = 4\del{1 + c(\bar \rho,\epsilon)}^{1/2}M_n^{1/2}.
	\end{equation}
	\begin{comment}
	If $\rho \leq \kappa_1(r)$, then
	\begin{equation}
	\sup_{t \geq 0} \enVert{w_{\tilde \rho}(t)}_{-n} \leq
	...
	\end{equation}
	\end{comment}
\end{lemma}

\begin{proof}
	Take \eqref{eq:mu-n 1} with $t_0 = 0, t_1 = t$, multiply by $e^{-{\rho} t}$ to get
	\begin{multline} \label{eq:mu-n 01}
		\enVert{\mu_{\rho}(t)}_{-n} \leq M_n \enVert{\mu_0}_{-n}
		+ \del{1 + c(\bar \rho,\epsilon)} D_n(r)M_n  \int_{0}^{t} e^{-\rho(t-s)}(t-s)^{-1/2}\enVert{\mu_{\rho}(s)}_{-n}\dif s\\
		+ \del{1 + c(\bar \rho,\epsilon)}M_n\int_{0}^{t}e^{-\rho(t-s)} (t-s)^{-1/2}\del{E_{\rho}(s)^{1/2} + \enVert{\nu_{\rho}(s)}_{-n}}\dif s.
	\end{multline}
	We first use H\"older's inequality to estimate
	\begin{equation} \label{eq:mu-n 01a}
		\begin{split}
			\int_0^t e^{-\rho(t-s)} (t-s)^{-1/2}E_{\rho}(s)^{1/2}\dif s
			&\leq \enVert{\dpd{w_\rho}{x}}_0^{1/2}\int_0^t e^{-\rho(t-s)} (t-s)^{-1/2}E_{\rho}(s)^{1/4}\dif s\\
			&\leq \enVert{\dpd{w_\rho}{x}}_0^{1/2}\del{\int_0^t e^{-\frac{4}{3}\rho(t-s)} (t-s)^{-2/3}\dif s}^{3/4}
			\del{\int_0^t E_{\rho}(s)\dif s}^{1/4}.
		\end{split}
	\end{equation}
	Using the substitution $s \mapsto t-\frac{s}{\rho}$, we find
	\begin{equation}
		\label{eq:mu-n 01b}
		\begin{split}
			\int_0^t e^{-\frac{4}{3}\rho(t-s)} (t-s)^{-2/3}\dif s
			&= \rho^{-1/3}\int_0^t e^{-\frac{4}{3}s} s^{-2/3}\dif s\\
			&\leq \rho^{-1/3}\del{\int_0^1 s^{-2/3}\dif s + \int_1^\infty e^{-\frac{4}{3}s} \dif s}
			\leq 4\rho^{-1/3}
		\end{split}
	\end{equation}
	and also
	\begin{equation} \label{eq:mu-n 01c}
		\int_0^t e^{-\rho(t-s)} (t-s)^{-1/2}\dif s
		\leq 3\rho^{-1/2}.
	\end{equation}
	Using \eqref{eq:mu-n 01a}, \eqref{eq:mu-n 01b}, and \eqref{eq:mu-n 01c} in \eqref{eq:mu-n 01}, we get
	\begin{multline} \label{eq:mu-n 02}
		\enVert{\mu_{\rho}(t)}_{-n} \leq M_n \enVert{\mu_0}_{-n}
		+ 3\rho^{-1/2}\del{1 + c(\bar \rho,\epsilon)}D_n(r)M_n  \sup_{0 \leq \tau \leq T}\enVert{\mu_{\rho}(\tau)}_{-n}\\
		+ 3\rho^{-1/2}\del{1 + c(\bar \rho,\epsilon)}M_n \sup_{0 \leq \tau \leq T}\enVert{\nu_{\rho}(\tau)}_{-n}
		+ 4\rho^{-1/4}\del{1 + c(\bar \rho,\epsilon)}M_n \enVert{\dpd{w_\rho}{x}}_0^{1/2}\del{\int_{0}^{T}E_{\rho}(s)\dif s}^{1/4}.
	\end{multline}
	By the assumption \eqref{eq:kappa_0}, \eqref{eq:mu-n 02} simplifies to
	\begin{multline} \label{eq:mu-n 03}
		\enVert{\mu_{\rho}(t)}_{-n} \leq M_n \enVert{\mu_0}_{-n}
		+ \frac{1}{2}  \sup_{0 \leq \tau \leq T}\enVert{\mu_{\rho}(\tau)}_{-n}\\
		+ \frac{1}{2} \sup_{0 \leq \tau \leq T}\enVert{\nu_{\rho}(\tau)}_{-n}
		+ 2\del{1 + c(\bar \rho,\epsilon)}^{1/2}M_n^{1/2} \enVert{\dpd{w_\rho}{x}}_0^{1/2}\del{\int_{0}^{T}E_{\rho}(s)\dif s}^{1/4}.
	\end{multline}
	Take the supremum and rearrange to deduce \eqref{eq:mu-n 0}.

	\begin{comment}
	\nextstep
	By \eqref{eq:w_r est 1} with $t_0 = 0$, using the estimates that follow in Step 1 of the proof of Lemma \ref{lem:w_rho int}, we deduce
	\begin{equation} \label{eq:w_r est 0-1}
	\enVert{\tilde w_r(\cdot,t)}_n \leq A_n\int_{0}^t (t-s)^{-1/2}\enVert{\tilde w_r(\cdot,s)}_{n} \dif s
	+ B_n\int_{0}^t \del{\abs{\widetilde{\s{Q}_r}(s)} + \enVert{\tilde f_r(\cdot,s)}_n}\dif s
	\quad \forall t \geq  0.
	\end{equation}
	for some constants $A_n$ and $B_n$.
	Here we recall that $\tilde w_r(x,t) = w_r(x,T-t)$, etc.
	We rewrite \eqref{eq:w_r est 0-1} using the transformations $t \mapsto T-t$ and $s \mapsto T-s$.
	We have
	\begin{equation} \label{eq:w_r est 0-2}
	\enVert{w_r(\cdot,t)}_n \leq A_n\int_{t}^{T} (s-t)^{-1/2}\enVert{w_r(\cdot,s)}_{n} \dif s
	+ B_n\int_{t}^T \del{\abs{\s{Q}_r(s)} + \enVert{f_r(\cdot,s)}_n}\dif s
	\quad \forall t \geq  0.
	\end{equation}
	Multiply both sides of \eqref{eq:w_r est 0-2} by $e^{(r-\tilde \rho)t}$ to get
	\end{comment}
\end{proof}

From now on we make the following assumption:
\begin{assumption}
	\label{as:r bigger than kappa}
	We assume $r \geq 2\max\cbr{\kappa(r),\kappa_1(r),\kappa_0(r)}$ with $\kappa(r)$ defined in \eqref{eq:kappa}, $\kappa_1(r)$ defined in \eqref{eq:rho upper bound}, and $\kappa_0(r)$ defined in \eqref{eq:kappa_0}.
\end{assumption}
Importantly, Assumption \ref{as:r bigger than kappa} can always be obtained by choosing $r$ large enough, because $\kappa(r)$, $\kappa_1(r)$, and $\kappa_0(r)$ are all decreasing functions of $r$.
\begin{remark}
	\label{rmk:kappa if n=0}
	When $n = 0$, Remark \ref{rmk:constants for n=0} shows us that $\kappa(r),\kappa_1(r),$ and $\kappa_0(r)$ no longer depend on $r$.
	In fact, they have the following formulas, more or less explicit:
	\begin{equation} 
		\kappa(r) := 32\del{1 + c(\bar \rho,\epsilon)\bar P \epsilon}^2\bar Q^2\ln(8),
	\end{equation}
	\begin{equation}
		\kappa_1(r) := 
		32B_0^2\del{\bar Q + \epsilon P(0) + 1}^2 \ln(2 \tilde A_0),
	\end{equation}
	\begin{equation} 
		\kappa_0(r) := 
		36\del{1 + c(\bar \rho,\epsilon)}^2\bar Q^2.
	\end{equation}
	Only the constant $B_0$ and $\tilde A_0$ from the proof of Lemma \ref{lem:w_rho int} are left undefined, but upon inspection of the proof we can see that $\tilde A_0$ and $B_0$ are constants no greater than, say, 10.
	Therefore \eqref{eq:r big for uniqueness} is surely an overestimate.
\end{remark}

\begin{corollary}
	[Summary of this subsection] \label{cor:Xn Xnstar est}
	Let $(w,\mu)$ be a solution of \eqref{eq:linearized}.
	Under Assumption \ref{as:r bigger than kappa}, we have the following a priori estimates:
	\begin{equation*}
		\begin{split}
			\int_0^T \enVert{\mu_{r/2}(t)}_{-n}^2 \dif t &\leq \enVert{\mu_0}_{-n}^2
			+ \int_0^T \del{\enVert{V_4}_0^2 E_{r/2}(t) + \enVert{\nu_{r/2}(t)}_{-n}^2}\dif t,\\
			\del{\int_0^T \abs{\s{Q}_{r/2}(t)}^2 \dif t}^{1/2} &\leq \hat D_n(r)\del{\enVert{\mu_0}_{-n} + \del{\int_0^T \enVert{\nu_{r/2}(t)}_{-n}^2 \dif t}^{1/2}
				+ \del{\int_0^T E_{r/2}(t) \dif t}^{1/2}},\\
			\int_0^T \enVert{\dpd{w_{r/2}}{x}(\cdot,t)}_n^2\dif t
			&\leq \hat D_n(r)^2\enVert{\mu_0}_{-n}^2 + \hat D_n(r)^2\int_{0}^{T}  \del{E_{r/2}(t) + \enVert{\nu_{r/2}(t)}_{-n}^2
				+ \enVert{f_{r/2}(\cdot,t)}_n^2}\dif t,\\
			\sup_{0 \leq t \leq T}\enVert{\mu_{{r/2}}(t)}_{-n} &\leq 2M_n \enVert{\mu_0}_{-n}
			+ \sup_{0 \leq \tau \leq T}\enVert{\nu_{{r/2}}(\tau)}_{-n}
			+ C_n \enVert{\dpd{w_{r/2}}{x}}_0^{1/2}\del{\int_{0}^{T}E_{{r/2}}(t)\dif t}^{1/4},
		\end{split}
	\end{equation*}
	where $\hat D_n(r) = c(\bar \rho,\epsilon)\del{1 + D_n(r)}$, $C_n = 4\del{1 + c(\bar \rho,\epsilon)}^{1/2}M_n^{1/2}$, and $D_n(r)$ is the constant appearing in Equation \eqref{eq:dHda_n}.
\end{corollary}

\begin{proof}
	It suffices to observe that the hypotheses of Lemmas \ref{lem:mu-n int}, \ref{lem:w_rho int}, and \ref{lem:mu-n0} are all satisfied with $\rho = r/2$.
\end{proof}

\subsection{Energy estimates} \label{sec:energy estimates}
In some mean field games, known as ``potential mean field games," the Nash equilibrium can be computed by minimizing a certain energy functional \cite{lasry07,cardaliaguet2015weak,cardaliaguet2014mean,cardaliaguet2015second}.
Because of a formal resemblance, we keep the name ``energy estimates" for the estimates derived in this subsection.
We divide our results into two lemmas.
The first deals with differences of solutions to System \eqref{eq:mfg infty}, in which case we assume \eqref{eq:V_i for differences} with $f = \nu = 0$, and the second deals with the case \eqref{eq:V_i for derivative}, with no restriction on $f,\nu$.
Although it is tempting to view the former as a special case of the latter, there are technical points in the proof in which it is not convenient to do so, and thus the proofs are treated separately.
Nevertheless, their basic outline is similar: differentiate the duality pairing $\ip{w}{\mu}$ with respect to time and use the PDE system to write an identity, then use the assumption on the uniform convexity of $H$ to derive an estimate of the integral $\int_0^T E_{r/2}(t) \dif t$.
(Recall that $E_{r/2}$ is defined by \eqref{eq:energy}.)

\begin{lemma}[Energy estimates, differences]
	\label{lem:energy differences}
	Let $(u,m,Q^*)$ and $(\hat u,\hat m,\hat Q^*)$ be solutions to System \eqref{eq:mfg infty} with initial conditions $m_0$ and $\hat m_0$, respectively.
	\begin{enumerate}
		\item Assume that $\epsilon$ satisfies the smallness condition
		\begin{equation} \label{eq:epsilon small differences}
			4C_H \hat D_n(r)\del{C_H\del{P(0)+1} + \bar P}\epsilon \leq 1.
		\end{equation}
		where, as in \eqref{eq:dHdQ estimate}, $\bar P = \max\cbr{\frac{\bar \rho - 1}{\bar \rho - 2},1}$.
		Then 
		\begin{multline} \label{eq:energy differences final}
			\int_0^T \int_{\s{D}}  e^{-r t} \abs{\dpd{u}{x}-\dpd{\hat u}{x}}^2\del{m(\dif x,t) + \hat m(\dif x,t)} \dif t \\
			\leq  \enVert{\hat m_0 - m_0}_{-n}^2
			+ 2C_H \enVert{\hat u(\cdot,0) - u(\cdot,0)}_n \enVert{\hat m_0 - m_0}_{-n},
		\end{multline}
		\item Assume instead that the demand schedule is linear, i.e.~$P(q) = 1-q$, and that $\epsilon < 2$.
		Then we have
		\begin{equation} \label{eq:energy differences linear}
			\int_0^T \int_{\s{D}}  e^{-r t} \abs{\dpd{u}{x}-\dpd{\hat u}{x}}^2\del{m(\dif x,t) + \hat m(\dif x,t)} \dif t 
			\leq  8\enVert{\hat u(\cdot,0) - u(\cdot,0)}_n \enVert{\hat m_0 - m_0}_{-n}.
		\end{equation}
	\end{enumerate}
\end{lemma}

\begin{proof}
	\firststep
	\textit{For a small parameter $\epsilon$.}
	In this first step, we make no further assumptions on the demand schedule $P$ but instead assume condition \eqref{eq:epsilon small differences} holds.
	Multiply $\eqref{eq:mfg infty}_i$(ii) by $u-\hat u$ and integrate by parts, then subtract.
	(See \cite[Theorem 2.4]{lasry07}.)
	After rearranging we get
	{\small \begin{equation} \label{eq:energy differences1}
			\begin{split}
				&\intcc{\int_{\s{D}} e^{-r t}(u(x,t)-\hat u(x,t))(m-\hat m)(\dif x,t)}_0^T
				\\
				&= \int_0^T \int_{\s{D}} e^{-r t}\del{H\del{\epsilon,\hat Q^*(t),\dpd{\hat u}{x}}-H\del{\epsilon,Q^*(t),\dpd{u}{x}}  -\dpd{H}{a}\del{\epsilon,Q^*(t),\dpd{u}{x}} \del{\dpd{\hat u}{x}-\dpd{u}{x}}}m(\dif x,t) \dif t\\
				& \ \ + \int_0^T \int_{\s{D}}  e^{-r t}\del{H\del{\epsilon,Q^*(t),\dpd{u}{x}}-H\del{\epsilon,\hat Q^*(t),\dpd{\hat u}{x}}  -\dpd{H}{a}\del{\epsilon,\hat Q^*(t),\dpd{\hat u}{x}} \del{\dpd{u}{x}-\dpd{\hat u}{x}}}\hat m(\dif x,t)  \dif t.
			\end{split}
	\end{equation}}
	By Equation \eqref{eq:C_H}, we deduce
	\begin{multline} \label{eq:energy differences2}
		\frac{1}{C_H}\int_0^T \int_{\s{D}}  e^{-r t} \abs{\dpd{u}{x}-\dpd{\hat u}{x}}^2\del{m(\dif x,t) + \hat m(\dif x,t)} \dif t\\
		\leq \int_0^T \int_{\s{D}} e^{-r t}\del{H\del{\epsilon,Q^*(t),\dpd{\hat u}{x}}-H\del{\epsilon,\hat Q^*(t),\dpd{\hat u}{x}} }m(\dif x,t) \dif t\\
		\ \ + \int_0^T \int_{\s{D}}  e^{-r t}\del{H\del{\epsilon,\hat Q^*(t),\dpd{u}{x}}-H\del{\epsilon,Q^*(t),\dpd{u}{x}}}\hat m(\dif x,t)  \dif t
		\\
		+ \intcc{\int_{\s{D}} e^{-r t}(u(x,t)-\hat u(x,t))(m-\hat m)(\dif x,t)}_0^T.
	\end{multline}
	Since $u,\hat u$ are bounded and $\int_{\s{D}} m_i(\dif x,T) \leq 1$ for all $T$, it follows that
	\begin{equation*}
		\lim_{T \to \infty}\int_{\s{D}} e^{-rT}(u(T,x)-\hat u(T,x))(m-\hat m)(\dif x,T) = 0.
	\end{equation*}
	We can rewrite the remaining terms on the right-hand side using the fundamental theorem of calculus.
	Thus \eqref{eq:energy differences2} becomes, after letting $T \to \infty$,
	\begin{equation} \label{eq:energy differences3}
		\frac{1}{C_H}\int_0^\infty \int_{\s{D}}  e^{-r t} \abs{\dpd{u}{x}-\dpd{\hat u}{x}}^2\del{m(\dif x,t) + \hat m(\dif x,t)} \dif t \\
		\leq  I_0 + I_1 + I_2,
	\end{equation}
	where $I_0 := \abs{\int_{\s{D}} (u(0,x)-\hat u(0,x))(m-\hat m)(\dif x,0)}$,
	\begin{equation*}
		\begin{aligned}
			&I_1 := \int_0^1\int_0^\infty \int_{\s{D}} e^{-r t} \dpd{H}{Q}\del{\epsilon,Q_{s}^*(t),\dpd{\hat u}{x}}\del{Q^*(t)-\hat Q^*(t)}(m-\hat m)(\dif x,t)\dif t\dif s, \quad \text{and}\\
			&I_2 := \int_0^1 \int_0^1\int_0^\infty \int_{\s{D}} e^{-r t} \dmpd{H}{Q}{a}\del{\epsilon,Q_{s}^*(t),\dpd{u_{\tilde s}}{x}}\del{\dpd{\hat u}{x}-\dpd{u}{x}}\del{Q^*(t)-\hat Q^*(t)}m(\dif x,t)\dif t\dif s \dif \tilde{s},\\
			&\text{where} \quad Q_{s}^*(t) := sQ^*(t) + (1-s)\hat Q^*(t),
			\
			u_s:=s\hat u + (1-s)u.
		\end{aligned}
	\end{equation*}
	By using Corollary \ref{cor:dHdQ} and \eqref{eq:main estimates}, we can estimate
	\begin{equation}
		\begin{aligned}
			&\abs{\dpd{H}{Q}\del{\epsilon,Q_{s}^*(t),\dpd{\hat u}{x}}} \leq \del{P(0)+1}\epsilon,
			\\
			&\abs{\dmpd{H}{Q}{a}\del{\epsilon,Q_{s}^*(t),\dpd{u_{\tilde s}}{x}}\del{\dpd{\hat u}{x}-\dpd{u}{x}}} \leq \bar P \epsilon, \quad \forall s,\tilde s \in [0,1],
		\end{aligned}
	\end{equation}
	where $\bar P := \max\cbr{\abs{\frac{\bar \rho - 1}{\bar \rho - 2}},1}$ is defined in Corollary \ref{cor:dHdQ}.
	Thus
	\begin{equation} \label{eq:epsilon energy estimates}
		\begin{aligned}
			\abs{I_1} &\leq \del{P(0)+1}\epsilon\int_0^\infty  e^{-r t} \abs{Q^*(t)-\hat Q^*(t)}\enVert{m(t)-\hat m(t)}_{-n}\dif t,\\
			\abs{I_2} &\leq \bar P \epsilon\int_0^\infty \int_{\s{D}} e^{-r t} \abs{\dpd{\hat u}{x}-\dpd{u}{x}}\abs{Q^*(t)-\hat Q^*(t)}m(\dif x,t)\dif t.
		\end{aligned}
	\end{equation}
	Recalling the definitions $w = \hat u - u, \mu = \hat m - m,$ and $\s{Q} = \hat Q - Q$, using the Cauchy-Schwartz inequality and the fact that $m$ is a sub-probability measure, we deduce the following from \eqref{eq:epsilon energy estimates}:
	\begin{equation} \label{eq:epsilon energy estimates1}
		\begin{aligned}
			\abs{I_1} &\leq \del{P(0)+1}\epsilon\del{\int_0^\infty  \abs{\s{Q}_{r/2}(t)}^2\dif t}^{1/2}
			\del{\int_0^\infty  \enVert{\mu_{r/2}(t)}_{-n}^2\dif t}^{1/2},\\
			\abs{I_2} &\leq \bar P \epsilon\del{\int_0^\infty   \abs{\s{Q}_{r/2}(t)}^2 \dif t}^{1/2}
			\del{\int_0^\infty E_{r/2}(t)\dif t}^{1/2}.
		\end{aligned}
	\end{equation}
	We now apply Corollary \ref{cor:Xn Xnstar est} and Assumption \ref{as:Vi}; here we can assume $\nu = 0$ and $f = 0$.
	Thus \eqref{eq:epsilon energy estimates1} implies
	\begin{equation} \label{eq:epsilon energy estimates2}
		\begin{aligned}
			\abs{I_1} &\leq 2\hat D_n(r)C_H\del{P(0)+1}\epsilon \del{\enVert{\mu_0}_{-n}^2 + \int_0^\infty  E_{r/2}(t)\dif t},\\
			\abs{I_2} &\leq 2\hat D_n(r)\bar P \epsilon \del{\enVert{\mu_0}_{-n}^2 + \int_0^\infty  E_{r/2}(t)\dif t}.
		\end{aligned}
	\end{equation}
	Plugging \eqref{eq:epsilon energy estimates2} into \eqref{eq:energy differences3}, we deduce
	\begin{equation} \label{eq:energy differences}
		\int_0^T \int_{\s{D}}  e^{-r t} \abs{\dpd{u}{x}-\dpd{\hat u}{x}}^2\del{m(\dif x,t) + \hat m(\dif x,t)} \dif t 
		\leq C_H\hat C \epsilon \del{\enVert{\mu_0}_{-n}^2 + \int_0^\infty  E_{r/2}(t)\dif t}
		+ C_H I_0,
	\end{equation}
	where $\hat C = 2\hat D_n(r)\del{C_H\del{P(0)+1} + \bar P}$.
	Equation \eqref{eq:epsilon small differences} can be written
	\begin{equation}
		2C_H \hat C\epsilon \leq 1.
	\end{equation}
	Since the left-hand side of \eqref{eq:energy differences} dominates $\int_0^\infty  E_{r/2}(t)\dif t$, we use \eqref{eq:epsilon small differences} and rearrange to deduce \eqref{eq:energy differences final}.
	
	\nextstep \textit{For a linear demand schedule.}
	Now we consider the case where $P(q) = 1-q$ and $\epsilon < 2$.
	In this case the same series of computations (cf.~the proof of Theorem \ref{thm:uniqueness lin demand}, see also Equation \eqref{eq:lin demand ibp} below) now leads to
	\begin{multline} \label{eq:lin energy differences1}
		\frac{1}{4}\int_0^{T'} \int_{\s{D}} \del{\dpd{w_{r/2}}{x} + \epsilon \s{Q}_{r/2}}^2 \dif(\hat m + m)(t) \dif t
		+ \epsilon\int_0^{T'} \s{Q}_{r/2}(t)^2 \dif t\\
		\leq e^{-rT'}\enVert{w(\cdot,T')}_0\enVert{\mu(T')}_{-0} + \enVert{w(\cdot,0)}_n \enVert{\mu_0}_{-n}.
	\end{multline}
	\begin{comment}
	Since $\epsilon \leq 2$ we have
	\begin{multline}
		\abs{\frac{1}{2}\int_0^{T'} \int_{\s{D}} \dpd{w_{r/2}}{x} \epsilon\s{Q}_{r/2}\dif(\hat m + m)(t) \dif t}\\
		\leq \frac{1}{8}\int_0^{T'} \int_{\s{D}} \del{\dpd{w_{r/2}}{x}}^2\dif(\hat m + m)(t) \dif t
		+ \frac{\epsilon^2}{2} \int_0^{T'} \s{Q}_{r/2}(t)^2 \dif t\\
		\leq \frac{1}{8}\int_0^{T'} \int_{\s{D}} \del{\dpd{w_{r/2}}{x}}^2\dif(\hat m + m)(t) \dif t
		+ \epsilon \int_0^{T'} \s{Q}_{r/2}(t)^2 \dif t,
	\end{multline}
	\end{comment}
	Let $T' \to T$, rearrange the square term in \eqref{eq:lin energy differences1} and perform standard estimates to deduce
	\begin{equation} \label{eq:lin energy differences2}
		\int_0^{T} \int_{\s{D}} \del{\dpd{w_{r/2}}{x}}^2 \dif(\hat m + m)(t) \dif t
		\leq 8\enVert{w(\cdot,0)}_n \enVert{\mu_0}_{-n},
	\end{equation}
	which is the same as \eqref{eq:energy differences linear}.
	
\end{proof}

\begin{lemma}[Energy estimates, all other cases]
	\label{lem:energy estimates}
	Let $(w,\mu)$ be a solution of the system \eqref{eq:linearized}, and assume that $V_1,\ldots,V_5,f,\nu$ satisfy  \eqref{eq:V_i for derivative}.
	\begin{enumerate}
		\item Assume that $\epsilon$ is sufficiently small, namely
		\begin{equation} \label{eq:epsilon small linearized}
		\epsilon 4\hat D_n(r)^2\del{\tilde D_n(r) + \bar P}C_H^2 \leq (4C_H)^{-1}.
	\end{equation}
	Then
	\begin{equation} \label{eq:energy linearized}
		\int_0^{T} E_{r/2}(t) \dif t \leq 4C_H\enVert{w(\cdot,0)}_n \enVert{\mu_0}_{-n} 
		+ 4\enVert{\mu_0}_{-n}^2
		+ \hat C\int_0^{T} \del{\enVert{f_{r/2}(\cdot,t)}_n^2 + \enVert{\nu_{r/2}(t)}_{-n}^2} \dif t,
	\end{equation}
	where $\hat C = 4C_H^2\del{\hat D_n(r)^2 + C_H^2} + 1$.
	\item Assume instead that the demand schedule $P$ is linear, i.e.~$P(q) = 1-q$, and that $\epsilon < 2$.
	Then
	\begin{equation} \label{eq:energy linearized lin}
		\int_0^T E_{r/2}(s)\dif s \leq 16 \enVert{w(\cdot,0)}_n \enVert{\mu_0}_{-n}	+ \enVert{\mu_0}_{-n}^2	
		\\
		+ \hat C\int_0^T \del{ \enVert{\nu_{r/2}(t)}_{-n}^2 + \enVert{f_{r/2}(\cdot,t)}_n^2} \dif t,
	\end{equation}
	where $\hat C = \del{32\max\cbr{\hat D_n(r),C_H}^2 + 17}$.
	\end{enumerate}
	
\end{lemma}

\begin{proof}
	Note that the case when $(w,\mu)$ is a difference of two solutions to System \eqref{eq:mfg infty}, so that \eqref{eq:V_i for differences} holds with $f = \nu = 0$, is already proved in Section \ref{sec:uniqueness}.
	
\firststep	\textit{For a small parameter $\epsilon$.}
In this first step, we make no further assumptions on the demand schedule $P$ but instead assume condition \eqref{eq:epsilon small linearized} holds.
Note that when 

Differentiate $e^{-r t} \int_{\s{D}} w\mu$ with respect to $t$ and integrate by parts to get
\begin{multline} \label{eq:energy linearized1}
\dod{}{t}\del{\int_{\s{D}} w_{r/2}\mu_{r/2}} =
\int_{\s{D}} f_{r/2}\mu_{r/2} - \s{Q}_{r/2}(t)\int_{\s{D}} V_2\mu_{r/2}\\
+ \int_{\s{D}} V_4 \abs{\dpd{w_{r/2}}{x}}^2 m + \s{Q}_{r/2}(t) \int_{\s{D}} V_5 \dpd{w_{r/2}}{x} m
+ \int_{\s{D}} \dpd{w_{r/2}}{x} \nu_{r/2}.
\end{multline}
Let $T' \in (0,T)$ and integrate \eqref{eq:energy linearized1} from $0$ to $T'$.
Recalling that $V_4 \geq C_H^{-1}$ from Assumption \ref{as:Vi}, we get
\begin{multline} \label{eq:energy linearized2}
	C_H^{-1} \int_0^{T'} E_{r/2}(t) \dif t \leq \left. \ip{w_{r/2}(\cdot,t)}{\mu_{r/2}(t)}  \right|_0^{T'} 
	+ \int_0^{T'} \enVert{f_{r/2}(\cdot,t)}_n \enVert{\mu_{r/2}(t)}_{-n} \dif t \\
	+ \int_0^{T'} \abs{\s{Q}_{r/2}(t)} \del{\enVert{V_2}_n \enVert{\mu_{r/2}(t)}_{-n} +  \enVert{V_5}_0 E_{r/2}(t)^{1/2}} \dif t
	+ \int_0^{T'} \enVert{\dpd{w_{r/2}}{x}(\cdot,t)}_n \enVert{\nu_{r/2}(t)}_{-n} \dif t.
\end{multline}
Then let $T' \to T$ and recall that by assumption \eqref{eq:w(T) T =infty}, $\lim_{t \to T}\ip{w_{r/2}(\cdot,t)}{\mu_{r/2}(t)} = 0$.
Thus,
\begin{multline} \label{eq:energy linearized3}
	C_H^{-1} \int_0^{T} E_{r/2}(t) \dif t \leq \enVert{w(\cdot,0)}_n \enVert{\mu_0}_{-n}
	+ \int_0^{T} \enVert{f_{r/2}(\cdot,t)}_n \enVert{\mu_{r/2}(t)}_{-n} \dif t \\
	+ \int_0^{T} \abs{\s{Q}_{r/2}(t)} \del{\enVert{V_2}_n \enVert{\mu_{r/2}(t)}_{-n} +  \enVert{V_5}_0 E_{r/2}(t)^{1/2}} \dif t
	+ \int_0^{T} \enVert{\dpd{w_{r/2}}{x}(\cdot,t)}_n \enVert{\nu_{r/2}(t)}_{-n} \dif t.
\end{multline}
Now using Corollary \ref{cor:Xn Xnstar est}, recalling $\enVert{V_4}_0 \leq C_H$ (Assumption \ref{as:Vi}),
\begin{comment}
we have
\begin{multline*}
	\del{\hat D_n(r)^2 + C_H^2}^{-1}\int_0^T \del{\enVert{\mu_{r/2}(t)}_{-n}^2 + \enVert{\dpd{w_{r/2}}{x}(\cdot,t)}_n^2}\dif t\\
	\leq \enVert{\mu_0}_{-n}^2 +  \int_{0}^{T} E_{r/2}(t)\dif t
	+ \int_{0}^{T}  \del{\enVert{\nu_{r/2}(t)}_{-n}^2
		+ \enVert{f_{r/2}(\cdot,t)}_n^2}\dif t
\end{multline*}
On the other hand,
\begin{multline*}
	\int_0^{T} \del{\enVert{f_{r/2}(\cdot,t)}_n \enVert{\mu_{r/2}(t)}_{-n} + \enVert{\dpd{w_{r/2}}{x}(\cdot,t)}_n \enVert{\nu_{r/2}(t)}_{-n}} \dif t\\
	\leq (2C_H)^{-1}\del{\hat D_n(r)^2 + C_H^2}^{-1}\int_0^T \del{\enVert{\mu_{r/2}(t)}_{-n}^2 + \enVert{\dpd{w_{r/2}}{x}(\cdot,t)}_n^2}\dif t\\
	+ C_H\del{\hat D_n(r)^2 + C_H^2}\int_0^{T} \del{\enVert{f_{r/2}(\cdot,t)}_n^2 + \enVert{\nu_{r/2}(t)}_{-n}^2} \dif t.
\end{multline*}
Thus,
\end{comment}
we derive
\begin{multline*}
	\int_0^{T} \del{\enVert{f_{r/2}(\cdot,t)}_n \enVert{\mu_{r/2}(t)}_{-n} + \enVert{\dpd{w_{r/2}}{x}(\cdot,t)}_n \enVert{\nu_{r/2}(t)}_{-n}} \dif t\\
	\leq (2C_H)^{-1}\int_{0}^{T} E_{r/2}(t)\dif t
	+ (2C_H)^{-1}\enVert{\mu_0}_{-n}^2 
	+ C_1\int_0^{T} \del{\enVert{f_{r/2}(\cdot,t)}_n^2 + \enVert{\nu_{r/2}(t)}_{-n}^2} \dif t,
\end{multline*}
where $C_1 := C_H\del{\hat D_n(r)^2 + C_H^2}$.
Thus \eqref{eq:energy linearized3} yields
\begin{multline} \label{eq:energy linearized4}
	(2C_H)^{-1} \int_0^{T} E_{r/2}(t) \dif t \leq \enVert{w(\cdot,0)}_n \enVert{\mu_0}_{-n} 
	+ C_1\int_0^{T} \del{\enVert{f_{r/2}(\cdot,t)}_n^2 + \enVert{\nu_{r/2}(t)}_{-n}^2} \dif t \\
	+ (2C_H)^{-1}\enVert{\mu_0}_{-n}^2
	+ \int_0^{T} \abs{\s{Q}_{r/2}(t)} \del{\enVert{V_2}_n \enVert{\mu_{r/2}(t)}_{-n} +  \enVert{V_5}_0 E_{r/2}(t)^{1/2}} \dif t.
\end{multline}
Also, again using Corollary \ref{cor:Xn Xnstar est} and also Assumption \ref{as:Vi}, we get
\begin{comment}
\begin{equation*}
	\int_0^T \abs{\s{Q}_{r/2}(t)}^2 \dif t \leq 3\hat D_n(r)^2\del{\enVert{\mu_0}_{-n}^2 + \int_0^T \enVert{\nu_{r/2}(t)}_{-n}^2 \dif t 
		+ \int_0^T E_{r/2}(t) \dif t}
\end{equation*}
\begin{equation*}
	\int_0^T \enVert{\mu_{r/2}(t)}_{-n}^2 \dif t \leq \enVert{\mu_0}_{-n}^2
	+ \int_0^T \enVert{\nu_{r/2}(t)}_{-n}^2\dif t
	+ C_H^2\int_0^T E_{r/2}(t)\dif t
\end{equation*}
\begin{multline*}
	\int_0^{T} \abs{\s{Q}_{r/2}(t)} \del{\enVert{V_2}_n \enVert{\mu_{r/2}(t)}_{-n} +  \enVert{V_5}_0 E_{r/2}(t)^{1/2}} \dif t\\
	\leq \del{\enVert{V_2}_n + \enVert{V_5}_0}\int_0^{T} \abs{\s{Q}_{r/2}(t)}^2  \dif t
	+ \enVert{V_2}_n\int_0^{T} \enVert{\mu_{r/2}(t)}_{-n}^2 \dif t
	+ \enVert{V_5}_0\int_0^{T} E_{r/2}(t) \dif t\\
	\leq 4\hat D_n(r)^2\del{\enVert{V_2}_n + \enVert{V_5}_0}\del{\enVert{\mu_0}_{-n}^2
	+ \int_0^T \enVert{\nu_{r/2}(t)}_{-n}^2\dif t
	+ C_H^2\int_0^T E_{r/2}(t)\dif t}
\end{multline*}
\end{comment}
\begin{multline*}
	\int_0^{T} \abs{\s{Q}_{r/2}(t)} \del{\enVert{V_2}_n \enVert{\mu_{r/2}(t)}_{-n} +  \enVert{V_5}_0 E_{r/2}(t)^{1/2}} \dif t\\
	\leq \epsilon 4\hat D_n(r)^2\del{\tilde D_n(r) + \bar P}\del{\enVert{\mu_0}_{-n}^2
		+ \int_0^T \enVert{\nu_{r/2}(t)}_{-n}^2\dif t
		+ C_H^2\int_0^T E_{r/2}(t)\dif t}
\end{multline*}
where $\bar P = \max\cbr{\abs{\frac{\bar \rho - 1}{2 - \bar \rho}},1}$.
Then by \eqref{eq:epsilon small linearized}, Equation \eqref{eq:energy linearized4} yields
\begin{multline} \label{eq:energy linearized5}
	(4C_H)^{-1} \int_0^{T} E_{r/2}(t) \dif t \leq \enVert{w(\cdot,0)}_n \enVert{\mu_0}_{-n} 
	+ C_1\int_0^{T} \del{\enVert{f_{r/2}(\cdot,t)}_n^2 + \enVert{\nu_{r/2}(t)}_{-n}^2} \dif t \\
	+ C_H^{-1}\enVert{\mu_0}_{-n}^2
	+ (4C_H)^{-1}\int_0^T \enVert{\nu_{r/2}(t)}_{-n}^2\dif t.
\end{multline}
We rearrange \eqref{eq:energy linearized5} to conclude with \eqref{eq:energy linearized}.

\nextstep
\textit{For a linear demand schedule.}
Now we consider the case where $P(q) = 1-q$ and $\epsilon < 2$, so that the system has the form \eqref{eq:mfg lin dem}.
After doing integration by parts and canceling like terms, we get
\begin{multline} \label{eq:lin demand ibp}
	\del{2 + \epsilon(t) \int_{\s{D}} \dif m(t)}^{-1}\int_{\s{D}} \abs{\dpd{w}{x}}^2 \dif m(t)
	+ 2\del{2 + \epsilon(t) \int_{\s{D}} \dif m(t)}^{-1}\del{\int_{\s{D}} q^*(\cdot,t) \dif \mu(t)}^2\\
	= e^{rt}\dod{}{t}\del{e^{-rt}\int_{\s{D}} w(\cdot,t)\dif \mu(t)}
	+ 2\del{2 + \epsilon(t) \int_{\s{D}} \dif m(t)}^{-1}\int_{\s{D}} \dif \nu(t) \int_{\s{D}} q^*(\cdot,t) \dif \mu(t)\\
	+ \del{2 + \epsilon(t) \int_{\s{D}} \dif m(t)}^{-1}\int_{\s{D}} \dif \nu(t)\int_{\s{D}} \dpd{w}{x} \dif m(t)
	- \int_{\s{D}} \dpd{w}{x}\dif \nu(t) 
	+  \int_{\s{D}} f(\cdot,t)\dif \mu(t),
\end{multline}
from which we deduce
\begin{multline} \label{eq:dEdt}
	\frac{1}{2}\del{2 + \epsilon(t) \int_{\s{D}} \dif m(t)}^{-1}\int_{\s{D}} \abs{\dpd{w}{x}}^2 \dif m(t)
	+ \del{2 + \epsilon(t) \int_{\s{D}} \dif m(t)}^{-1}\del{\int_{\s{D}} q^*(\cdot,t) \dif \mu(t)}^2\\
	\leq e^{rt}\dod{}{t}\del{e^{-rt}\int_{\s{D}} w(\cdot,t)\dif \mu(t)}
	+ \frac{\epsilon^2 + 2\epsilon}{2}\del{2 + \epsilon(t) \int_{\s{D}} \dif m(t)}^{-1}\del{\int_{\s{D}} \dif \nu(t) }^2\\
	- \int_{\s{D}} \dpd{w}{x}\dif \nu(t) 
	+  \int_{\s{D}} f(\cdot,t)\dif \mu(t).
\end{multline}
Multiply \eqref{eq:dEdt} by $e^{-r t}$, integrate from $0$ to $T'$ and let $T' \to T$ to get
\begin{multline} \label{eq:energy estimate1}
	\int_0^T E_{r/2}(s)\dif s \leq 2(2+\epsilon) \enVert{w(\cdot,0)}_n \enVert{\mu_0}_{-n}		
	+ \frac{\epsilon(2+\epsilon)^2}{2}\int_0^T \enVert{\nu_{r/2}(s)}_{-n}^2 \dif s\\
	+ 2(2+\epsilon)\int_0^T \del{\enVert{\dpd{w_{r/2}}{x}(\cdot,s)}_n\enVert{\nu_{r/2}(s)}_{-n} + \enVert{f_{r/2}(\cdot,s)}_n \enVert{\mu_{r/2}(s)}_{-n}} \dif s.
\end{multline}
\begin{comment}
\begin{multline*}
	2(2+\epsilon)\int_0^T \del{\enVert{\dpd{w_{r/2}}{x}(\cdot,t)}_n\enVert{\nu_{r/2}(t)}_{-n} + \enVert{f_{r/2}(\cdot,t)}_n \enVert{\mu_{r/2}(t)}_{-n}} \dif t\\
	\leq \int_0^T \del{a^{-1}\enVert{\dpd{w_{r/2}}{x}(\cdot,t)}_n^2 + a(2+\epsilon)^2\enVert{\nu_{r/2}(t)}_{-n}^2 + b(2+\epsilon)^2\enVert{f_{r/2}(\cdot,t)}_n^2 +  b^{-1}\enVert{\mu_{r/2}(t)}_{-n}^2} \dif t\\
	\leq a^{-1}\hat D_n(r)^2\del{\enVert{\mu_0}_{-n}^2 + \int_{0}^{T}  \del{E_{r/2}(t) + \enVert{\nu_{r/2}(t)}_{-n}^2
			+ \enVert{f_{r/2}(\cdot,t)}_n^2}\dif t}\\
		+ b^{-1}C_H^2\del{\enVert{\mu_0}_{-n}^2
	+ \int_0^T \del{ E_{r/2}(t) + \enVert{\nu_{r/2}(t)}_{-n}^2}\dif t}\\
	+ \int_0^T \del{ a(2+\epsilon)^2\enVert{\nu_{r/2}(t)}_{-n}^2 + b(2+\epsilon)^2\enVert{f_{r/2}(\cdot,t)}_n^2} \dif t
\end{multline*}
\end{comment}
\begin{comment}
\begin{multline*}
	2(2+\epsilon)\int_0^T \del{\enVert{\dpd{w_{r/2}}{x}(\cdot,t)}_n\enVert{\nu_{r/2}(t)}_{-n} + \enVert{f_{r/2}(\cdot,t)}_n \enVert{\mu_{r/2}(t)}_{-n}} \dif t\\
	\leq \frac{1}{4}\del{\enVert{\mu_0}_{-n}^2 + \int_{0}^{T}  \del{E_{r/2}(t) + \enVert{\nu_{r/2}(t)}_{-n}^2
			+ \enVert{f_{r/2}(\cdot,t)}_n^2}\dif t}\\
	+ \frac{1}{4}\del{\enVert{\mu_0}_{-n}^2
		+ \int_0^T \del{ E_{r/2}(t) + \enVert{\nu_{r/2}(t)}_{-n}^2}\dif t}\\
	+ \int_0^T \del{ 4\hat D_n(r)^2(2+\epsilon)^2\enVert{\nu_{r/2}(t)}_{-n}^2 + 4C_H^2(2+\epsilon)^2\enVert{f_{r/2}(\cdot,t)}_n^2} \dif t
\end{multline*}
\end{comment}
\begin{multline*}
	4(2+\epsilon)\int_0^T \del{\enVert{\dpd{w_{r/2}}{x}(\cdot,t)}_n\enVert{\nu_{r/2}(t)}_{-n} + \enVert{f_{r/2}(\cdot,t)}_n \enVert{\mu_{r/2}(t)}_{-n}} \dif t\\
	\leq \enVert{\mu_0}_{-n}^2 + \int_{0}^{T}  E_{r/2}(t)\dif t\\
	+ \del{8\max\cbr{\hat D_n(r),C_H}^2(2+\epsilon)^2 + 1}\int_0^T \del{ \enVert{\nu_{r/2}(t)}_{-n}^2 + \enVert{f_{r/2}(\cdot,t)}_n^2} \dif t
\end{multline*}
Using Corollary \ref{cor:Xn Xnstar est} and rearranging \eqref{eq:energy estimate1}, we deduce
\begin{multline} \label{eq:energy estimate1-1}
	\int_0^T E_{r/2}(s)\dif s \leq 4(2+\epsilon) \enVert{w(\cdot,0)}_n \enVert{\mu_0}_{-n}	+ \enVert{\mu_0}_{-n}^2	
	+ \epsilon(2+\epsilon)^2\int_0^T \enVert{\nu_{r/2}(s)}_{-n}^2 \dif s\\
	+ \del{8\max\cbr{\hat D_n(r),C_H}^2(2+\epsilon)^2 + 1}\int_0^T \del{ \enVert{\nu_{r/2}(t)}_{-n}^2 + \enVert{f_{r/2}(\cdot,t)}_n^2} \dif t,
\end{multline}
which can be rewritten as \eqref{eq:energy linearized lin}, using $\epsilon < 2$.

\end{proof}

We now introduce the following condition on $\epsilon$:
\begin{assumption}
	\label{as:epsilon small enough}
	We assume either that
	\begin{equation} \label{eq:epsilon small}
		\epsilon \max\cbr{16\hat D_n(r)^2\del{\tilde D_n(r) + \bar P}C_H^3,4C_H \hat D_n(r)\del{C_H\del{P(0)+1} + \bar P}} \leq 1.
	\end{equation}
	where $\bar P = \max\cbr{\frac{\bar \rho - 1}{\bar \rho - 2},1}$, or else $P(q) = 1-q$ and $\epsilon < 2$.
\end{assumption}

\begin{corollary}
	\label{cor:energy estimate2}
	Let $(w,\mu)$ be a solution of \eqref{eq:linearized}, where either \eqref{eq:V_i for derivative} or \eqref{eq:V_i for differences} holds.
	Define
	\begin{align}
		\label{eq:J_n}
		J_n(\rho) &:= \enVert{w(\cdot,0)}_n \enVert{\mu_0}_{-n}		
		+ \enVert{\mu_0}_{-n}^2
		+ \int_0^T \del{\enVert{f_{\rho}(\cdot,s)}_n^2 +  \enVert{\nu_{\rho}(s)}_{-n}^2} \dif s,\\
		\label{eq:K_n}
		K_n(\rho) &:= \enVert{\mu_0}_{-n} + \sup_{0 \leq \tau \leq T} \enVert{\nu_\rho(\tau)}_{-n} + \enVert{\dpd{w_\rho}{x}}_0^{1/2}J_n(\rho)^{1/4}.
	\end{align}
	Let Assumptions \ref{as:r bigger than kappa} and \ref{as:epsilon small enough} hold.
	Then there exists a  constant $C$, depending on the data but not on $T$, such that the following three estimates hold:
	\begin{align}\label{eq:energy estimate2}
		\int_0^T E_{r/2}(s)\dif s &\leq C J_n(r/2),\\
		\label{eq:mu-n int2}
		\int_0^T \enVert{\mu_{r/2}(t)}_{-n}^2 \dif t &\leq C J_n(r/2),\\
		\label{eq:int dw_rdx2}
		\int_0^T \enVert{\dpd{w_{r/2}}{x}(\cdot,t)}_{n}^2 \dif t &\leq C J_n(r/2),\\
		\label{eq:int sQ 2}
		\int_0^T \abs{\s{Q}_{r/2}(t)}^2 \dif t &\leq CJ_n(r/2),\\
		\label{eq:mu_n-0 2}
		\sup_{t \in [0,T]} \enVert{\mu_n(t)}_{-n} &\leq C K_n(r/2).
	\end{align}
\end{corollary}
\begin{proof}
	By Assumption \ref{as:r bigger than kappa}, taking \eqref{eq:w(T) T =infty} into account, we can apply Lemmas \ref{lem:mu-n int} and \ref{lem:w_rho int} with $\rho = r/2$.
	Apply Lemmas \ref{lem:energy differences} and \ref{lem:energy estimates}, we deduce \eqref{eq:energy estimate2}.
	Then Equations \eqref{eq:mu-n int2}, \eqref{eq:int dw_rdx2}, \eqref{eq:int sQ 2}, and \eqref{eq:mu_n-0 2} follow from applying Lemma \ref{lem:mu-n int}, Lemma  \ref{lem:w_rho int}, Corollary \ref{cor:int sQ estimate}, and Lemma \ref{lem:mu-n0}, respectively, using Equation \eqref{eq:energy estimate2}.
\end{proof}

\subsection{H\"older estimates} \label{sec:holder estimates}

Recall that $Y_{1 + \alpha} := \s{C}^{1+\alpha}_\diamond(\s{D})$ is the space of all $\phi \in \s{C}^{1+\alpha}(\overline{\s{D}})$ with the compatibility condition $\phi(0) = 0$.
Set $\psi(x) = 1 - e^{-x}$.
For $n \geq 2$ we will define $Y_{n+\alpha}$ to be the space of all $\phi \in \s{C}^{1+\alpha}_\diamond(\s{D})$ such that $\psi^{j-1} \phi \in \s{C}^{j+\alpha}_\diamond(\s{D})$ for $j = 2,\ldots,n$, with norm given by
\begin{equation} 
	\enVert{\phi}_{Y_{n+\alpha}} = \sum_{j=1}^n \enVert{\psi^{j-1}\phi}_{\s{C}^{j+\alpha}}.
\end{equation}
This defines a Banach space.
The following two lemmas provide estimates on solutions to parabolic equations in the spaces $Y_{n+\alpha}$ for $n = 1,2,3$.
\begin{lemma} \label{lem:uniform t Holder ux}
	Let $u$ be a the solution of 
	\begin{equation}
		\label{eq:linear hj}
		\dpd{u}{t} + \lambda u - \frac{\sigma^2}{2}\dpd[2]{u}{x} + V(x,t)\dpd{u}{x} = F, \quad u(0,t) = 0, \quad u(x,0) = u_0(x)
	\end{equation}
	where $\lambda$ is any positive constant, $F$ is a bounded continuous function, and $u_0 \in \s{C}^{1+\alpha}_\diamond(\s{D})$ (i.e.~$u_0 \in \s{C}^{1+\alpha}(\overline{\s{D}})$ with $u_0(0) = 0$).
	Then
	\begin{equation} \label{eq:uniform t Holder ux}
		\enVert{u}_{\s{C}^{\alpha,\alpha/2}\del{\overline{\s{D}}  \times \intcc{0,T}}} + \enVert{\dpd{u}{x}}_{\s{C}^{\alpha,\alpha/2}\del{\overline{\s{D}} \times \intcc{0,T}}} \leq C\del{\enVert{V}_0,\alpha,\lambda}\del{\enVert{F}_0 + \enVert{u_0}_{\s{C}^{1+\alpha}}},
	\end{equation}
	where $C\del{\enVert{V}_0,\alpha,\lambda}$ is independent of $T$.
	\begin{comment}
	\begin{equation} 
	\begin{aligned}
	\sup_t \del{\enVert{u(\cdot,t)}_{\s{C}^{\alpha}} +  \enVert{\dpd{u}{x}(\cdot,t)}_{\s{C}^{\alpha}}} &\leq C\del{\enVert{V}_0,\alpha,\lambda}\del{\enVert{F}_0 + \enVert{u_0}_{\s{C}^{1+\alpha}}},\\
	\sup_{t_1\neq t_2} \frac{\enVert{u(\cdot,t_1)-u(\cdot,t_2)}_{\s{C}^{\alpha}} +  \enVert{\dpd{u}{x}(\cdot,t_1)-\dpd{u}{x}(\cdot,t_2)}_{\s{C}^{\alpha}}}{\abs{t_1-t_2}} &\leq C\del{\enVert{V}_0,\alpha,\lambda}\del{\enVert{F}_0 + \enVert{u_0}_{\s{C}^{1+\alpha}}},\\
	\end{aligned}.
	\end{equation}
	\end{comment}
\end{lemma}

\begin{proof}
	See \cite[Lemma 2.7]{graber2022parameter}.
\end{proof}

\begin{lemma} \label{lem:3rd order holder}
	Let $u$ be a solution of \eqref{eq:linear hj}, in which $F, \pd{F}{x}, V, \pd{V}{x} \in \s{C}^{\alpha,\alpha/2}(\overline{\s{D}} \times [0,T])$.
	Assume also that $u_0 \in Y_{n+\alpha}$ for $n=2$ or $n = 3$;
	that is, assume $\psi^{j-1} u_0^{(j)} \in \s{C}^{1+\alpha}_\diamond(\s{D})$ for $j = 1,\ldots,n$.
	Then
	\begin{multline}
		\label{eq:2nd order holder psi}
		\enVert{\psi u}_{\s{C}^{2+\alpha,1+\alpha/2}(\overline{\s{D}} \times [0,T])} \\ \leq C\del{\enVert{V}_{\s{C}^{\alpha,\alpha/2}},\lambda,\sigma,\alpha}\del{\enVert{\psi u_0}_{\s{C}^{2+\alpha}}
			+ \enVert{u_0}_{\s{C}^{1+\alpha}}
			+ \enVert{\psi F}_{\s{C}^{\alpha,\alpha/2}}
			+ \enVert{F}_0}, \quad \text{and}
	\end{multline}
	\begin{multline}
		\label{eq:3rd order holder}
		\enVert{\psi^2 \dpd{u}{x}}_{\s{C}^{2+\alpha,1+\alpha/2}(\overline{\s{D}} \times [0,T])} \leq C\del{\enVert{V}_{\s{C}^{\alpha,\alpha/2}},\enVert{\dpd{V}{x}}_{\s{C}^{\alpha,\alpha/2}},\lambda,\sigma,\alpha}\\
		\times \del{\enVert{\psi^2 u_0'}_{\s{C}^{2+\alpha}}  + \enVert{\psi u_0}_{\s{C}^{2+\alpha}} + \enVert{\psi^2 \dpd{F}{x}}_{\s{C}^{\alpha,\alpha/2}} 
			+ \enVert{\psi F}_{\s{C}^{\alpha,\alpha/2}} + \enVert{F}_0}.
	\end{multline}
\end{lemma}

\begin{proof}
	Multiply \eqref{eq:linear hj} by $\psi$ to see that $v(x,t) = \psi(x)u(x,t)$ is the solution to
	\begin{multline}
		\dpd{v}{t} + \lambda v - \frac{\sigma^2}{2}\dpd[2]{v}{x}
		+ V\dpd{v}{x} = \psi F - \sigma^2 \psi' \dpd{u}{x} - \frac{\sigma^2}{2}\psi'' u
		+ V\psi' u,\\
		v(0,t) = 0, \quad v(x,0) = \psi(x)u_0(x).
	\end{multline}
	Note that the compatibility conditions of order 0 and 1 are satisfied.
	Indeed, the condition of order 0 is trivial: $\psi(0)u_0(0) = 0$.
	The condition of order 1 is
	\begin{multline*}
		\lambda \psi(0)u_0(0) - \frac{\sigma^2}{2}\dod[2]{}{x}(\psi u_0)(0)
		+ V(0,0)\dod{}{x}(\psi u_0)(0)\\
		= \psi(0) F(0,0) - \sigma^2 \psi'(0) u_0'(0) - \frac{\sigma^2}{2}\psi''(0) u_0(0)
		+ V(0,0)\psi'(0) u_0(0),
	\end{multline*}
	which can be verified by expanding the derivatives and using the fact that $\psi(0) = 0$.
	Now observe that
	\begin{equation}
		\begin{split}
			\enVert{\sigma^2 \psi' \dpd{u}{x} - \frac{\sigma^2}{2}\psi'' u
				+ V\psi' u}_{\s{C}^{\alpha,\alpha/2}}
			&\leq C\del{\enVert{V}_{\s{C}^{\alpha,\alpha/2}} + 1}\del{\enVert{u}_{\s{C}^{\alpha,\alpha/2}} + \enVert{\dpd{u}{x}}_{\s{C}^{\alpha,\alpha/2}}}\\
			&\leq C\del{\enVert{V}_{\s{C}^{\alpha,\alpha/2}} + 1}\del{\enVert{F}_0 + \enVert{u_0}_{\s{C}^{1+\alpha}}},
		\end{split}
	\end{equation}
	where $C$ depends on $\enVert{V}_0, \alpha,$ and $\lambda$ as in Lemma \ref{lem:uniform t Holder ux}.
	Here we have used the fact that $\enVert{\psi^{(n)}}_0 = 1$ for all $n$.
	From Lemma \ref{lem:C2alpha estimates} we have
	\begin{multline}
		\enVert{v}_{\s{C}^{2+\alpha,1+\alpha/2}(\overline{\s{D}} \times [0,T])}\\
		\leq C\del{\enVert{V}_{\s{C}^{\alpha,\alpha/2}},\lambda,\sigma,\alpha}\del{\enVert{\psi u_0}_{\s{C}^{2+\alpha}} + \enVert{u_0}_{\s{C}^{1+\alpha}} + \enVert{\psi F}_{\s{C}^{\alpha,\alpha/2}} + \enVert{F}_0
			+ \enVert{\dpd{v}{x}}_{\s{C}^{\alpha,\alpha/2}}},
	\end{multline}
	and Equation \eqref{eq:2nd order holder psi} follows from interpolation.

	To derive Equation \eqref{eq:3rd order holder}, take the derivative with respect to $x$ of \eqref{eq:linear hj} and multiply by $\psi^2$.
	Rearrange to see that $w(x,t) = \psi(x)^2\pd{u}{x}(x,t)$ is the (weak) solution to
	\begin{multline*}
		\dpd{w}{t} + \lambda w - \frac{\sigma^2}{2}\dpd[2]{w}{x} = \psi^2 \dpd{F}{x} - \del{\psi^2 \dpd{V}{x} + \sigma^2(\psi')^2 + \sigma^2\psi\psi''}\dpd{u}{x} - \del{\psi^2 V + 2\sigma^2\psi\psi'} \dpd[2]{u}{x}, \\ w(0,t) = 0, \quad w(x,0) = \psi(x)^2u_0'(x).
	\end{multline*}
	Notice that, thanks to the fact that $\psi(0) = 0$, the compatibility conditions of order 0 and 1 are satisfied, by the same reasoning as above.
	We also have, using Lemma \ref{lem:uniform t Holder ux},
	\begin{equation}
		\begin{split}
			&\hspace{-1cm}\enVert{\psi^2 \dpd{F}{x} - \del{\psi^2 \dpd{V}{x} + \sigma^2(\psi')^2 + \sigma^2\psi\psi''}\dpd{u}{x} - \del{\psi^2 V + 2\sigma^2\psi\psi'} \dpd[2]{u}{x}
			}_{\s{C}^{\alpha,\alpha/2}}\\
			&\leq \enVert{\psi^2 \dpd{F}{x}}_{\s{C}^{\alpha,\alpha/2}}
			+ C\enVert{\dpd{u}{x}}_{\s{C}^{\alpha,\alpha/2}}
			+ C\enVert{\psi\dpd[2]{u}{x}}_{\s{C}^{\alpha,\alpha/2}}\\
			&\leq \enVert{\psi^2 \dpd{F}{x}}_{\s{C}^{\alpha,\alpha/2}}
			+ C\enVert{\dpd{u}{x}}_{\s{C}^{\alpha,\alpha/2}}
			+ C\enVert{u}_{\s{C}^{\alpha,\alpha/2}}
			+ C\enVert{\psi u}_{\s{C}^{2+\alpha,1+\alpha/2}}\\
			&\leq \enVert{\psi^2 \dpd{F}{x}}_{\s{C}^{\alpha,\alpha/2}}
			+ C\enVert{F}_0
			+ \enVert{\psi u}_{\s{C}^{2+\alpha,1+\alpha/2}},
		\end{split}
	\end{equation}
	where $C$ depends on $\enVert{V}_{\s{C}^{\alpha,\alpha/2}}$ and $\enVert{\pd{V}{x}}_{\s{C}^{\alpha,\alpha/2}}$.
	By Lemma \ref{lem:C2alpha estimates} and Equation \eqref{eq:2nd order holder psi}, we deduce \eqref{eq:3rd order holder}.
\end{proof}

Lemmas \ref{lem:uniform t Holder ux} and \ref{lem:3rd order holder} have the following consequence in the case $F = 0$:
\begin{corollary}
	\label{cor:holder est}
	Let $u$ be the solution of \eqref{eq:linear hj} where $\lambda$ is any positive constant and where $F = 0$.
	Then
	\begin{equation}
		\sup_{t \geq 0}\enVert{u(\cdot,t)}_{Y_{n+\alpha}}
		+ \sup_{t_1 \neq t_2} \frac{\enVert{u(\cdot,t_1) - u(\cdot,t_2)}_{Y_{n+\alpha}}}{\abs{t_1-t_2}^{\alpha/2}}
		\leq C\enVert{u_0}_{Y_{n+\alpha}},
	\end{equation}
	where $C$ depends on $\alpha,\lambda,\sigma$, and on either $\enVert{V}_0$ (if $n = 1$), $\enVert{V}_{\s{C}^{\alpha,\alpha/2}}$ (if $n = 2$), or $\enVert{V}_{\s{C}^{\alpha,\alpha/2}} + \enVert{\pd{V}{x}}_{\s{C}^{\alpha,\alpha/2}}$ (if $n = 3$).
\end{corollary}

Next we wish to establish estimates on the Fokker-Planck equation in the spaces $Y_{n + \alpha}^*$, denoting the dual of $Y_{n+\alpha}$, with regularity in time as well.
Note that $\enVert{\cdot}_n \leq \enVert{\cdot}_{Y_{n+\alpha}}$ and thus $\enVert{\cdot}_{Y_{n + \alpha}^*} \leq \enVert{\cdot}_{-n}$.

\begin{lemma} \label{lem:fp holder}
	Let $(w,\mu)$ be a solution of \eqref{eq:linearized}.
	Suppose Assumption \ref{as:r bigger than kappa} holds.
	Then
	\begin{equation} \label{eq:fp holder}
		\enVert{\mu_{r/2}}_{\s{C}^{\alpha/2}\del{[0,T]; Y_{n + \alpha}^*}} \leq C(\alpha,r,\sigma)J_n(r/2)^{1/2}, \quad n = 1,2,
	\end{equation}
	where $J_n$ is defined in \eqref{eq:J_n}.
\end{lemma}

\begin{proof}
	\firststep
	Let $\lambda > 0$ be such that $\lambda < r/2$.
	Fix $t_1 > 0$, let $\phi_{t_1} \in Y_{n+\alpha}$ with $\enVert{\phi_{t_1}}_{Y_{n+\alpha}} \leq 1$, and for any $\lambda > 0$ let $\phi^{(\lambda)}$ denote the solution of
	\begin{equation} \label{eq:phi equation}
		-\dpd{\phi}{t} + \lambda \phi - \frac{\sigma^2}{2}\dpd[2]{\phi}{x} - V_3(x,t)\dpd{\phi}{x} = 0, \quad \phi(0,t) = 0, \quad \phi(x,t_1) = \phi_{t_1}(x).
	\end{equation}
	Note that we have the relation
	\begin{equation}
		\phi^{(\lambda_1 + \lambda_2)}(x,t) = e^{\lambda_2(t-t_1)}\phi^{(\lambda_1)}(x,t).
	\end{equation}
	Now $\enVert{q^*}_{\s{C}^{\alpha,\alpha/2}}$ and $\enVert{\pd{q^*}{x}}_{\s{C}^{\alpha,\alpha/2}}$ can be estimated using the norm $\enVert{u}_{\s{C}^{2+\alpha,1+\alpha/2}}$, which in turn is estimated by the a priori estimates in Theorem \ref{thm:exist mfg infty}.
	By Corollary \ref{cor:holder est} we therefore have
	\begin{equation} \label{eq:phi holder}
		\sup_{t \in [0,t_1]}\enVert{\phi^{(\lambda)}(t)}_{Y_{n+\alpha}} + \sup_{t_0 \in [0,t_1)}\frac{\enVert{\phi^{(\lambda)}(t_1) - \phi^{(\lambda)}(t_0)}_{Y_{n+\alpha}}}{(t_1-t_0)^{\alpha/2}}
		\leq C\del{\alpha,\lambda,\sigma}.
	\end{equation}
	For any $t_0 \in [0,t_1]$ we have, using integration by parts,
	\begin{multline} \label{eq:duality}
		\int_{\s{D}} \phi_{t_1}^{({r/2})}(x)\mu_{{r/2}}(x,t_1)\dif x =
		\int_{\s{D}} \phi^{({r/2})}(x,t_0)\mu_{{r/2}}(x,t_0)\dif x \\ -\int_{t_0}^{t_1}\int_{\s{D}}\dpd{\phi^{({r/2})}}{x}\del{\del{V_4(x,t)\dpd{w_{r/2}}{x}(x,t) + V_5(x,t)\s{Q}_{r/2}(t)}m + \nu_{r/2}}\dif x\dif t.
	\end{multline}
	Applying \eqref{eq:phi holder} and Corollary \ref{cor:energy estimate2}, using the identity $\phi^{({r/2})} = e^{(\lambda - {r/2})(t_1-t)}\phi^{(\lambda)}$, we have
	\begin{equation} \label{eq:fp holder 1}
		\begin{split}
			&\abs{\int_{t_0}^{t_1}\int_{\s{D}}\dpd{\phi^{({r/2})}}{x}\del{\del{V_4(x,t)\dpd{w_{r/2}}{x}(x,t) + V_5(x,t)\s{Q}_{r/2}(t)}m + \nu_{r/2}}\dif x\dif t}\\
			&\leq C\del{\alpha,\lambda,\sigma}\int_{t_0}^{t_1} e^{(\lambda - {r/2})(t_1-t)}\del{E_{r/2}(t)^{1/2} + \enVert{\mu_{r/2}(t)}_{-n}  + \enVert{\nu_{r/2}(t)}_{-n}}\dif t\\
			&\leq C\del{\alpha,\lambda,\sigma}\del{\int_{0}^{t_1-t_0} e^{-2({r/2} - \lambda)t} \dif t}^{1/2}\del{\int_{t_0}^{t_1} \del{E_{r/2}(t) + \enVert{\mu_{r/2}(t)}_{-n}^2  + \enVert{\nu_{r/2}(t)}_{-n}^2}\dif t}^{1/2}\\
			&\leq C\del{\alpha,\lambda,\sigma}\min\cbr{({r/2} - \lambda)^{-1/2},(t_1-t_0)^{1/2}}J_n({r/2})^{1/2}.
		\end{split}
	\end{equation}
	Using \eqref{eq:fp holder 1} in \eqref{eq:duality} with $t_0 = 0$ and taking the supremum over all $\phi_{t_1}$ we deduce the bound
	\begin{equation} \label{eq:fp holder 2}
		\enVert{\mu_{r/2}(t_1)}_{Y_{n + \alpha}^*}
		\leq C\del{\alpha,\lambda,\sigma}({r/2} - \lambda)^{-1/2}\del{\enVert{\mu_0}_{Y_{n + \alpha}^*} +  J_n({r/2})^{1/2}} \quad \forall t_1 \geq 0.
	\end{equation}
	On the other hand, subtracting $\int \phi_{t_1}\mu_{r/2}(t_1)$ from both sides of \eqref{eq:duality}, we have
	\begin{multline} \label{eq:duality2}
		\int_{\s{D}} \phi_{t_1}(x)\del{\mu_{{r/2}}(x,t_1) - \mu_{{r/2}}(x,t_0)}\dif x =
		\int_{\s{D}} \del{\phi^{({r/2})}(x,t_0) - \phi^{({r/2})}(x,t_1)}\mu_{r/2}(x,t_0)\dif x \\ -\int_{t_0}^{t_1}\int_{\s{D}}\dpd{\phi^{({r/2})}}{x}\del{\del{V_4(x,t)\dpd{w}{x}(x,t) + V_5(x,t)\s{Q}(t)}m + \nu}\dif x\dif t.
	\end{multline}
	Combining \eqref{eq:phi holder}, \eqref{eq:fp holder 1}, and \eqref{eq:fp holder 2} to estimate the right-hand side of \eqref{eq:duality2}, we deduce that
	\begin{equation} \label{eq:fp holder 3}
		\enVert{\mu_{r/2}(t_1) - \mu_{r/2}(t_0)}_{Y_{n + \alpha}^*} \leq C\del{\alpha,\lambda,\sigma}\del{\enVert{\mu_0}_{Y_{n + \alpha}^*} + J_n({r/2})^{1/2}}(t_1-t_0)^{\alpha/2}.
	\end{equation}
	It suffices to take $\lambda = r/4$.
	Then recalling that $\enVert{\mu_0}_{Y_{n + \alpha}^*} \leq \enVert{\mu_{0}}_{-n} \leq J_n({r/2})^{1/2}$, we see that \eqref{eq:fp holder 2} and \eqref{eq:fp holder 3} imply \eqref{eq:fp holder}.
\end{proof}

\subsection{H\"older regularity of the mass function}

Let $(w,\mu)$ solve \eqref{eq:linearized}.
Our goal is to prove the H\"older regularity of the following
functional:
\begin{equation}
	\eta_{\rho}(t) = e^{-\rho t}\ip{1}{\mu_{\rho}(t)} = e^{-\rho t}\int_{\s{D}} \mu(x,t)\dif x.
\end{equation}
This will allow us to estimate $\s{Q}_\rho$ in a H\"older space.

We introduce the space $\sr{M}_{\alpha}^{-n}$, in analogy to the space $\sr{M}_\alpha$ defined in Section \ref{sec:prelim}.
For any $\mu \in X_n^*$ define the \emph{mass function}
\begin{equation}
	\label{eq:eta function -n}
	\eta^h[\mu](t) := \ip{\mu}{1} - \int_{\s{D}} \ip{\del{S(x-\cdot,t) - S(x+\cdot,t)}}{\mu} \dif x,
\end{equation}
cf.~\eqref{eq:eta function}.
By Proposition \ref{pr:potentials}, we deduce that
\begin{equation}
	\enVert{\int_{\s{D}} \del{S(x-\cdot,t) - S(x+\cdot,t)}\dif x}_n \leq 2M_n \enVert{1}_n = 2M_n,
\end{equation}
and thus we can write \eqref{eq:eta function -n} as
\begin{equation}
	\label{eq:eta function -n alt}
	\eta^h[\mu](t) := \ip{\mu}{1} -  \ip{\int_{\s{D}}\del{S(x-\cdot,t) - S(x+\cdot,t)}\dif x}{\mu},
\end{equation}
from which we also deduce
\begin{equation}
	\abs{\eta^h[\mu](t)} \leq C_n\enVert{\mu}_{-n} \quad \forall t \geq 0.
\end{equation}
Now we define $\sr{M}_{\alpha}^{-n}$ to be the set of all $\mu \in X_n^*$ such that $\eta^h[\mu]$ is $\alpha$-H\"older continuous.
It is a Banach space endowed with the norm
\begin{equation}
	\enVert{\mu}_{\sr{M}_{\alpha}^{-n}} = \enVert{\mu}_{-n} + \enVert{\eta^h[\mu]}_{\s{C}^\alpha\del{\intco{0,\infty}}}.
\end{equation}

\begin{lemma} \label{lem:eta_rho holder}
	Let $(w,\mu)$ solve \eqref{eq:linearized},
	and suppose Assumption \ref{as:r bigger than kappa} holds.
	Assume that $\alpha \leq 2/5$.
	There exists a constant $C$ depending only on the data but not on $T$ such that
	\begin{equation}
		\label{eq:eta_rho holder}
		\enVert{\eta_{r/2}}_{\s{C}^{\alpha/2}([0,T])} \leq C\del{\enVert{\mu_0}_{\sr{M}_{\alpha/2}^{-n}} + \tilde K_n({r/2})},
	\end{equation} 
	where
	\begin{equation} \label{eq:tilde K_n}
		\tilde K_n({r/2}) = K_n({r/2}) + \enVert{\dpd{w_{{r/2}}}{x}}_0^{2/3}J_n({r/2})^{1/6}.
	\end{equation}	
	and where $J_n({r/2})$ and $K_n({r/2})$ are defined in \eqref{eq:J_n} and \eqref{eq:K_n}, respectively.
\end{lemma}
\begin{proof}
	Observe that $\enVert{1}_n = 1$ for all $n$ and $\enVert{\xi}_n \leq 2$ for $n = 0,1,2$, so we have the bounds
	\begin{equation}
		\abs{\eta_{{r/2}}(t)} \leq \enVert{\mu_{{r/2}}(t)}_{-n}, \quad
		\abs{\zeta_{{r/2}}(t)} \leq \enVert{\mu_{{r/2}}(t)}_{-n}.
	\end{equation}
	It remains to prove estimates on the H\"older seminorms.
	
	\firststep
	By Duhamel's Principle, we can write
	\begin{equation} \label{eq:mu=sum I_j}
		\mu(x,t) = I_1(x,t) + I_2(x,t) + I_3(x,t) + I_4(x,t)
	\end{equation}
	where
	\begin{gather*}
		I_1(x,t) = -\int_{0}^t \int_0^\infty \del{S(x-y,t-s)-S(x+y,t-s)}(V_3\mu)_y(y,s)\dif y \dif s,\\
		I_2(x,t) = -\int_{0}^t \int_0^\infty \del{S(x-y,t-s)-S(x+y,t-s)}(bm)_y(y,s)\dif y \dif s,\\
		I_3(x,t) =- \int_{0}^t \int_0^\infty \del{S(x-y,t-s)-S(x+y,t-s)}\nu_y(y,s)\dif y \dif s,\\
		I_4(x,t) = \int_0^\infty \del{S(x-y,t)-S(x+y,t)}\mu_0(y)\dif y,\\
		b(x,t) = V_4(x,t)\dpd{w}{x}(x,t) + V_5(x,t)\s{Q}(t).
	\end{gather*}
	Using integration by parts, we deduce
	\begin{equation*}
		\eta_{{r/2}}(t) = \sum_{j=1}^4 \eta^j(t)
	\end{equation*}
	where
	\begin{equation}
		\begin{aligned}
			\eta^1(t) &= 2\int_0^t \int_0^\infty e^{-{r/2} (t-s)}S(x,t-s)V_3(x,s)\mu_{{r/2}}(x,s)\dif x\dif s,\\
			\eta^2(t) &= 2\int_0^t \int_0^\infty e^{-{r/2} (t-s)}S(x,t-s)b_{{r/2}}(x,s) m(x,s)\dif x\dif s,\\
			\eta^3(t) &= 2\int_0^t \int_0^\infty S(x,t-s)e^{-{r/2} s}\nu_{{r/2}}(x,s)\dif x\dif s,\\
			\eta^4(t) &= \frac{2}{\sqrt{\pi}}e^{-{r/2} t}\int_0^\infty \int_{\del{2\sigma^2 t}^{1/2}x}^\infty e^{-x^2}\mu_0(y)\dif y\dif x
		\end{aligned}
	\end{equation}
	where we follow the usual convention defining $b_{r/2}(x,t) = e^{-{r/2} t}b(x,t)$.
	We use much the same arguments as in Lemma \ref{lem:fp eta holder} to establish H\"older estimates.
	
	\nextstep
	For the first term, we write
	\begin{multline}
		\eta^1(t_1) - \eta^1(t_0) 
		= 2\int_{t_0}^{t_1} \int_0^\infty e^{-{r/2} (t-s)}S(x,t_1-s)V_3(x,s)\mu_{{r/2}}(x,s)\dif x\dif s\\
		+ 2\int_{0}^{t_0} \int_0^\infty \int_{t_0}^{t_1}\od{}{t}\sbr{e^{-{r/2} (t-s)}S(x,t-s)}V_3(x,s)\mu_{{r/2}}(x,s)\dif t\dif x\dif s.
	\end{multline}
	Use Corollary \ref{cor:energy estimate2} and Assumption \ref{as:Vi} to get
	\begin{multline}
		\abs{\eta^1(t_1) - \eta^1(t_0)}
		\leq C(n,r)K_n({r/2})\int_{t_0}^{t_1}  e^{-{r/2} (t-s)}\enVert{S(\cdot,t_1-s)}_n\dif s\\
		+ C(n,r)K_n({r/2})\int_{0}^{t_0} \int_{t_0}^{t_1}\enVert{\od{}{t}\sbr{e^{-{r/2} (t-s)}S(x,t-s)}}_n \dif t\dif s.
	\end{multline}
	Use Lemma \ref{lem:m_n} to get
	\begin{equation}
		\int_{t_0}^{t_1}  e^{-{r/2} (t-s)}\enVert{S(\cdot,t_1-s)}_n\dif s
		\leq C(n)\int_{t_0}^{t_1} (t_1-s)^{-1/2}\dif s = C(n)(t_1-t_0)^{1/2}.
	\end{equation}
	On the other hand, from Lemma \ref{lem:m_n} we also have
	\begin{equation}
		\sup_{x,t} \abs{x^n t^{3/2} \dpd[n+2]{S}{x}(x,t)} < \infty
		\quad
		\Rightarrow \enVert{\dpd[2]{S}{x}(x,t)}_n \leq C(n)t^{-3/2}
	\end{equation}
	for any $n$. We use this to deduce
	\begin{multline} \label{eq:ddt S int}
		\int_{0}^{t_0} \int_{t_0}^{t_1}\enVert{\od{}{t}\sbr{e^{-{r/2} (t-s)}S(x,t-s)}}_n \dif t\dif s\\
		= \int_{0}^{t_0} \int_{t_0}^{t_1}e^{-{r/2} (t-s)}\enVert{\frac{\sigma^2}{2}\dpd[2]{S}{x}(x,t-s) - {r/2} S(x,t-s)}_n \dif t\dif s\\
		\leq C(n,\sigma)\int_0^{t_0} \int_{t_0}^{t_1} (t-s)^{-3/2}\dif t \dif s
		+ C(n){r/2}\int_0^{t_0} \int_{t_0}^{t_1} (t-s)^{-1/2}\dif t \dif s
		\leq C(n,\sigma,{r})(t_1-t_0)^{1/2}
	\end{multline}
	so long as $t_1 - t_0 \leq 1$.
	These estimates combine to give
	\begin{equation} \label{eq:eta_1}
		\abs{\eta^1(t_1) - \eta^1(t_0)}
		\leq C({r},\sigma,n)K_n({r/2})(t_1-t_0)^{1/2} \quad \forall 0 \leq t_0 \leq t_1 \leq t_0 + 1.
	\end{equation}
	By the very same argument, we also have
	\begin{equation} \label{eq:eta_3}
		\abs{\eta^3(t_1) - \eta^3(t_0)}
		\leq C({r},\sigma,n)K_n({r/2})(t_1-t_0)^{1/2} \quad \forall 0 \leq t_0 \leq t_1 \leq t_0 + 1.
	\end{equation}

	\nextstep
	Next we write
	\begin{multline}
		\eta^2(t_1) - \eta^2(t_0) 
		= -2\int_{t_0}^{t_1} \int_0^\infty e^{-{r/2} (t-s)}S(x,t_1-s)b_{{r/2}}(x,s)m(x,s)\dif x\dif s\\
		-2\int_{0}^{t_0} \int_0^\infty \int_{t_0}^{t_1}\od{}{t}\sbr{e^{-{r/2} (t-s)}S(x,t-s)}b_{{r/2}}(x,s)m(x,s)\dif t\dif x\dif s.
	\end{multline}
	Recall that $b = V_4\pd{w}{x} + V_5\s{Q}$, and recall also the formula \eqref{eq:linearized}(iii) for $\s{Q}$.
	Applying Lemma \ref{lem:mu-n0}, we have
	\begin{equation}
		\begin{split}
			\int_0^\infty \abs{b_{{r/2}}(x,s)}m(x,s)\dif x
			&\leq C(n)K_n({r/2}) + C\int_0^\infty \abs{\dpd{w_{{r/2}}}{x}(x,s)}m(x,s)\dif x\\
			&\leq C(n)K_n({r/2}) + C\enVert{\dpd{w_{{r/2}}}{x}}_0^{2/3}E_{r/2}(s)^{1/6}.
		\end{split}
	\end{equation}
	Using the same reasoning as in the previous step, we deduce
	\begin{multline}
		\abs{\eta^2(t_1) - \eta^2(t_0)}
		\leq C({r},\sigma,n)K_n({r/2})(t_1-t_0)^{1/2}
		+ C(n)\enVert{\dpd{w_{{r/2}}}{x}}_0^{2/3}\int_{t_0}^{t_1} (t-s)^{-1/2}E_{r/2}(s)^{1/6}\dif s\\
		+ C(n,\sigma,{r})\enVert{\dpd{w_{{r/2}}}{x}}_0^{2/3}\int_0^{t_0}\int_{t_0}^{t_1} \del{(t-s)^{-3/2} + (t-s)^{-1/2}}E_{r/2}(s)^{1/6}\dif s,
	\end{multline}
	for $0 \leq t_0  \leq t_1 \leq t_0 + 1$.
	(Cf.~Equation \eqref{eq:ddt S int}.)
	By H\"older's inequality, we compute
	\begin{equation}
		\begin{split}
			\int_{t_0}^{t_1}(t-s)^{-1/2}E_{r/2}(s)^{1/6}\dif s
			&\leq C(t_1-t_0)^{1/3}
			\del{\int_{t_0}^{t_1} E_{r/2}(s)\dif s}^{1/6},\\
			\int_0^{t_0}\int_{t_0}^{t_1} \del{(t-s)^{-3/2} + (t-s)^{-1/2}}E_{r/2}(s)^{1/6}\dif s
			&\leq C(t_1-t_0)^{1/5}\del{\int_{t_0}^{t_1} E_{r/2}(s)\dif s}^{1/6}.
		\end{split}
	\end{equation}
	Combining this with Corollary \ref{cor:energy estimate2}, we have
	\begin{equation} \label{eq:eta_2}
		\abs{\eta^2(t_1) - \eta^2(t_0)}
		\leq C({r},\sigma,n)\tilde K_n({r/2})(t_1-t_0)^{1/5},
	\end{equation}
	where $\tilde K_n({r/2})$ is defined in \eqref{eq:tilde K_n}.
	
	\nextstep
	For the last term $\eta^4(t)$, we use the definition of $\sr{M}_\alpha^{-n}$ and the mass function \eqref{eq:eta function -n} to see that
	\begin{equation} \label{eq:eta_4}
		\eta^4(t) = e^{-{r/2} t}\eta^h[\mu_0](t),
	\end{equation}
	and so, because $t \mapsto e^{-{r/2} t}$ is globally Lipschitz with constant ${r/2}$ on the interval $\intco{0,\infty}$, we deduce
	\begin{equation}
		\enVert{\eta^4}_{\s{C}^{\alpha/2}([0,T])} \leq \max\cbr{1,{r/2}}\enVert{\mu_0}_{\sr{M}_{\alpha/2}^{-n}}.
	\end{equation}
	Putting together \eqref{eq:eta_1}, \eqref{eq:eta_2}, \eqref{eq:eta_3}, and \eqref{eq:eta_4}, we deduce \eqref{eq:eta_rho holder}.
	
	\begin{comment}
	\nextstep
	It remains to prove \eqref{eq:zeta_rho holder}.
	Multiply the Fokker-Planck equation \eqref{eq:linearized}(ii) by $\xi(x)$ and integrate by parts to get
	\begin{multline} \label{eq:xi mu_rho}
	\dod{}{t}\int_0^\infty \xi(x)\mu_{r/2}(x,t)\dif x
	= \int_0^\infty \del{\frac{\sigma^2}{2}\xi''(x)
	- \xi'(x)q^*(x,t) - {r/2} \xi(x)}\mu_{r/2}(x,t)\dif x\\
	+ \frac{1}{2}\int_{\s{D}} \xi'(x)\dpd{w_{r/2}}{x}(x,t)m(\dif x,t)
	+ \frac{\epsilon}{2}\s{Q}_{r/2}(t)\int_{\s{D}} \xi'(x)m(\dif x,t)
	+ \int_{\s{D}} \xi'(x)\nu_{r/2}(x,t)\dif x.
	\end{multline}
	Using $\enVert{\xi}_{\s{C}^4} \leq 2$ we deduce
	\begin{equation} \label{eq:xi mu_rho2}
	\abs{\dod{}{t}\int_0^\infty \xi(x)\mu_{r/2}(x,t)\dif x}
	\leq C_n\max\cbr{{r/2},1}\enVert{\mu_{r/2}(t)}_{-n}
	+ E_{r/2}(t)^{1/2}
	+ \epsilon\abs{\s{Q}_{r/2}(t)}
	+ 2\enVert{\nu_{r/2}(t)}_{-n}.
	\end{equation}
	Applying \eqref{eq:abs Q} to \eqref{eq:xi mu_rho2} we get
	\begin{equation} \label{eq:xi mu_rho3}
	\abs{\dod{}{t}\int_0^\infty \xi(x)\mu_{r/2}(x,t)\dif x}
	\leq C_n\max\cbr{{r/2},1}\enVert{\mu_{r/2}(t)}_{-n}
	+ \del{1 + \frac{\epsilon}{2}}E_{r/2}(t)^{1/2}
	+ (2 + \epsilon)\enVert{\nu_{r/2}(t)}_{-n}.
	\end{equation}
	Recalling the definition of $J_n({r/2})$ from \eqref{eq:J_n}, we deduce
	\begin{equation}
	\abs{\int_0^\infty \xi(x)\del{\mu_{r/2}(x,t_1) - \mu_{r/2}(x,t_0)}\dif x}
	\leq C_n\max\cbr{{r/2},1}J_n({r/2})\abs{t_1-t_0}^{1/2} \quad \forall t_0,t_1 \geq 0,
	\end{equation}
	from which \eqref{eq:zeta_rho holder} follows.
	\end{comment}
\end{proof}

\begin{corollary} \label{cor:sQ Holder}
	Let $(w,\mu)$ be a solution of \eqref{eq:linearized} and suppose Assumption \ref{as:r bigger than kappa} holds.
	Assume $\alpha \leq 2/5$.
	Then there exists a constant $C$, depending only on the data but not on $T$, such that, for $n=1,2$,
	\begin{equation} \label{eq:sQ Holder}
		\enVert{\s{Q}_{r/2}}_{\s{C}^{\alpha/2}([0,T])} \leq C\del{\tilde J_n(r/2) + \enVert{\dpd{w_{r/2}}{x}}_{\s{C}^{\alpha,\alpha/2}}},
	\end{equation}
	where
	\begin{equation}
		\label{eq:tilde J_n}
		\tilde J_n(r/2) := \enVert{\ip{\nu_{r/2}}{1}}_{\s{C}^{\alpha/2}([0,T])} + \enVert{\mu_0}_{\sr{M}_{\alpha/2}^{-n}} + J_n({r/2})^{1/2} + \tilde K_n({r/2}),
	\end{equation}
	and $J_n$ and $\tilde K_n$ are defined in \eqref{eq:J_n} and \eqref{eq:tilde K_n}, respectively.
\end{corollary}

\begin{proof}
	Multiplying Equation \eqref{eq:linearized}(iii) by $e^{-r t/2}$, we have
	\begin{equation} \label{eq:sQ(t)}
		\s{Q}_{r/2}(t) = g_1(t)\del{g_2(t) + g_3(t) + g_4(t)}
	\end{equation}
	where
	\begin{equation}
		\begin{split}
			g_1(t) &:= \del{1 + \frac{\epsilon}{2}\int_{\s{D}} \dif m(t)}^{-1},\\
			g_2(t) &:= - \int_{\s{D}} \dif \nu_{r/2}(t),\\
			g_3(t) &:= \int_{\s{D}} q^*(\cdot,t)\dif \mu_{r/2}(t),\\
			g_4(t) &:= - \frac{1}{2}\int_{\s{D}} \dpd{w_{r/2}}{x}(\cdot,t)\dif m(t).
		\end{split}
	\end{equation}
	Using the fact that $m(t)$ is a positive measure-valued process together with the H\"older regularity deduced from Lemma \ref{lem:fp eta holder}, we have
	\begin{equation} \label{eq:g1}
		\enVert{g_1}_{\s{C}^{\alpha/2}([0,T])} \leq C.
	\end{equation}
	On the other hand,
	\begin{equation} \label{eq:g2}
		\enVert{g_2}_{\s{C}^{\alpha/2}([0,T])} =  \enVert{\ip{\nu_{r/2}}{1}}_{\s{C}^{\alpha/2}([0,T])},
	\end{equation}
	which is taken as given.
	Next, we analyze $g_3$.
	Set
	\begin{equation*}
		\tilde q(x,t) = q^*(x,t) - q^*(0,t),
	\end{equation*}
	so that
	\begin{equation} \label{eq:g3}
		g_3(t) = \int_{\s{D}} \tilde q(\cdot,t)\dif \mu_{r/2}(t) + q^*(0,t)\eta_{r/2}(t) =: g_{3,1}(t) + g_{3,2}(t).
	\end{equation}
	Observe that, since $\tilde q(0,t) = 0$ by construction, we have $\tilde q \in \s{C}^{\alpha/2}([0,T];Y_{n+\alpha})$, where by computing the derivatives of $q^*$ we deduce
	\begin{equation}
		\enVert{\tilde q}_{\s{C}^{\alpha/2}([0,T];Y_{n+\alpha})}
		\leq \enVert{u}_{\s{C}^{2+\alpha,1+\alpha/2}} + \enVert{\psi \dpd{u}{x}}_{\s{C}^{2+\alpha,1+\alpha/2}} \leq C, \quad n = 1,2.
	\end{equation}
	Therefore, using Lemma \ref{lem:fp holder}, we get
	\begin{equation} \label{eq:g31}
		\enVert{g_{3,1}}_{\s{C}^{\alpha/2}([0,T])}
		\leq C\enVert{\mu_{r/2}}_{\s{C}^{\alpha/2}\del{[0,T];Y_{n+\alpha}^*}} \leq C J_n({r/2}), \quad n=1,2.
	\end{equation}
	On the other hand,
	\begin{equation}
		\enVert{q^*(0,\cdot)}_{\s{C}^{\alpha/2}([0,T])} \leq C
	\end{equation}
	by the H\"older regularity of $\pd{u}{x}$.
	By Lemma \ref{lem:eta_rho holder}, we deduce
	\begin{equation} \label{eq:g32}
		\enVert{g_{3,2}}_{\s{C}^{\alpha/2}([0,T])} \leq C\del{\enVert{\mu_0}_{\sr{M}_{\alpha/2}^{-n}} + \tilde K_n({r/2})}.
	\end{equation}
	Finally, we analyze $g_4$ in a similar way.
	Write
	\begin{equation}
		g_4(t) =  -\frac{1}{2}\int_{\s{D}} \del{\dpd{w_{r/2}}{x}(\cdot,t) - \dpd{w_{r/2}}{x}(0,t)}\dif m(t)  -\frac{1}{2}\dpd{w_{r/2}}{x}(0,t)\int_{\s{D}} \dif m(t).
	\end{equation}
	Using Lemmas \ref{lem:fp} and \ref{lem:fp eta holder} applies to the solution $m$ of System \eqref{eq:mfg lin dem}, we deduce
	\begin{equation} \label{eq:g4}
		\enVert{g_4}_{\s{C}^{\alpha/2}([0,T])}
		\leq \enVert{\dpd{w_{r/2}}{x}}_{\s{C}^{\alpha,\alpha/2}}\del{\enVert{m}_{\s{C}^{\alpha/2}\del{[0,T];\s{C}_\diamond^\alpha(\s{D})^*} } + \enVert{\int_{\s{D}} \dif m(\cdot)}_{\s{C}^{\alpha/2}([0,T])}} \leq C\enVert{\dpd{w_{r/2}}{x}}_{\s{C}^{\alpha,\alpha/2}}.
	\end{equation}
	Combining \eqref{eq:g1}, \eqref{eq:g2}, \eqref{eq:g3}, \eqref{eq:g31}, \eqref{eq:g32}, and \eqref{eq:g4}, we obtain \eqref{eq:sQ Holder}.
\end{proof}

\subsection{Full regularity of $w$}

Multiply Equation \eqref{eq:linearized}(i) by $e^{-\rho t}$ to see that $w_\rho$ satisfies
\begin{equation} \label{eq:w_rho}
	\dpd{w_\rho}{t} + \dfrac{\sigma^2}{2}\dpd[2]{w_\rho}{x} + V_1(x,t)\dpd{w_\rho}{x} + V_2(x,t)\s{Q}_\rho(t) - (r-\rho)w_\rho = f_\rho.
\end{equation}
In this section we will derive an estimate on $w_{r/2}$ in classical H\"older spaces.
In particular, let us define $Z_\alpha(T)$ to be the set of all $w \in \s{C}^{2+\alpha,1+\alpha/2}(\overline{\s{D}} \times [0,T])$ such that $\psi \dpd{w}{x} \in \s{C}^{2+\alpha,1+\alpha/2}(\overline{\s{D}} \times [0,T])$ as well.
(As usual, when $T = \infty$ we replace $[0,T]$ with $\intco{0,\infty}$.)
It is a Banach space with norm
\begin{equation*}
	\enVert{w}_{Z_\alpha} = \enVert{w}_{\s{C}^{2+\alpha,1+\alpha/2}(\overline{\s{D}} \times [0,T])}
	+
	\enVert{\psi \dpd{w}{x}}_{\s{C}^{2+\alpha,1+\alpha/2}(\overline{\s{D}} \times [0,T])}.
\end{equation*}

\begin{theorem} \label{thm:full w}
	Let $(w,\mu)$ be a solution of \eqref{eq:linearized}, with $V_1,\ldots,V_5$ satisfying either \eqref{eq:V_i for derivative} or \eqref{eq:V_i for differences}.
	Then there is a constant $C\del{r,\sigma,\alpha}$, not depending on $T$, such that
	\begin{equation} \label{eq:full w}
		\enVert{w}_{Z_\alpha} 
		\leq C\del{r,\sigma,\alpha} \del{ \enVert{\mu_0}_{\sr{M}_{\alpha/2}^{-2}}
			+ N(f) + N(\nu)
		},
	\end{equation}
	where		
	\begin{equation} \label{eq:f norms}
		N(f) := \enVert{f_{r/2}}_{\s{C}^{\alpha,\alpha/2}} + \enVert{\psi\dpd{f_{r/2}}{x}}_{\s{C}^{\alpha,\alpha/2}} + \del{\int_0^T \enVert{f_{r/2}(\cdot,s)}_2^2 \dif s}^{1/2}
	\end{equation}
	and
	\begin{equation}
		\label{eq:nu norms}
		N^*(\nu) := \enVert{\ip{\nu_{r/2}}{1}}_{\s{C}^{\alpha/2}([0,T])} +
		\del{\int_0^T \enVert{\nu_{r/2}(s)}_{-2}^2 \dif s}^{1/2}
		+ 
		\sup_{0 \leq \tau \leq T} \enVert{\nu_{r/2}(\tau)}_{-2}.
	\end{equation}
\end{theorem}

\begin{proof}
	\firststep
	We will first apply the maximum principle to find a bound on $w_{r/2}$.
	Let
	\begin{equation*}
		\tilde w(x,t) = w_r(x,t) - \int_t^T \del{\enVert{V_2}_0 \abs{\s{Q}_r(s)} + \enVert{f_r(\cdot,s)}_0}\dif s
	\end{equation*}
	and differentiate to see that
	\begin{equation}
		-\dpd{\tilde w}{t} - \dfrac{\sigma^2}{2}\dpd[2]{\tilde w}{x} + V_1(x,t)\dpd{\tilde w}{x}
		\leq 0.	
	\end{equation}
	By the maximum principle, using the fact that $\tilde w(0,t) \leq 0$ for all $t$ and $\tilde w(x,T) = 0$ for all $x$, we have
	\begin{equation}
		\tilde w(x,t) \leq 0
		\quad \Rightarrow \quad
		w_r(x,t) \leq \int_t^T \del{\enVert{V_2}_0 \abs{\s{Q}_r(s)} + \enVert{f_r(\cdot,s)}_0}\dif s.
	\end{equation}
	Multiply by $e^{rt/2}$ and use the Cauchy-Schwartz inequality to get
	\begin{equation}
		\begin{split}
			w_{r/2}(x,t) &\leq \int_t^T e^{\frac{r}{2}(t-s)}\del{\enVert{V_2}_0 \abs{\s{Q}_{r/2}(s)} + \enVert{f_{r/2}(\cdot,s)}_0}\dif s\\
			&\leq r^{-1/2}\del{\int_t^T \del{\enVert{V_2}_0 \abs{\s{Q}_{r/2}(s)} + \enVert{f_{r/2}(\cdot,s)}_0}^2\dif s}^{1/2}\\
			&\leq C(r)\del{\int_0^T \abs{\s{Q}_{r/2}(s)}^2 \dif s}^{1/2} + C(r)\del{\int_0^T \enVert{f_{r/2}(\cdot,s)}_0^2 \dif s}^{1/2}.
		\end{split}
	\end{equation}
	Applying Corollary \ref{cor:energy estimate2} \ref{cor:int sQ estimate}, we see that
	\begin{equation}
		w_{r/2}(x,t) \leq C J_2(r)^{1/2}.
	\end{equation}
	By the same argument applied that $-w$, we deduce
	\begin{equation} \label{eq:w bound}
		\enVert{w_{r/2}}_0 \leq C J_2(r)^{1/2}.
	\end{equation}
	
	\nextstep
	If we apply Lemma \ref{lem:3rd order holder} to \eqref{eq:w_rho} with $\rho = r/2$, we obtain an estimate
	\begin{multline} \label{eq:full w1}
		\enVert{w_{r/2}}_{Z_\alpha} \leq C\del{\enVert{V_1}_{\s{C}^{\alpha,\alpha/2}},\enVert{\dpd{V_1}{x}}_{\s{C}^{\alpha,\alpha/2}},r,\sigma,\alpha}\\
		\times \del{\enVert{f_{r/2} - V_2 \s{Q}_{r/2}}_{\s{C}^{\alpha,\alpha/2}} + \enVert{\psi\del{\dpd{f_{r/2}}{x} - \dpd{V_2}{x} \s{Q}_{r/2}}}_{\s{C}^{\alpha,\alpha/2}}}.
	\end{multline}
	The H\"older norms of $V_1,V_2, \pd{V_1}{x}$ and $\pd{V_2}{x}$ are already estimated by the estimates \eqref{eq:main estimates} from Theorem \ref{thm:exist mfg infty}.
	Moreover, $f_{r/2}$ is given.
	Using Equation \eqref{eq:sQ Holder} from Corollary \ref{cor:sQ Holder} in \eqref{eq:full w1}, we obtain
	\begin{equation} \label{eq:full w2}
		\enVert{w_{r/2}}_{Z_\alpha} 
		\leq C\del{r,\sigma,\alpha} \del{\enVert{f_{r/2}}_{\s{C}^{\alpha,\alpha/2}} + \enVert{\psi\dpd{f_{r/2}}{x}}_{\s{C}^{\alpha,\alpha/2}}
			+ \tilde J_2({r/2}) + \enVert{\dpd{w_{r/2}}{x}}_{\s{C}^{\alpha,\alpha/2}}}.
	\end{equation}
	By using the interpolation inequality
	\begin{equation*}
		\enVert{\dpd{w_{r/2}}{x}}_{\s{C}^{\alpha,\alpha/2}} \leq \varepsilon \enVert{w_{r/2}}_{\s{C}^{2+\alpha,1+\alpha/2}(\overline{\s{D}} \times [0,T])}
		+ C(\varepsilon)\enVert{w_{r/2}}_0
	\end{equation*} 
	and applying \eqref{eq:w bound}, using the fact that $J_2(r/2)^{1/2} \leq \tilde J_2(r/2)$, estimate \eqref{eq:full w2} yields
	\begin{equation} \label{eq:full w2-2}
		\enVert{w_{r/2}}_{Z_\alpha} 
		\leq C\del{r,\sigma,\alpha} \del{\enVert{f_{r/2}}_{\s{C}^{\alpha,\alpha/2}} + \enVert{\psi\dpd{f_{r/2}}{x}}_{\s{C}^{\alpha,\alpha/2}}
			+ \tilde J_2({r/2})}.
	\end{equation}
	We now return to the definition of $\tilde J_2$, Equation \eqref{eq:tilde J_n}, which can be written
	\begin{multline} \label{eq:tilde J_n1}
		\tilde J_2(r/2) := \enVert{\mu_0}_{\sr{M}_{\alpha/2}^{-2}}  + \enVert{\ip{\nu_{r/2}}{1}}_{\s{C}^{\alpha/2}([0,T])} + \sup_{0 \leq \tau \leq T} \enVert{\nu_{r/2}(\tau)}_{-2}\\
		+ J_2({r/2})^{1/2} + \enVert{\dpd{w_{r/2}}{x}}_0^{1/2}J_2({r/2})^{1/4} + \enVert{\dpd{w_{{r/2}}}{x}}_0^{2/3}J_2({r/2})^{1/6}.
	\end{multline}
	Now since $\enVert{\pd{w_{r/2}}{x}}_{0}$ is dominated by $\enVert{w_{r/2}}_{\s{C}^{2+\alpha,1+\alpha/2}(\overline{\s{D}} \times [0,T])}$, we apply Young's inequality to \eqref{eq:tilde J_n} to get
	\begin{multline} \label{eq:tilde J_n2}
		\tilde J_2({r/2}) \leq  \enVert{\mu_0}_{\sr{M}_{\alpha/2}^{-2}} +  \enVert{\ip{\nu}{1}}_{\s{C}^{\alpha/2}([0,T])} + \sup_{0 \leq \tau \leq T} \enVert{\nu_\rho(\tau)}_{-2}\\
		+ C(\varepsilon)J_2({r/2})^{1/2}
		+ \varepsilon \enVert{w_{r/2}}_{\s{C}^{2+\alpha,1+\alpha/2}(\overline{\s{D}} \times [0,T])}.
	\end{multline}
	Applying \eqref{eq:tilde J_n} to \eqref{eq:full w2-2}, we derive,
	using the definition of $J_2$ in \eqref{eq:J_n},
	\begin{equation} \label{eq:full w3}
		\enVert{w_{r/2}}_{Z_\alpha} 
		\leq C\del{r,\sigma,\alpha} \del{\enVert{\mu_0}_{\sr{M}_{\alpha/2}^{-2}} + \enVert{\mu_0}_{-2} +
			\enVert{w(\cdot,0)}_2^{1/2} \enVert{\mu_0}_{-2}^{1/2}
			+ N^*(\nu) + N(f)}.
	\end{equation}
	where $N(f)$ and $N^*(\nu)$ are defined in \eqref{eq:f norms} and \eqref{eq:nu norms}, respectively.
	Using the fact that $\enVert{w(\cdot,0)}_2$ is dominated by $\enVert{w}_{Z_\alpha}$, we apply Young's inequality to \eqref{eq:full w3} and rearrange to deduce \eqref{eq:full w}.
\end{proof}

\subsection{An existence theorem for the linearized system}

Before formulating the main result of this section, let us collect assumptions on $r$ and $\epsilon$ so that all of the a priori estimates of this section hold.
We will formulate two alternatives, one for a linear demand schedule, and one for a more general case where $\epsilon$ must be small.
\begin{assumption}[$r$ big, $\epsilon$ small]
	\label{as:r big, epsilon small}
	Let $r^*$ be a number large enough to satisfy Assumption \ref{as:r bigger than kappa}, Equation \eqref{eq:r assm n} for $n = 2$, and
	\begin{equation*} 
		2\sqrt{\frac{2}{\sigma^2 r}}H(0,0,0) < P(0) \quad \forall r \geq r^*.
	\end{equation*}
	Let $\epsilon^* > 0$ be small enough to satisfy \eqref{eq:epsilon small} and
	\begin{equation*}
		M = 2\sqrt{\frac{2}{\sigma^2 r^*}}H(0,0,0) < P\del{\epsilon^* \bar Q}.
	\end{equation*} 
	We assume that $r \geq r^*$ and $0 < \epsilon \leq \epsilon^*$.
\end{assumption}
We remark that Assumption \ref{as:r big, epsilon small} implies Assumption \ref{as:H must be smooth}; see Remark \ref{rmk:H must be smooth}.

An alternative assumption is as follows.
\begin{assumption}[$r$ big, $P$ linear]
	\label{as:r big, P linear}
	We assume that $P(q) = 1-q$ and that $0 < \epsilon < 2$.
	Let $r^*$ be a number large enough to satisfy Assumption \ref{as:r bigger than kappa}, Equation \eqref{eq:r assm n} for $n = 2$, and
	\begin{equation*} 
		2\sqrt{\frac{2}{\sigma^2 r}}H(0,0,0) < 1 - \frac{\epsilon}{2} \quad \forall r \geq r^*.
	\end{equation*}
	We assume that $r \geq r^*$.
\end{assumption}

\begin{theorem} \label{thm:linearized exist T}
	Let $T > 0$ be a fixed time horizon.
	Suppose that Assumption \ref{as:r big, epsilon small} or \ref{as:r big, P linear} holds.
	Then System \eqref{eq:linearized} has a unique solution $(w,\mu)$ satisfying $w(x,T) = 0$ with regularity
	\begin{itemize}
		\item $w_{r/2} \in Z_\alpha(T)$,
		\item $\mu_{r/2} \in \s{C}^{\alpha/2}\del{[0,T]; Y_{n + \alpha}^*} \cap L^\infty\del{(0,T); X_{2}^*} =: \tilde Z_\alpha(T)$.
	\end{itemize}
	There exists a constant $C(r,\sigma,\alpha)$, not depending on $T$, such that
	\begin{multline} \label{eq:linearized a priori estimates}
		\enVert{w_{r/2}}_{Z_\alpha} + \del{\int_0^T \enVert{\dpd{w_{r/2}}{x}(\cdot,t)}_2^2 \dif t}^{1/2}
		+ \enVert{\mu_{r/2}}_{\s{C}^{\alpha/2}\del{[0,T]; Y_{2 + \alpha}^*}}
		+ N^*(\mu)\\
		\leq C\del{r,\sigma,\alpha} \del{\enVert{\mu_0}_{\sr{M}_{\alpha/2}^{-2}}
			+ N(f) + N^*(\nu)
		},
	\end{multline}
	where $N(f)$ and $N^*(\nu)$ are defined in \eqref{eq:f norms} and \eqref{eq:nu norms}, respectively.
\end{theorem}

\begin{proof}
	First we assume the data are smooth.
	Then existence of solutions follows from the Leray-Schauder fixed point theorem, along the same lines as in the proof of Theorem \ref{thm:exist mfg T}.
	The a priori estimates \eqref{eq:linearized a priori estimates} follow from Lemmas \ref{thm:full w} and \ref{lem:fp holder} (Equations \eqref{eq:full w} and \eqref{eq:fp holder}).
	A similar argument is also found in \cite[Lemma 3.3.1]{cardaliaguet2019master}.
	To see that the solution is unique, note that the system is linear, so the a priori bounds also imply uniqueness.
\end{proof}

\begin{theorem} \label{thm:linearized exist infty}
	Suppose that Assumption \ref{as:r big, epsilon small} or \ref{as:r big, P linear} holds.
	Then System \eqref{eq:linearized} has a unique solution $(w,\mu)$ satisfying
	\begin{itemize}
		\item $w_{r/2} \in Z_\alpha(\infty)$,
		\item $\mu_{r/2} \in \s{C}^{\alpha/2}\del{\intco{0,\infty}; Y_{2 + \alpha}^*} \cap L^\infty\del{(0,\infty); X_{2}^*} =: \tilde Z_\alpha(\infty)$,
		\item $\lim_{t \to \infty}\enVert{w_{r/2}(\cdot,t)}_2 = \lim_{t \to \infty}\enVert{\dpd{w_{r/2}}{x}(\cdot,t)}_2 = 0$.
	\end{itemize}
	The estimate \eqref{eq:linearized a priori estimates} holds with $T = \infty$.
\end{theorem}

\begin{proof}
	For each $T > 0$, let $(w^T,\mu^T)$ be the solution to the finite time horizon problem on $[0,T]$ given by Theorem \ref{thm:linearized exist T}.
	We extend $(w^T,\mu^T)$ in time such that $w^T(x,t) = 0$ for all $t \geq T$ and such that
	the a priori estimate \eqref{eq:linearized a priori estimates} implies that $(w^T_{r/2},\mu^T_{r/2})$ is bounded in $(Z_\alpha(\infty),\tilde Z_\alpha(\infty))$.
	Then by standard compactness arguments there exists a subsequence $T_n \to \infty$ such that $(w^{T_n}_{r/2},\mu^{T_n}_{r/2})$ converges to some $(w_{r/2},\mu_{r/2})$ in $Z_\beta(\infty) \times \s{C}^{\beta/2}\del{[0,T];Y_{2+\alpha}^*}$ for $\beta < \alpha$.
	Moreover, $(w,\mu)$ satisfies \eqref{eq:linearized a priori estimates} with $T = \infty$.
	Passing to the limit in the system satisfied by $(w^{T_n},\mu^{T_n})$, we deduce that $(w,\mu)$ satisfies System \eqref{eq:linearized}.
	It follows that $(w,\mu)$ is a solution.
	To see that $\lim_{t \to \infty}\enVert{w_{r/2}(\cdot,t)}_2 = \lim_{t \to \infty}\enVert{\dpd{w_{r/2}}{x}(\cdot,t)}_2 = 0$, we observe that since $w_{r/2}^{T_n} \to w_{r/2}$ in $Z_\beta(\infty)$, it follows that $w_{r/2}^{T_n} \to w_{r/2}$ and $\dpd{w_{r/2}^{T_n}}{x} \to \dpd{w_{r/2}}{x}$ in $\s{C}\del{\intco{0,\infty};X_2}$.
	Then the fact that $w^{T_n}(x,t) = 0$ for all $t \geq T_n$ implies the desired limit.
	Finally, the a priori bounds together with the linearity of the system imply uniqueness.
\end{proof}

	\section{The solution to the Master Equation}
	
	\label{sec:solution master}
	
	For each $m_0$, define $U(m_0,x) = u(x,0)$ where $(u,m)$ is the solution of \eqref{eq:mfg infty} given initial condition $m_0$.
	We refer to $U(m_0,x)$ as the \emph{master field}.
	We will prove that it satisfies the master equation \eqref{eq:master equation}.
	All the hypotheses of Theorem \ref{thm:exist mfg infty} plus Assumption \ref{as:r big, epsilon small} or \ref{as:r big, P linear} are in force.

	\subsection{Continuity and differentiability of the master field}
	
	In this subsection we show that $U(m_0,x)$ is Lipschitz continuous and differentiable with respect to the measure variable $m_0$.
	To do this, we appeal to the estimates found in Section \ref{sec:sensitivity}.
	
	\begin{theorem} \label{thm:U lipschitz}
		There exists a constant $C$ such that
		\begin{equation}
			\enVert{U(\tilde m_0,\cdot) - U(m_0,\cdot)}_{Y_{n+1+\alpha}}
			\leq C\enVert{\tilde m_0 - m_0}_{\sr{M}_{\alpha/2}^{-2}} \quad \forall \tilde m_0, m_0 \in \sr{M}_{\alpha/2}^{-2}.
		\end{equation}
	\end{theorem}

	\begin{proof}
		We may assume that $\tilde m_0,m_0 \in \sr{M}_{\alpha/2}$; then by density of this set in $\sr{M}_{\alpha/2}^{-2}$, we deduce the result.
		We have $U(m_0,x) = u(x,0)$ and $U(\tilde m_0,x) = \tilde u(x,0)$, where $(u,m)$ is the solution of \eqref{eq:mfg lin dem} given initial condition $m_0$ and $(\tilde u,\tilde m)$ is the solution of \eqref{eq:mfg infty} given initial condition $\tilde m_0$.
		Let $(w,\mu) = (\tilde u - u,\tilde m - m)$.
		Then $(w,\mu)$ solves the linearized system \eqref{eq:linearized} with $f = 0$, $\nu = 0$, and $V_1,\ldots,V_5$ defined in \eqref{eq:V_i for differences}.
		Observe that
		\begin{equation}
			\enVert{U(\tilde m_0,\cdot) - U(m_0,\cdot)}_{Y_{3+\alpha}}
			= \enVert{w(\cdot,0)}_{Y_{3+\alpha}}
			\leq \enVert{w}_{Z_\alpha}.
		\end{equation}
		We conclude by appealing to Theorem \ref{thm:full w}.
	\end{proof}

	Before proving that $U$ is differentiable with respect to $m$, we provide a candidate for the derivative in the following lemma.
	\begin{lemma} \label{lem:derivative of U}
		Let $f = 0$ and $\nu = 0$.
		There exists a map $K(m_0,x,y)$ such that $K$ is thrice differentiable with respect to $x$ and twice differentiable with respect to $y$, such that
		\begin{equation} \label{eq:K estimates}
			\enVert{\dpd[\ell]{K}{y}(m_0,\cdot,y)}_{Y_{3+\alpha}} \leq C\max\cbr{\abs{y}^{-\alpha-\ell},1},
		\end{equation}
		and such that if $(w,\mu)$ is the solution of System \eqref{eq:linearized}, then
		\begin{equation} \label{eq:w(x,0)=}
			w(x,0) = \ip{K(m_0,x,\cdot)}{\mu_0}.
		\end{equation}
	Moreover, $K$ and its derivatives in $(x,y)$ are continuous with respect to the topology on $\s{M}_+(\s{D}) \times \s{D} \times \s{D}$.
	\end{lemma}

	\begin{proof}
		The proof is very similar to that of \cite[Corollary 3.4.2]{cardaliaguet2019master}: for $\ell = 0,1,2$ and $y > 0$ let the pair
		 $(w^{(\ell)}(\cdot,\cdot,y),\mu^{(\ell)}(\cdot,\cdot,y))$ be the solution of \eqref{eq:linearized} with $f = 0, \nu = 0$, $V_1,\ldots,V_5$ given by \eqref{eq:V_i for derivative} and initial condition $\mu_0 = D^\ell \delta_y$, where $\delta_y$ is the Dirac delta mass concentrated at $y$ and $D^\ell \delta_y$ is its $\ell$th derivative in the sense of distributions.
		Then set $K(m_0,x,y) = w^{(0)}(x,0,y)$.
		Notice that by the density of empirical measures, \eqref{eq:w(x,0)=} follows for any solution $(w,\mu)$ of System \eqref{eq:linearized}.
		Moreover, one can check by induction that
		\begin{equation*}
			\dpd[\ell]{K}{y}(m_0,x,y) = (-1)^\ell w^{(\ell)}(x,0,y).
		\end{equation*}
		To prove \eqref{eq:K estimates}, we use the estimates \eqref{eq:linearized a priori estimates} from Theorem \ref{thm:linearized exist infty}, which imply in particular that
		\begin{equation}
			\enVert{w^{(\ell)}(\cdot,0,y)}_{Y_{3+\alpha}} \leq C\enVert{D^\ell \delta_y}_{\sr{M}_{\alpha/2}^{-2}}.
		\end{equation}
		It remains only to estimate $D^\ell \delta_y$ in $\sr{M}_{\alpha/2}^{-2}$.
		First, we see that
		\begin{equation}
			\ip{\phi}{D^\ell \delta_y} = \dod[\ell]{\phi}{y}(y) \leq \enVert{\phi}_2 \max\cbr{\abs{y}^{-\ell},1} \quad \forall \phi \in X_{2}
			\quad
			\Rightarrow
			\quad
			\enVert{D^\ell \delta_y}_{-2} \leq \max\cbr{\abs{y}^{-\ell},1}.
		\end{equation}
		Next, we plug $\mu = D^\ell \delta_y$ into \eqref{eq:eta function -n} to get
		\begin{equation} \label{eq:etahDelldelta}
			\eta^h[D^\ell \delta_y](t) =  \int_{0}^\infty \del{(-1)^{\ell+1}\dpd[\ell]{S}{x}(x-y,t) + \dpd[\ell]{S}{x}(x+y,t)}\dif x
			= - 2\dpd[\ell-1]{S}{y}(y,t),
		\end{equation}
		where we define
		\begin{equation}
			\dpd[-1]{S}{y}(y,t) = \int_0^y S(x,t)\dif x.
		\end{equation}
		Taking the derivative with respect to $t$ in \eqref{eq:etahDelldelta}, we get
		\begin{equation}
			\dod{}{t}\eta^h[D^\ell \delta_y](t) = -\sigma^2 \dpd[\ell+1]{S}{y}(y,t).
		\end{equation}
		Let $p > 1$.
		Applying Lemma \ref{lem:m_n} we estimate 
		\begin{equation}
			\int_0^\infty \abs{\dod{}{t}\eta^h[D^\ell \delta_y](t)}^p \dif t \leq C(\sigma)\int_0^\infty t^{-p(\ell + 2)/2}P_{\ell + 1}\del{\frac{\abs{y}}{\sqrt{t}}}^p\exp\cbr{-\frac{py^2}{2\sigma^2 t}}\dif t.
		\end{equation}
		Here $P_{\ell + 1}$ is a polynomial of degree $\ell + 1$.
		Using the substitution $s = \frac{py^2}{t}$ we obtain
		\begin{equation*}
			\del{\int_0^\infty \abs{\dod{}{t}\eta^h[D^\ell \delta_y](t)}^p \dif t}^{1/p} \leq C(\sigma,p)y^{\frac{2}{p} - (\ell + 2)}\del{\int_0^\infty s^{\frac{p(\ell + 2)}{2}-2}P_{\ell+1}(s^{1/2})^p\exp\cbr{-\frac{s}{2\sigma^2}} \dif s}^{1/p},
		\end{equation*}
		where the integral on the right-hand side converges; hence
		\begin{equation}
			\del{\int_0^\infty \abs{\dod{}{t}\eta^h[D^\ell \delta_y](t)}^p \dif t}^{1/p} \leq C(\sigma,p)y^{-\frac{2}{p'} - \ell}, \quad p' := p/(p-1).
		\end{equation}
		By H\"older's inequality,
		\begin{equation}
			\abs{\eta^h[D^\ell \delta_y](t_1) - \eta^h[D^\ell \delta_y](t_2)} \leq C(\sigma,p)y^{-\frac{2}{p'} - \ell}\abs{t_1-t_2}^{1/p'}
		\end{equation}
		Cf.~the proof of Lemma \ref{pr:inverse moment}.
		Choosing $p = \frac{2}{2-\alpha}$, we now deduce
		\begin{equation}
			\enVert{\eta^h[D^\ell \delta_y]}_{\s{C}^{\alpha/2}\del{\intco{0,\infty}}} \leq C(\sigma,\alpha)y^{-\alpha - \ell}.
		\end{equation}
		Therefore,
		\begin{equation}
			\enVert{D^\ell \delta_y}_{\sr{M}_{\alpha/2}^{-2}} 
			\leq C(\sigma,\alpha)\max\cbr{\abs{y}^{-\alpha-\ell},1},
		\end{equation}
		from which we deduce \eqref{eq:K estimates}.
		The remaining details are the same as in \cite[Corollary 3.4.2]{cardaliaguet2019master}.
	\end{proof}

	\begin{lemma} \label{lem:differentiability of U}
		Let $(u,m)$ and $(\hat u,\hat m)$ be the solutions of \eqref{eq:mfg infty} with initial conditions $m_0$ and $\hat m_0$, respectively.
		Let $(w,\mu)$ be the solution of \eqref{eq:linearized} with initial condition $\hat m_0 - m_0$.
		Then
		\begin{equation} \label{eq:diffble estimate}
			\enVert{\hat u(\cdot,0) - u(\cdot,0) - w(\cdot,0)}_{Y_{3+\alpha}} \leq C\enVert{\hat m_0 - m_0}_{\sr{M}_{\alpha/2}^{-2}}^2.
		\end{equation}
		As a corollary, $U(m_0,x)$ is differentiable with respect to $m_0$ with
		\begin{equation} \label{eq:dUdm = K}
			\vd{U}{m}(m_0,x,y) = K(m_0,x,y),
		\end{equation}
		where $K$ is defined in Lemma \ref{lem:derivative of U}.
		Moreover, \eqref{eq:diffble estimate} reads
		\begin{equation} \label{eq:diffble estimate2}
			\enVert{U(\hat m_0,\cdot) - U(m_0,\cdot) - \int_{\s{D}} \vd{U}{m}(m_0,x,y) \dif (\hat m_0 - m_0)(y)}_{Y_{3+\alpha}} \leq C\enVert{\hat m_0 - m_0}_{\sr{M}_{\alpha/2}^{-2}}^2.
		\end{equation}
	\end{lemma}

	\begin{proof}
		Let $f$ and $\nu$ be defined by \eqref{eq:f nu}.
		We follow the same steps as in \cite[Chapter 3]{cardaliaguet2019master} to find an estimate
		\begin{equation}
			N(f) + N^*(\nu) \leq C\enVert{\hat m_0 - m_0}_{\sr{M}_{\alpha/2}^{-2}}^2.
		\end{equation}
		By Theorem \ref{thm:linearized exist infty}, this proves \eqref{eq:diffble estimate}.
		Combined with Lemma \ref{lem:derivative of U}, we deduce \eqref{eq:dUdm = K} and \eqref{eq:diffble estimate2}.
	\end{proof}
	\subsection{The master field satisfies the master equation}
	
	In this subsection we prove Theorem \ref{thm:main result}.
	\begin{theorem}
		Suppose that Assumption \ref{as:r big, epsilon small} or \ref{as:r big, P linear} holds.
		For all $m_0 \in \s{M}^{2+\alpha}$, $x \in \s{D}$, System \eqref{eq:master equation}-\eqref{eq:clearing condition} is satisfied.
	Moreover, $U$ is the unique continuously differentiable function satisfying
	\begin{equation} \label{eq:dUdm estimate}
		\enVert{\dpd[\ell]{}{y}\vd{U}{m}(m_0,\cdot,y)}_{Y_{3+\alpha}} \leq C\max\cbr{\abs{y}^{-\alpha-\ell},1},
	\end{equation}
	such that System \eqref{eq:master equation}-\eqref{eq:clearing condition} holds for all $m_0 \in \s{M}^{2+\alpha}$, $x \in \s{D}$.
	\end{theorem}
	\begin{proof}
		Let $(u,m)$ be the solution to the mean field game system with initial condition $m_0 \in \s{M}^{2+\alpha}$.
	Set $m_s := sm(t) + (1-s)m_0$ for $0 \leq s \leq 1$.
	Then for any $t > 0$ we have
	\begin{multline} \label{eq:u(x,t)-u(x,0)}
		u(x,t) - u(x,0) =
		U(m(t),x) - U(m_0,x) = \int_0^1 \int_{\s{D}} \vd{U}{m}(m_s,x,y)\dif \del{m(t)-m_0}(y)\dif s\\
		= \int_0^1 \int_0^t \int_{\s{D}} \del{\frac{\sigma^2}{2}\dpd[2]{}{y}\vd{U}{m}(m_s,x,y)
			+ \dpd{H}{a}\del{\epsilon,Q^*(\tau),\dpd{u}{x}(y,\tau)}\dpd{}{y}\vd{U}{m}(m_s,x,y)}\dif m(\tau)(y)\dif \tau\dif s,
	\end{multline}
	using the Fokker-Planck equation satisfied by $m$.
	To see that the last integral converges, first note that \eqref{eq:dUdm estimate} holds by Lemmas \ref{lem:derivative of U} and \ref{lem:differentiability of U}.
	Then we note that by the assumption $m_0 \in \s{M}^{2+\alpha}$ together with Lemma \ref{lem:Malpha fp bound},
	\begin{equation}
		\int_{\s{D}} \del{1+ x^{-(2+\alpha)}}m(\dif x,t) \leq Ce^{Ct}\int_{\s{D}} \del{1+ x^{-(2+\alpha)}}m_0(\dif x).
	\end{equation}
	Combining this with \eqref{eq:dUdm estimate}, we deduce that \eqref{eq:u(x,t)-u(x,0)} holds.
	Now divide by $t$ and let $t \to 0$ to get
	\begin{equation}
		\dpd{u}{t}(x,0) = \int_{\s{D}} \del{\frac{\sigma^2}{2}\dpd[2]{}{y}\vd{U}{m}(m_0,x,y)
			+ \dpd{H}{a}\del{\epsilon,Q^*(0),\dpd{u}{x}(y,0)}\dpd{}{y}\vd{U}{m}(m_0,x,y)}\dif m_0(y).
	\end{equation}
	By substituting for $\dpd{u}{t}(x,0)$ using Equation \eqref{eq:mfg infty}(i), we get
	\begin{multline}
		-\dfrac{\sigma^2}{2}\dpd[2]{u}{x}(x,0) - H\del{\epsilon,Q^*(0),\dpd{u}{x}(x,0)} + ru(x,0)\\
		= \int_{\s{D}} \del{\frac{\sigma^2}{2}\dpd[2]{}{y}\vd{U}{m}(m_0,x,y)
			+ \dpd{H}{a}\del{\epsilon,Q^*(0),\dpd{u}{x}(y,0)}\dpd{}{y}\vd{U}{m}(m_0,x,y)}\dif m_0(y),
	\end{multline}
	which becomes Equation \eqref{eq:master equation} after defining $Q^* = Q^*(0)$.
	Equation \eqref{eq:clearing condition} follows from \eqref{eq:mfg infty}(iii).
	
	To see that $U$ is unique, we follow the same argument as in \cite{cardaliaguet2019master}.
	By using the Leray-Schauder fixed point theorem and the estimates we have established, it is straightforward to show the existence of a solution to the Fokker-Planck equation
	\begin{equation*}
		\dpd{m}{t} - \frac{\sigma^2}{2}\dpd[2]{m}{x} + \dpd{}{x}\del{\dpd{H}{a}\del{\epsilon,Q^*(m(t)),\dpd{U}{x}(m(t),x)}} = 0,
	\end{equation*}
	where $Q^*(m)$ is defined using \eqref{eq:clearing condition}.
	Set $u(x,t) = U(m(t),x)$.
	Using condition \eqref{eq:dUdm estimate} together with Lemma \ref{lem:Malpha fp bound}, as above, we can differentiate $u$ with respect to time.
	Then using the fact that \eqref{eq:master equation} holds, we deduce that $(u,m)$ is the solution of \eqref{eq:mfg infty}, which is unique.
	It follows that $U(m,x)$ is uniquely determined.
	\end{proof}

	\bibliographystyle{siam}
	\bibliography{mybib}
	\appendix

\section{Proofs of Results from Section \ref{sec:fokkerplanck}}
\label{ap:fokkerplanck proofs}

\begin{proof}[Proof of Lemma \ref{lem:fp}]
	\emph{Uniqueness:}
	Let us start by observing that uniqueness of weak solutions holds. This follows from a proof by duality, cf.~\cite[Proposition B.1]{graber2020commodities} and \cite[Corollary 3.5]{porretta2015weak}, which also provide the basic estimate \eqref{eq:TV}.
	
	\emph{Existence:} We thus turn our attention to existence and estimates. By linearity we can assume that $m_0(\s{D}) = 1$, i.e.~$m_0$ is a probability measure, without loss of generality.
	
	Assume for now that $b$ is infinitely smooth and bounded, and that $m_0 \in \s{P}_1(\s{D})$ is in fact a smooth density such that $m_0 \in \s{C}^\infty_c(\overline{\s{D}})$.
	Then classical theory \cite[Theorems IV.5.2, IV.9.1]{ladyzhenskaia1968linear} implies that \eqref{eq:fp} has a smooth solution $m$ whose derivatives are also in $L^p$ for arbitrarily large $p$.
	We have the following probabilistic interpretation: for any continuous function on $\overline{\s{D}}$ satisfying $$\abs{\phi(x)} \leq C\del{1+\abs{x}},$$
	we have
	\begin{equation} \label{eq:prob interp}
		\int_0^\infty \phi(x) m(x,t)\dif x = \bb{E}\intcc{\phi(X_t)\mathds{1}_{t < \tau}}
	\end{equation}
	where $X_t$ is the diffusion process given by
	\begin{equation} \label{eq:sde}
		\dif X_t = -b(X_t,t)\dif t + \sigma \dif W_t, \ X_0 \sim m_0,
	\end{equation}
	$W_t$ is a standard Brownian motion with respect to a filtered probability space $(\Omega,\bb{P},\s{F}_t)$, and
	\begin{equation}\label{eq:stopping time}
		\tau := \min\cbr{\inf\{t \geq 0 : X_t \leq 0\},T}.
	\end{equation}
	In particular the complementary mass function $\bar\eta(t)$ can be written
	\begin{equation} \label{eq:eta characterized}
		\bar\eta(t) = \bb{P}(t < \tau).
	\end{equation}
	The continuity of this function follows from probabilistic arguments, which can be found in \cite{hambly2017stochastic} and \cite{graber2020commodities}.
	
	It remains to establish \eqref{eq:m holder in time}.
	Fix $t_1,t_2 \in [0,T]$ with $t_1 < t_2$.
	Pick any $\phi:\overline{\s{D}} \to \bb{R}$ that is $\alpha$-H\"older continuous (or Lipschitz, in the case $\alpha = 1$) such that $\phi(0) = 0$ and $[\phi]_{\s{C}^\alpha} \leq 1$.
	Let $X_t$ be a solution to \eqref{eq:sde}.
	Then by \eqref{eq:prob interp} we have
	\begin{equation}
		\begin{split}
			\abs{\int_0^\infty \phi(x)\del{m(x,t_1)-m(x,t_2)}\dif x}
			&= \abs{\bb{E}\intcc{\phi(X_{t_1})\mathds{1}_{t_1 < \tau}-\phi(X_{t_2})\mathds{1}_{t_2 < \tau}}}\\
			&\leq \bb{E}\intcc{\abs{X_{t_1}}^\alpha\mathds{1}_{t_1 < \tau \leq t_2} + \abs{X_{t_1}-X_{t_2}}^\alpha\mathds{1}_{t_2 < \tau}}\\
			&= \bb{E}\intcc{\abs{-\int_{t_1}^\tau b(X_t,t)\dif t + \sigma(W_\tau - W_{t_1})}^\alpha\mathds{1}_{t_1 < \tau \leq t_2}}\\
			&\quad + \bb{E}\intcc{\abs{-\int_{t_1}^{t_2} b(X_t,t)\dif t + \sigma(W_{t_2} - W_{t_1})}^\alpha\mathds{1}_{t_2 < \tau}}\\
			&\leq \bb{E}\intcc{\enVert{b}_\infty^\alpha \abs{\tau-t_1}^\alpha \mathds{1}_{t_1 < \tau \leq t_2}}
			+ \sigma^\alpha \bb{E}\intcc{\abs{W_\tau - W_{t_1}}^\alpha\mathds{1}_{t_1 < \tau \leq t_2}}\\
			&\quad + \enVert{b}_\infty^\alpha \abs{t_2-t_1}^\alpha + \sigma^\alpha \bb{E}\intcc{\abs{W_{t_2} - W_{t_1}}^\alpha}\\
			&\leq 2\enVert{b}_\infty^\alpha \abs{t_2-t_1}^\alpha + 2\sigma^\alpha\abs{t_2-t_1}^{\alpha/2}.
		\end{split}
	\end{equation}
	Taking $t_1 = t$ and $t_2 = 0$, we get
	\begin{equation}
		\begin{aligned}
			\abs{\int_0^\infty \phi(x)m(x,t)\dif x}
			&\leq \abs{\int_0^\infty \phi(x)m_0(x)\dif x} + 2\enVert{b}_\infty^\alpha t^\alpha + 2\sigma^\alpha t^{\alpha/2}\\
			&\leq \int_0^\infty x^\alpha m_0(x)\dif x + 2\enVert{b}_\infty^\alpha t^\alpha + 2\sigma^\alpha t^{\alpha/2}.
		\end{aligned}
	\end{equation}
	
	Finally, to get existence for general data, let $b_n$ be a sequence of smooth functions converging uniformly to $b$ and let $m_{0,n}$ be a sequence of measures with smooth densities converging to $m_0$ in $\s{M}_{1,+}$.
	Letting $m_n$ be the solution corresponding to $b_n,m_{0,n}$, we have that $m_n$ is uniformly H\"older continuous in the $\Lip_{\diamond}(\s{D})^*$ metric, hence by Arzel\'a-Ascoli we have a subsequence converging to $m$ in $\s{C}^0([0,T];\s{M}_{1,+})$.
	We deduce that $m$ is a weak solution, i.e.~it satisfies \eqref{eq:fp weak}.
\end{proof}

\begin{proof}[Proof of Lemma \ref{lem:Malpha bound}]
	For each $n \in \bb{N}$ define
	\begin{equation} \label{eq:phin}
		\phi_n(x) = \begin{cases}
			n^\alpha x &\text{if}~0 < x \leq n^{-1},\\
			x^{-\alpha} &\text{if}~x > n^{-1}.
		\end{cases}
	\end{equation}
	Set $\Phi_n^{(0)}(x) = \phi_n(x)$, and inductively define
	\begin{equation}
		\Phi_n^{(j)}(x) = \int_0^x \Phi_n^{(j-1)}(t)\dif t,
		\quad j = 1,2,3,\ldots.
	\end{equation}
	By induction we have that
	\begin{equation} \label{eq:Phin estimate}
		\abs{\Phi_n^{(j)}(x)} \leq C(j,\alpha)x^{j-\alpha} \quad \forall x > 0.
	\end{equation}
	
	Since $\phi_n$ is a bounded, continuous function, we have
	\begin{equation} \label{eq:phin m}
		\begin{split}
			\int_{\s{D}} \phi_n(x) m(\dif x,t)
			&= \int_{\s{D}} \phi_n(x) \int_{\s{D}} \del{S(x-y,t) - S(x+y,t)}m_0(\dif y) \dif x\\
			&= \int_{\s{D}} \int_{\s{D}}\phi_n(x) \del{S(x-y,t) - S(x+y,t)}\dif x \ m_0(\dif y) 
		\end{split}
	\end{equation}
	using Fubini's Theorem.
	Our goal now is to prove that
	\begin{equation} \label{eq:Malpha bound n}
		\int_{\s{D}}\phi_n(x) \del{S(x-y,t) - S(x+y,t)}\dif x
		\leq C(\alpha)y^{-\alpha} \quad \forall y > 0.
	\end{equation}
	By plugging \eqref{eq:Malpha bound n} into \eqref{eq:phin m} and then applying the Monotone Convergence Theorem, \eqref{eq:Malpha bound} follows.
	
	To prove \eqref{eq:Malpha bound n}, start by noting
	\begin{equation} \label{eq:Malpha bound n1}
		\begin{split}
			\int_{\s{D}}\phi_n(x) \del{S(x-y,t) - S(x+y,t)}\dif x
			&\leq \int_0^{\infty} \phi_n(x) S(x-y,t)\dif x\\
			&\leq \int_0^{y/2} \phi_n(x) S(x-y,t)\dif x
			+ (y/2)^{-\alpha}\int_{y/2}^\infty \phi_n(x) S(x-y,t)\dif x\\
			&\leq \int_0^{y/2} \phi_n(x) S(x-y,t)\dif x + 2^{\alpha}y^{-\alpha},
		\end{split}
	\end{equation}
	using the fact that $S(\cdot,t)$ is a probability density.
	Now for any $j > \alpha - 1$, integrate by parts $j$ times to get
	\begin{multline} \label{eq:Malpha bound n2}
		\int_0^{y/2} \phi_n(x) S(x-y,t)\dif x\\
		= \sum_{i=0}^{j-1} (-1)^i \Phi_n^{(i+1)}(y/2)\dpd[i]{S}{x}(-y/2,t)
		+ (-1)^j \int_0^{y/2} \Phi_n^{(j)}(x)\dpd[j]{S}{x}(x-y,t)\dif x.
	\end{multline}
	Applying \eqref{eq:Phin estimate} and Lemma \ref{lem:m_n} to Equation \eqref{eq:Malpha bound n2}, we obtain
	\begin{equation} \label{eq:Malpha bound n3}
		\begin{split}
			\int_0^{y/2} \phi_n(x) S(x-y,t)\dif x
			&\leq C(j,\alpha)\del{\sum_{i=0}^{j-1} y^{i+1-\alpha}y^{-(i+1)}
				+ \int_0^{y/2} x^{j-\alpha}\abs{x-y}^{-(j+1)}\dif x}\\
			&\leq C(j,\alpha)y^{-\alpha}.
		\end{split}
	\end{equation}
	Take $j = \floor{\alpha}$ and combine \eqref{eq:Malpha bound n3} with \eqref{eq:Malpha bound n1} to obtain \eqref{eq:Malpha bound n}, which completes the proof.
\end{proof}

\begin{proof}[Proof of \ref{lem:fp eta holder}]
	First, note that $\abs{\eta(t)} \leq \enVert{m_0}_{TV} \leq \enVert{m_0}_{\sr{M}_\alpha}$ for all $t \geq 0$, using Lemma \ref{lem:fp}.
	Thus, it suffices to prove estimates of the H\"older constant for $\eta$.
	We will assume the data are sufficiently regular so that the solution is smooth.
	The claim then follows from a density argument.
	
	We have, by Duhamel's principle,
	\begin{multline}
		m(x,t) = \int_0^\infty \del{S(x-y,t)-S(x+y,t)}m_0(y)\dif y\\
		+ \int_0^t\int_0^\infty \del{S(x-y,t-s)-S(x+y,t-s)}(bm)_y(y,s)\dif y\dif s,
	\end{multline}
	which becomes
	\begin{multline} \label{eq:fp formula}
		m(x,t) 
		= \int_0^\infty \del{S(x-y,t)-S(x+y,t)}m_0(y)\dif y\\
		+ \int_0^t\int_0^\infty \del{\dpd{S}{x}(x-y,t-s) + \dpd{S}{x}(x+y,t-s)}(bm)(y,s)\dif y\dif s,
	\end{multline}
	using integration by parts.
	Integrating in $x$ and using Fubini's Theorem, we get
	\begin{equation}
		\eta(t) = \eta^h[m_0](t) + \eta_2(t),
	\end{equation}
	where $\eta^h[m_0](t)$ is defined in \eqref{eq:eta function} and
	\begin{equation} \label{eq:eta holder2-1}
		\begin{aligned}
			\eta_2(t) &=  \int_0^t \int_0^\infty \int_{0}^\infty \del{\dpd{S}{x}(x-y,t-s) + \dpd{S}{x}(x+y,t-s)}(bm)(y,s)\dif x\dif y\dif s\\
			&= -2\int_0^t \int_0^\infty S(y,t-s)(bm)\del{y,s}\dif y\dif s.
		\end{aligned}
	\end{equation}
	By definition of the norm in $\sr{M}_\alpha$,
	\begin{equation}
		\label{eq:eta holder1}
		\enVert{\eta^h[m_0]}_{\s{C}^{\alpha}([0,T])} \leq \enVert{m_0}_{\sr{M}_\alpha}.
	\end{equation}
	It remains to derive H\"older estimates for $\eta_2$.
	Let $t_2 > t_1 \geq 0$.
	Then $\eta_2(t_2) - \eta_2(t_1) = -2\del{I_1 + I_2}$ where
	\begin{equation}
		\begin{aligned}
			I_1 &= \int_{t_1}^{t_2} \int_0^\infty S(y,t_2-s)b(y,s)m(\dif y,s)\dif s,\\
			I_2 &= \int_{0}^{t_1} \int_0^\infty \del{S(y,t_2-s)-S(y,t_1-s)}b(y,s)m(\dif y,s)\dif s.
		\end{aligned}
	\end{equation}
	In the first place, we have
	\begin{equation} \label{eq:int Sbm holder cty1}
		\abs{I_1} \leq (2\sigma^2 \pi)^{-1/2}\enVert{b}_\infty  \int_{t_1}^{t_2} (t_2-s)^{-1/2}\dif s = 2(2\sigma^2 \pi)^{-1/2}\enVert{b}_\infty (t_2-t_1)^{1/2}.
	\end{equation}
	In the second place, we write
	\begin{equation} \label{eq:int Sbm holder cty2}
		I_2 = \int_{0}^{t_1} \int_0^\infty \int_{t_1}^{t_2} \dpd{S}{t}(y,\tau-s) \dif \tau \ b(y,s)m(\dif y,s)\dif s.
	\end{equation}
	Since $\dpd{S}{t} = \dfrac{\sigma^2}{2}\dpd[2]{S}{x}$, Lemma \ref{lem:m_n} implies
	\begin{equation} \label{eq:int Sbm holder cty2-1}
		\abs{I_2} \leq C(\sigma)\enVert{b}_\infty \int_0^{t_1}\int_{t_1}^{t_2} (\tau-s)^{-3/2} \dif \tau\dif s.
	\end{equation}
	By Fubini's Theorem,
	\begin{equation}\label{eq:int Sbm holder cty2-2}
		\int_0^{t_1}\int_{t_1}^{t_2} (\tau-s)^{-3/2} \dif \tau\dif s
		= 2\int_{t_1}^{t_2} \del{(\tau-t_1)^{-1/2} - \tau^{-1/2}}\dif \tau
		\leq 4(t_2-t_1)^{1/2}.
	\end{equation}
	Combining \eqref{eq:int Sbm holder cty2-1} and \eqref{eq:int Sbm holder cty2-2}, we get
	\begin{equation}
		\label{eq:eta holder2}
		\abs{\eta_2(t_1) - \eta_2(t_2)} \leq C(\sigma)\enVert{b}_\infty\abs{t_1-t_2}^{1/2}.
	\end{equation}
	Equation \eqref{eq:eta holder} follows from combining \eqref{eq:eta holder1} and \eqref{eq:eta holder2}.
\end{proof}

\begin{proof}[Proof of Lemma \ref{lem:Malpha fp bound}]
	We start from Equation \eqref{eq:fp formula} and multiply by $\phi_n(x)$, which is defined in \eqref{eq:phin}.
	Then integrate and use Lemma \ref{lem:Malpha bound} to get
	\begin{multline} \label{eq:phin m fp}
		\int_{\s{D}} \phi_n(x)m(\dif x,t) \leq
		C(\alpha)\int_{\s{D}} \abs{x}^{-\alpha}m_0(\dif x) \\
		+ \enVert{b}_\infty \int_0^t\int_{\s{D}}\abs{\int_0^\infty \phi_n(x)\del{\dpd{S}{x}(x-y,t-s) + \dpd{S}{x}(x+y,t-s)}\dif x} \ m(\dif y,s)\dif s.
	\end{multline}
	Let $j = \floor{\alpha}$.
	Integrating by parts $j$ times as in the proof of Lemma \ref{lem:Malpha bound}, we get
	\begin{multline} 
		\int_0^{y/2} \phi_n(x)\dpd{S}{x}(x \pm y,t-s)\dif x
		= \sum_{i=1}^{j}(-1)^{i-1} \Phi_n^{(i)}(y/2)\dpd[i]{S}{x}(y/2 \pm y,t-s)
		\\
		+ (-1)^j \int_0^{y/2} \Phi_n^{(j)}(x)\dpd[j+1]{S}{x}(x \pm y,t-s)\dif x.
	\end{multline}
	Using Lemma \ref{lem:m_n} and Equation \eqref{eq:Phin estimate}, we deduce
	\begin{multline} \label{eq:phindSdx est1}
		\abs{\int_0^{y/2} \phi_n(x)\dpd{S}{x}(x \pm y,t-s)\dif x}
		\leq \sum_{i=1}^{j} C(i,\alpha) \abs{y}^{i-\alpha}\abs{y}^{-i}\sigma^{-1}(t-s)^{-1/2}
		\\
		+  C(j,\alpha)\int_0^{y/2} \abs{x}^{j-\alpha}\abs{x \pm y}^{-j}\sigma^{-1}(t-s)^{-1/2}\dif x.
	\end{multline}
	For $0 \leq x \leq y/2$ we have $\abs{x \pm y} \geq y/2$, and thus \eqref{eq:phindSdx est1} yields
	\begin{equation} \label{eq:phindSdx est2}
		\abs{\int_0^{y/2} \phi_n(x)\dpd{S}{x}(x \pm y,t-s)\dif x}
		\leq 
		C(j,\alpha)\sigma^{-1}\abs{y}^{-\alpha}(t-s)^{-1/2}.
	\end{equation}
	On the other hand, using Lemma \ref{lem:m_n} it follows that $\int_0^\infty t^{1/2}\abs{\pd{S}{x}(x,t)}\dif x \leq C$ for all $t$, and thus
	\begin{equation} \label{eq:phindSdx est3}
		\abs{\int_{y/2}^\infty \phi_n(x)\dpd{S}{x}(x \pm y,t-s)\dif x}
		\leq
		C\abs{y}^{-\alpha}(t-s)^{-1/2}.
	\end{equation}
	Combining \eqref{eq:phindSdx est2} and \eqref{eq:phindSdx est3} into \eqref{eq:phin m fp}, then letting $n \to \infty$, we derive
	\begin{equation} \label{eq:phin m fp2}
		\int_{\s{D}} \abs{x}^{-\alpha}m(\dif x,t) \leq
		C(\alpha)\int_{\s{D}} \abs{x}^{-\alpha}m_0(\dif x) 
		+ C(\alpha,\sigma)\enVert{b}_\infty \int_0^t\int_{\s{D}}\abs{y}^{-\alpha} m(\dif y,s)(t-s)^{-1/2}\dif s.
	\end{equation}
	For $\lambda > 0$ let
	\begin{equation}
		f_\lambda(t) = e^{-\lambda t}\int_{\s{D}} \abs{x}^{-\alpha}m(\dif x,t).
	\end{equation}
	Multiply \eqref{eq:phin m fp2} by $e^{-\lambda t}$ to derive
	\begin{equation}
		\label{eq:phin m fp3}
		\begin{split}
			f_\lambda(t) &\leq C(\alpha)f_\lambda(0) + C(\alpha,\sigma)\enVert{b}_\infty \int_0^t e^{-\lambda(t-s)}(t-s)^{-1/2}f_\lambda(s)\dif s\\
			&\leq C(\alpha)f_\lambda(0) + C(\alpha,\sigma)\enVert{b}_\infty \lambda^{-1/2} \sup_{\tau \geq 0} f_\lambda(\tau)
		\end{split}
	\end{equation}
	where by a change of variables we have computed
	\begin{equation*}
		\int_0^t e^{-\lambda(t-s)}(t-s)^{-1/2}\dif s
		=
		\lambda^{-1/2}\int_0^{\lambda t} e^{- s}s^{-1/2}\dif s
		\leq \lambda^{-1/2}\del{\int_0^1 s^{-1/2} \dif s + \int_1^\infty e^{-s}\dif s} \leq 3\lambda^{-1/2}.
	\end{equation*}
	(As usual, the value of $C(\alpha,\sigma)$ might have changed from line to line.)
	Let $\lambda = \del{2C(\alpha,\sigma)\enVert{b}_\infty}^2$.
	Take the supremum in \eqref{eq:phin m fp3} to deduce
	\begin{equation}
		\label{eq:phin m fp4}
		\sup_{t \geq 0} f_\lambda(t) \leq C(\alpha)f_\lambda(0) + \frac{1}{2}\sup_{t \geq 0} f_\lambda(t)
		\quad
		\Rightarrow
		\sup_{t \geq 0} f_\lambda(t) \leq 2C(\alpha)f_\lambda(0).
	\end{equation}
	Equation \eqref{eq:phin m fp4} implies \eqref{eq:inverse moment bound fp}, as desired.
\end{proof}

\section{Proofs of Results from Section \ref{sec:fwdbckwd}}
\label{ap:forward-backward proofs}

\subsection{Proofs of Results from Section \ref{sec:hamiltonian}}
\label{ap:hamiltonian proofs}

We will actually show that all of the results of this section hold on a larger domain.
Set $p_\infty := \lim_{q \to \infty} P(q)$.
Note that $p_\infty < 0$ because there exists a finite saturation point (Assumption \ref{as:P}).
Recall that the profit function $\pi$ is defined as
\begin{equation*}
	\pi(\epsilon,q,Q,a) = \begin{cases}
		q\del{P(\epsilon Q+q)-a} &\text{if}~q > 0,\\
		0 &\text{if}~q = 0.
	\end{cases}
\end{equation*}
In the following the domain of $\pi$ is defined to be $\intco{0,\infty}^3 \times (p_\infty,\infty)$.
Thus the domain of $H(\epsilon,Q,a) := \sup_{q \geq 0} \pi(\epsilon,q,Q,a)$ is $\intco{0,\infty}^2 \times (p_\infty,\infty)$.
All the statements about the regularity of $H$ hold on this larger domain.
This remark will be useful in Lemma \ref{lem:HJ existence} below.
\begin{proof}[Proof of Lemma \ref{lem:opt quant}]
	We first compute
	\begin{equation}
		\dpd{\pi}{q}(\epsilon,q,Q,a) := qP'(\epsilon Q+q) + P(\epsilon Q+q)-a
	\end{equation}
	and
	\begin{equation}
		\dpd[2]{\pi}{q}(\epsilon,q,Q,a) = qP''(\epsilon Q+q) + 2P'(\epsilon Q  +q) = -\del{q\frac{\rho(\epsilon Q + q)}{\epsilon Q + q}-2}P'(\epsilon Q+q).
	\end{equation}
	By Assumption \ref{as:prudence} we deduce
	\begin{equation} \label{eq:2nd der pi}
		\dpd[2]{\pi}{q}(\epsilon,q,Q,a) \leq -\del{\bar \rho - 2}P'(\epsilon Q+q)  < 0,
	\end{equation}
	i.e.~$\pi$ is strictly concave with respect to $q$.
	On the other hand, since $P' \leq 0$ we also have
	\begin{equation} 
		\limsup_{q \to \infty} \dpd{\pi}{q}(\epsilon,q,Q,a) \leq \lim_{q \to \infty} P(\epsilon Q+q) - a = p_\infty - a < 0.
	\end{equation}
	Thus if $\dpd{\pi}{q}(\epsilon,0,Q,a) = P(\epsilon Q) - a > 0$ there must exist a unique $q^* > 0$ such that $\dpd{\pi}{q}(\epsilon,q^*,Q,a) = 0$, and hence $q^*$ maximizes $\pi(\epsilon,\cdot,Q,a)$.
	We also compute
	\begin{equation}
		\begin{aligned}
			\dmpd{\pi}{Q}{q}(\epsilon,q,Q,a) &= \epsilon qP''(\epsilon Q + q) + \epsilon P'(\epsilon Q + q)\\
			&= -\epsilon\del{q\frac{\rho(\epsilon Q + q)}{\epsilon Q + q}-1}P'(\epsilon Q+q), \\ \dmpd{\pi}{a}{q}(\epsilon,q,Q,a) &= -1,\\
			\dmpd{\pi}{\epsilon}{q}(\epsilon,q,Q,a) &= QqP''(\epsilon Q+q) + QP'(\epsilon Q+q)\\
			&= -Q\del{q\frac{\rho(\epsilon Q + q)}{\epsilon Q + q}-1}P'(\epsilon Q+q)
		\end{aligned}
	\end{equation}
	By the implicit function theorem, we deduce that $q^*$ is differentiable function of $(\epsilon,Q,a)$ in the region where $P(\epsilon Q) > a$, with
	\begin{equation} \label{eq:dq}
		\begin{aligned}
			\dpd{q^*}{Q} &= - \frac{\md{\pi}{2}{Q}{}{q}{}(\epsilon,q^*,Q,a)}{\pd[2]{\pi}{q}(\epsilon,q^*,Q,a)} = -\epsilon\del{1-\del{2-\frac{q^*\rho(\epsilon Q+q^*)}{\epsilon Q + q^*}}^{-1}},\\ 
			\dpd{q^*}{a} &= \frac{1}{\pd[2]{\pi}{q}(\epsilon,q^*,Q,a)} < 0,\\
			\dpd{q^*}{\epsilon} &= - \frac{\md{\pi}{2}{\epsilon}{}{q}{}(\epsilon,q^*,Q,a)}{\pd[2]{\pi}{q}(\epsilon,q^*,Q,a)} = -Q\del{1-\del{2-\frac{q^*\rho(\epsilon Q+q^*)}{\epsilon Q + q^*}}^{-1}}
		\end{aligned}
	\end{equation}
	Note that \eqref{eq:qstarQ leq} follows from \eqref{eq:dq}.
	
	In this region we also compute
	\begin{equation} \label{eq:dHde}
		\dpd{H}{\epsilon} = \dpd{q^*}{\epsilon}\del{P(\epsilon Q + q^*) - a} + q^*P'(\epsilon Q + q^*)\del{Q + \dpd{q^*}{\epsilon}} = Q q^*P'(\epsilon Q + q^*).
	\end{equation}
	\begin{equation} \label{eq:dHdQ}
		\dpd{H}{Q} = \dpd{q^*}{Q}\del{P(\epsilon Q + q^*) - a} + q^*P'(\epsilon Q + q^*)\del{\epsilon + \dpd{q^*}{Q}} = \epsilon q^*P'(\epsilon Q + q^*).
	\end{equation}
	and
	\begin{equation} \label{eq:dHda}
		\dpd{H}{a} = \dpd{q^*}{a}\del{P(\epsilon Q + q^*) - a} + q^*\del{P'(\epsilon Q + q^*)\dpd{q^*}{a} - 1} = -q^*.
	\end{equation}
	On the other hand, if $P(\epsilon Q) \leq a$ it follows that the unique maximizer is $q^* = 0$.
	Because $P$ is continuous and monotone decreasing, the interior of this region is the set where $P(\epsilon Q) < a$, while its boundary is where $P(\epsilon Q) = a$.
	It remains to show that as $(\epsilon,Q,a)$ approaches this boundary set, the derivative of $q^*$ remains bounded.
	By \eqref{eq:dq} it is enough to show that $\pd[2]{\pi}{q}(\epsilon,q^*(Q,a),Q,a)$ remains bounded away from zero.
	For this we observe that as $(\epsilon,Q,a)$ approaches the set where $P(\epsilon Q) = a$, $q^*(\epsilon,Q,a) \to 0$ and thus $\pd[2]{\pi}{q}(\epsilon,q^*,Q,a) \to 2P'(\epsilon Q)$, which is bounded away from zero for bounded values of $Q$.
\end{proof}

\begin{proof}[Proof of Corollary \ref{cor:smoothness}]
	For $(\epsilon,Q,a) \in [0,\bar \epsilon] \times [0,\bar Q] \times [0,\bar a]$ we have that $a \leq \bar a < P(\bar \epsilon \bar Q) \leq P(\epsilon Q)$, since $P$ is decreasing.
	By by differentiating \eqref{eq:dHde}, \eqref{eq:dHdQ}, and \eqref{eq:dHda} in the proof of Lemma \ref{lem:opt quant}, and using \eqref{eq:dq}, we see that $H$ is $n$ times continuously differentiable in this region.
	These derivatives are Lipschitz on this domain because $P^{(n)}$ is locally Lipschitz by Assumption \ref{as:P}.
	In particular,
	\begin{equation}
		\dpd[2]{H}{a}(\epsilon,Q,a) = -\dpd{q^*}{a}(\epsilon,Q,a) > 0.
	\end{equation}
	The claim follows from compactness of the region.
\end{proof}

\begin{proof}[Proof of Corollary \ref{cor:dHdQ}]
	From Equation \eqref{eq:dHdQ} and the first-order condition for optimality, using the fact that $P' < 0$, we have
	\begin{equation}
		\abs{\dpd{H}{Q}} = -\epsilon q^* P'(\epsilon Q + q^*) = \epsilon (P(\epsilon Q + q^*) - a),
	\end{equation}
	from which the first estimate in \eqref{eq:dHdQ estimate} follows.
	The second estimate follows from \eqref{eq:qstarQ leq} and \eqref{eq:dHdQ}.
\end{proof}

\begin{proof}[Proof of Lemma \ref{lem:market clearing Cournot}]
	Let $f(Q) = Q - \int_{\s{D}} q^*(\epsilon,Q,\phi(x))\dif m(x)$.
	We claim that $f(Q^*) = 0$ for a unique $Q^* \geq 0$.
	Note that $f(0) \leq 0$ because $q^* \geq 0$.
	By Lemma \ref{lem:opt quant} and Assumption \ref{as:prudence} we have
	\begin{equation}
		f'(Q) = 1 - \int_{\s{D}} \pd{q^*}{Q}(\epsilon,Q,\phi(x))\dif m(x) \geq 1 - \epsilon\frac{\bar \rho-1}{2-\bar \rho}\int_{\s{D}} \dif m(x) \geq \frac{2+\epsilon - (1+\epsilon)\bar \rho}{2-\bar \rho} > 0
	\end{equation}
	if $\bar \rho \geq 1$; otherwise we get simply $f'(Q) \geq 1$.
	The claim follows, and we deduce \eqref{eq:market clearing Cournot}.
	To derive estimate \eqref{eq:Q upper bound}, we use the lower bound on $f'$ to deduce
	\begin{equation} \label{eq:Q upper 1}
		Q^* \leq c(\bar \rho,\epsilon)\del{f(Q^*) - f(0)}
		= c(\bar \rho,\epsilon)\int_{\s{D}} q^*(\epsilon,0,\phi(x))\dif m(x).
	\end{equation}
	Now because $\pi(\epsilon,q,0,a) = \pi(0,q,0,a)$ for all $\epsilon,a$, it follows that $q^*(\epsilon,0,\phi(x)) = q^*(0,0,\phi(x))$.
	Then, since $q^*$ is decreasing in the last variable and $\int_{\s{D}} \dif m(x) \leq 1$, we use \eqref{eq:Q upper 1} to deduce \eqref{eq:Q upper bound}.
	
	We now prove \eqref{eq:Q1-Q2}.
	Without loss of generality we will assume $Q_1^* \geq Q_2^*$.
	First, observe that
	\begin{equation} \label{eq:Q lips 1}
		Q_1^*-Q_2^* \leq c(\bar \rho,\epsilon)\del{\int_{\s{D}} q^*(\epsilon_1,Q_2^*,\phi_1(x))\dif m_1(x)
			-\int_{\s{D}} q^*(\epsilon_2,Q_2^*,\phi_2(x))\dif m_2(x)}.
	\end{equation}
	To see this, note that \eqref{eq:qstarQ leq} implies
	\begin{multline}
		Q_1^*-Q_2^* = \int_{\s{D}} q^*(\epsilon_1,Q_1^*,\phi_1(x))\dif m_1(x)
		- \int_{\s{D}} q^*(\epsilon_2,Q_2^*,\phi_2(x))\dif m_2(x)\\
		\leq -\epsilon_1\frac{1-\bar \rho}{2-\bar \rho}(Q_1^*-Q_2^*)\int_{\s{D}} \dif m_1(x)
		+ \int_{\s{D}} q^*(\epsilon_1,Q_2^*,\phi_1(x))\dif m_1(x)
		- \int_{\s{D}} q^*(\epsilon_2,Q_2^*,\phi_2(x))\dif m_2(x).
	\end{multline}
	Then one obtains \eqref{eq:Q lips 1} by rearranging and using the fact that $\int_{\s{D}} \dif m_1(x) \leq 1$ and $c(\bar \rho,\epsilon)$ is increasing in $\epsilon$.
	Next, appealing to \eqref{eq:Q upper bound} and the fact that $q^*$ is locally Lipschitz, recalling once more that $\int_{\s{D}} \dif m_1(x) \leq 1$, \eqref{eq:Q lips 1} becomes
	\begin{equation} \label{eq:Q lips 2}
		\begin{split}
			Q_1^*-Q_2^* &\leq C\del{\abs{\epsilon_1-\epsilon_2} + \int_{\s{D}} \abs{\phi_1(x)-\phi_2(x)}\dif m_1(x)} + C\int_{\s{D}} q^*(\epsilon_2,Q_2^*,\phi_2(x))\dif (m_1-m_2)(x)\\
			&\leq C\del{\abs{\epsilon_1-\epsilon_2} + \enVert{\phi_1-\phi_2}_\infty} + C\sup_x\abs{\dod{q^*}{a}(\epsilon_2,Q_2^*,\phi_2(x))\dod{\phi_2}{x}(x)}{\bf d}_1(m_1,m_2)\\
			&\quad \quad +
			q^*(\epsilon_2,Q_2^*,\phi_2(0))\int_{\s{D}} \dif (m_1-m_2)(x),
		\end{split}
	\end{equation}
	which implies \eqref{eq:Q1-Q2}.
\end{proof}

\begin{proof}[Proof of Corollary \ref{cor:q* a priori bound}]
	We use Lemma \ref{lem:opt quant} to get
	\begin{equation} \label{eq:q* a priori bound1}
		q^*(\epsilon,Q^*,\phi(x)) \leq q^*(\epsilon,0,0) + \epsilon \frac{\bar \rho -1}{2-\bar \rho}Q^*.
	\end{equation}
	We recall that $q^*(\epsilon,0,0) = q^*(0,0,0)$.
	To derive \eqref{eq:q* a priori bound}, it suffices to plug \eqref{eq:Q upper bound} into \eqref{eq:q* a priori bound1} and use the definition of $c(\bar \rho,\epsilon)$.
\end{proof}	

\subsection{Proofs of Results from Section \ref{sec:estimates on hj}}
\label{ap:hj proofs}

\begin{proof}[Proof of Lemma \ref{lem: a priori HJ}]
	First let $v = e^{-rt}u$.
	Then $v$ satisfies
	\begin{equation} \label{eq:HJint}
		\dpd{v}{t} + \dfrac{\sigma^2}{2}\dpd[2]{v}{x} + e^{-rt}H\del{\epsilon(t),Q^*(t),e^{rt}\dpd{v}{x}} = 0, \ x \in \s{D}, \ t > 0.
	\end{equation}
	Using the fact that $H \geq 0$ and $v(0,t) = 0$, the maximum principle (see \cite[Proposition 2.1]{souplet2006global}) implies
	\begin{equation} \label{eq:ugeq0}
		\min_{x \in \s{D}, 0 \leq t \leq T} v(x,t) = \min_{x \in \s{D}} v(x,T)
		= e^{-rT} \min_{x \in \s{D}} u(x,T) = 0
		\quad
		\Rightarrow v \geq 0
		\Rightarrow u \geq 0.
	\end{equation}
	It also follows that $u(0,t) = \min u$ and so $u_x(0,t) \geq 0$.
	
	We now use the fact that $H$ is decreasing in all variables to deduce
	\begin{equation} \label{eq:H bound}
		0 \leq H\del{\epsilon(t),Q^*(t),e^{rt}\dpd{v}{x}} \leq H(0,0,0)
	\end{equation}
	and thus
	\begin{equation} \label{eq:HJineq}
		-\dpd{v}{t} - \dfrac{\sigma^2}{2}\dpd[2]{v}{x} \leq e^{-rt}H(0,0,0).
	\end{equation}
	Set $\tilde v(x,t) = v(x,t) - \int_t^T e^{-rs}H(0,0,0)\dif s = v(x,t) + \frac{1}{r}H(0,0,0)\del{e^{-rT}-e^{-rt}}$.
	It follows that
	\begin{equation} \label{eq:HJineq2}
		-\dpd{\tilde v}{t} - \dfrac{\sigma^2}{2}\dpd[2]{\tilde v}{x} \leq 0
	\end{equation}
	and thus
	\begin{equation}
		\max_{x \in \s{D}, 0 \leq t \leq T} \tilde v(x,t) = \max_{t =T \ \text{or} \ x = 0} \tilde v(x,t)
		\leq c_1e^{-rT}
	\end{equation}
	since $\tilde v(x,T) = e^{-rT}u_T(x) \leq c_1e^{-rT}$ and $\tilde v(0,t) = \frac{1}{r}H(0,0,0)\del{e^{-rT}-e^{-rt}} \leq 0$.
	Together with \eqref{eq:ugeq0} we deduce that
	\begin{equation} \label{eq:uleq}
		0 \leq v(x,t) \leq \del{\frac{1}{r}H(0,0,0) + c_1}e^{-rt} 
		\ \Rightarrow \
		0 \leq u(x,t) \leq \frac{1}{r}H(0,0,0) + c_1.
	\end{equation}
	
	To get an estimate on $u_x$, we now use a Bernstein type argument, cf.~\cite[Section VI.3]{ladyzhenskaia1968linear}.
	Notice that
	\begin{equation} \label{eq:HJineq u}
		-\dpd{u}{t} - \dfrac{\sigma^2}{2}\dpd[2]{u}{x} \leq H(0,0,0).
	\end{equation}
	Set $\tilde u(x,t) = u(x,t) + M_\lambda e^{-\lambda x}$, where $M_\lambda  > 0$ and $\lambda > 0$ are defined below in \eqref{eq:M def1} and \eqref{eq:lambda}.		
	The constants $M_\lambda $ and $\lambda$ have to be chosen so that, for all $t \leq T$ and all $x \in [0,\ell]$ for $\ell > 0$ to be specified later, we have
	\begin{equation} \label{eq:M conditions}
		\begin{split}
			H(0,0,0) - \frac{\sigma^2}{2}\lambda^2 M_\lambda e^{-\lambda x} &\leq 0,\\
			c_3 &\leq M_\lambda \lambda e^{-\lambda x},\\
			\frac{1}{r}H(0,0,0) + c_1 + M_\lambda e^{-\lambda x} &\leq M_\lambda .
		\end{split}
	\end{equation}
	Then one can check that
	\begin{equation} \label{eq:HJineq u tilde}
		-\dpd{\tilde u}{t} - \dfrac{\sigma^2}{2}\dpd[2]{\tilde u}{x} \leq 0 \ \text{in} \ (0,\ell) \times (0,T),
	\end{equation}
	$\tilde u(x,T) \leq M_\lambda $ for all $x \in [0,\ell]$ (using the fact that $\max u_T' \leq c_3$), $\tilde u(\ell,t) \leq M_\lambda $ (using \eqref{eq:uleq}), and $\tilde u(0,t) = M_\lambda $.
	By the maximum principle, it follows that $\tilde u(x,t) \leq M_\lambda $ for all $x \in [0,\ell], t \in [0,T]$.
	This means $\tilde u(0,t) = \max_{0 \leq x \leq \ell} \tilde u(x,t)$, which implies $\tilde u_x(0,t) \leq 0$ and thus $u_x(0,t) \leq M_\lambda \lambda$.
	Finally, we can take the derivative of Equation \eqref{eq:HJ} to see that the maximum principle applies to $u_x$, and thus
	\begin{equation} \label{eq:uxleq}
		\max_{x \in \s{D}, 0 \leq t \leq T} u_x(x,t) = \max_{t =T \ \text{or} \ x = 0} u_x(x,t) \leq \max\cbr{M_\lambda \lambda, c_3}.
	\end{equation}
	To satisfy \eqref{eq:M conditions}, we choose
	\begin{equation} \label{eq:M def1}
		M_\lambda  = \max\cbr{\frac{2}{\sigma^2\lambda^2}H(0,0,0)e^{\lambda\ell},
			\frac{c_3}{\lambda}e^{\lambda\ell},
			\frac{1}{1-e^{-\lambda\ell}}\del{\frac{1}{r}H(0,0,0) + c_1}}.
	\end{equation}
	If we set $J = \max\cbr{\frac{2}{\sigma^2\lambda^2}H(0,0,0),
		\frac{c_3}{\lambda}}$, then \eqref{eq:M def1} becomes
	\begin{equation} \label{eq:M def2}
		M_\lambda  = \max\cbr{Je^{\lambda\ell},
			\frac{1}{1-e^{-\lambda\ell}}\del{\frac{1}{r}H(0,0,0) + c_1}}.
	\end{equation}
	To minimize the value of $M_\lambda \lambda$, we first choose the constant $\ell$ so as to minimize the maximum appearing in \eqref{eq:M def2}; it suffices to choose it so that the two maximands are equal, because the first is increasing in $\ell$ while the second is decreasing.
	This is achieved by setting
	\begin{equation}
		\ell = \frac{1}{\lambda}\ln\del{1 + \frac{1}{rJ}H(0,0,0) + \frac{c_1}{J}}
		\Rightarrow
		M_\lambda  = Je^{\lambda\ell} = J + \frac{1}{r}H(0,0,0) + c_1.
	\end{equation}
	We therefore have
	\begin{equation}
		M_\lambda \lambda = \max\cbr{\frac{2}{\sigma^2\lambda}H(0,0,0),c_3} + \frac{\lambda}{r}H(0,0,0) + \lambda c_1.
	\end{equation}
	The minimum possible value of the right-hand side is attained by setting
	\begin{equation}
		\label{eq:lambda}
		\lambda = \min\cbr{\frac{\sqrt{2rH(0,0,0)}}{\sigma\sqrt{H(0,0,0) + rc_1}},\frac{2}{\sigma^2 c_3}H(0,0,0)},
	\end{equation}
	and its minimum value is given by $M_\lambda \lambda = M$ where $M$ is defined in \eqref{eq:max ux}.

	Put together \eqref{eq:uleq}, and \eqref{eq:uxleq} to get \eqref{eq:u ux bounded}.		
\end{proof}

\subsection{Proof of Result from Section \ref{sec:estimates on coupling}}
\label{ap:coupling proofs}

\begin{proof}[Proof of Lemma \ref{lem:Q(t) Holder}]
	Estimate \eqref{eq:Q(t) upper bound} follows from Lemmas \ref{lem:market clearing Cournot}, \ref{lem:opt quant}, \ref{lem: a priori HJ}, and \ref{lem:fp}.
	Note that a direct application of Lemma \ref{lem:market clearing Cournot} would put the constant $c(\bar \rho,\epsilon(t))$ in place of $c(\bar \rho,\epsilon(0))$; however, $c(\bar \rho,\epsilon)$ defined in \eqref{eq:Q upper bound} is an increasing function of $\epsilon$, and since Assumption \ref{as:epsilon(t)} implies $\epsilon(t) \leq \epsilon(0)$, we have replaced $c(\bar \rho,\epsilon(t))$ with $c(\bar \rho,\epsilon(0))$ to get an upper bound that is uniform in time.
	
	Now we turn to estimate \eqref{eq:Q(t) Holder}.
	By Lemma \ref{lem:market clearing Cournot},
	there exists a constant $C = C(\epsilon(0),\bar \rho,M)$ such that
	\begin{multline*}
		\abs{Q^*(t_1)-Q^*(t_2)} \leq C\del{\abs{\epsilon(t_1)-\epsilon(t_2)} + \enVert{\dpd{u}{x}(\cdot,t_1)-\dpd{u}{x}(\cdot,t_2)}_\infty}\\
		+ C\del{\enVert{\dpd[2]{u}{x}}_\infty{\bf d}_1(m(t_1),m(t_2)) + \abs{\int_{\s{D}} \dif (m(t_1)-m(t_2))(x)}}
	\end{multline*}
	Now suppose $m_0 \in \sr{M}_{\alpha}$.
	Appealing to Lemma \ref{lem:fp} and also Assumption \ref{as:epsilon(t)}, we have
	\begin{multline*}
		\abs{Q^*(t_1)-Q^*(t_2)} \leq C\del{ \abs{t_1-t_2} + \enVert{\dpd{u}{x}}_{\s{C}^{\alpha,\alpha/2}}\abs{t_1-t_2}^{\alpha/2}}\\
		+ C\del{\enVert{\dpd[2]{u}{x}}_\infty\del{\enVert{b}_\infty + \sigma}\abs{t_1-t_2}^{1/2} + C(\sigma)\del{\enVert{m_0}_{\sr{M}_{\alpha/2}} + \enVert{b}_\infty}\abs{t_1-t_2}^{\alpha/2}}
	\end{multline*}
	for any $\abs{t_1-t_2} \leq 1$.
	Here $b = \pd{H}{a}\del{\epsilon,Q^*,\pd{u}{x}}$.
	By Lemmas \ref{lem:opt quant} and \ref{lem: a priori HJ} together with \eqref{eq:Q(t) upper bound}, we deduce there exists $C = C(\bar \rho,\epsilon(0),M)$ such that $\enVert{b}_\infty \leq C.$
	We deduce that there exists $C = C(\bar \rho,\epsilon(0),M,\sigma,\enVert{m_0}_{\sr{M}_{\alpha/2}})$ such that
	\begin{equation*}
		\abs{Q^*(t_1)-Q^*(t_2)} \leq C\del{\enVert{\dpd{u}{x}}_{\s{C}^{\alpha,\alpha/2}} + \enVert{\dpd[2]{u}{x}}_{\infty} + 1}\abs{t_1-t_2}^{\alpha/2},
		\quad \abs{t_1-t_2} \leq 1,
	\end{equation*}
	and since $Q^*$ is bounded according to \eqref{eq:Q(t) upper bound}, Equation \eqref{eq:Q(t) Holder} follows.
	\begin{comment}
	To prove \eqref{eq:Q(t) Holder delta} for a given $\delta > 0$, note that by interpolation,
	\begin{equation} \label{eq:interpolation1}
	\enVert{\dpd{u}{x}}_{\s{C}^{\alpha,\alpha/2}} + \enVert{\dpd[2]{u}{x}}_{\infty}
	\leq \delta\enVert{u}_{\s{C}^{2+\alpha,1+\alpha/2}} + C(\delta,\alpha)\enVert{u}_\infty,
	\end{equation}
	which we combine with Lemma \ref{lem: a priori HJ} to finish the proof.
	\end{comment}
\end{proof}

\subsection{Proofs of Results from Section \ref{sec:parabolic estimates}}
\label{ap:parabolic proofs}
\begin{proof}[Proof of Lemma \ref{lem:C2alpha estimates}]
	We begin by taking $u_0 = 0$.
	First we let $g(x,t) = e^{rt}f(x,t)$ and consider 
	\begin{equation} \label{eq:parabolic + int factor}
		\dpd{v}{t} - \frac{\sigma^2}{2}\dpd[2]{v}{x} = g, \ \forall x \in \s{D}, t > 0; v(0,t) = 0 \ \forall t > 0; \ v(x,0) = 0 \ \forall x \in \s{D}.
	\end{equation}
	By \cite[Theorem IV.6.1]{ladyzhenskaia1968linear}, \eqref{eq:parabolic + int factor} is uniquely solvable in $\s{C}^{2+\alpha,1+\alpha/2}(\overline{\s{D}} \times [0,T])$ for arbitrary $T > 0$.
	Also, by the maximum principle, we have
	\begin{equation}
		\abs{v(x,t)} \leq \frac{1}{r}e^{rt}\enVert{f}_{\s{C}^{0,0}} \ \forall x \in \overline{\s{D}}, t \in \intco{0,\infty}.
	\end{equation}
	To see this, first let $\tilde v(x,t) = v(x,t) - \frac{e^{rt}-1}{r}\enVert{f}_{\s{C}^{0,0}}$ and observe that
	\begin{equation}
		\dpd{\tilde v}{t} - \frac{\sigma^2}{2}\dpd[2]{\tilde v}{x} \leq 0, \ \tilde v(0,t) \leq 0, \ \tilde v(x,0) \leq 0.
	\end{equation}
	By the maximum principle, $\tilde v \leq 0$, which implies $v(x,t) \leq \frac{1}{r}e^{rt}\enVert{f}_{\s{C}^{0,0}}$.
	The opposite inequality is similarly proved.
	
	Now we let $u(x,t) = e^{-rt}v(x,t)$.
	Then $u$ satisfies \eqref{eq:parabolic + ru} and
	\begin{equation} \label{eq:max parabolic}
		\enVert{u}_{\s{C}^{0,0}(\overline{\s{D}} \times \intco{0,\infty})} \leq \frac{1}{r}\enVert{f}_{\s{C}^{0,0}}.
	\end{equation}
	Moreover, appealing again to \cite[Theorem IV.6.1]{ladyzhenskaia1968linear}, we have an estimate
	\begin{equation} \label{eq:holder parabolic}
		\intcc{\dpd[2]{u}{x}}_{\alpha,\alpha/2} + \intcc{\dpd{u}{t}}_{\alpha,\alpha/2} \leq C(\sigma)\del{\intcc{f}_{\alpha,\alpha/2} + r\intcc{u}_{\alpha,\alpha/2}},
	\end{equation}
	where $C(\sigma)$ does not depend on $T$.
	By interpolation, see \cite[Lemma II.3.2]{ladyzhenskaia1968linear}, we can find a constant $C(\alpha)$ such that for arbitrary $\delta > 0$  we have
	\begin{equation} \label{eq:interpolation parabolic}
		\intcc{u}_{\alpha,\alpha/2} \leq C(\alpha)\del{\delta^{-\alpha}
			\enVert{u}_{\s{C}^{0,0}(\overline{\s{D}})} + \delta^2\del{\intcc{\dpd[2]{u}{x}}_{\alpha,\alpha/2} + \intcc{\dpd{u}{t}}_{\alpha,\alpha/2}}}.
	\end{equation}
	Combining \eqref{eq:max parabolic}, \eqref{eq:holder parabolic}, and \eqref{eq:interpolation parabolic} with $\delta$ a sufficiently small multiple of $r^{-1/2}$, we deduce that \eqref{eq:global holder estimate1} holds for $u_0 = 0$.
	
	Now suppose $f = 0$ and let $u_0 \in \s{C}^{2+\alpha}(\overline{\s{D}})$ be given.
	Then appealing to \cite[Theorem IV.5.1]{ladyzhenskaia1968linear}, \eqref{eq:parabolic + ru} is uniquely solvable, and moreover by the potential estimates from \cite[Section IV.2]{ladyzhenskaia1968linear} we have
	\begin{equation} \label{eq:potential holder}
		\intcc{\dpd[2]{u}{x}}_{\alpha,\alpha/2} + \intcc{\dpd{u}{t}}_{\alpha,\alpha/2} \leq C(\sigma)\del{\enVert{u_0}_{\s{C}^{2+\alpha}(\overline{\s{D}})} + r\intcc{u}_{\alpha,\alpha/2}},
	\end{equation}
	where again $C(\sigma)$ does not depend on time.
	Using the maximum principle, we get $\enVert{u}_0 \leq \enVert{u_0}_0$.
	Arguing as before, we deduce \eqref{eq:global holder estimate1} for $f = 0$.
	
	The general case now follows from linearity.
	
	\subsection{Proofs of Results from Section \ref{sec:existence}} \label{ap:existence proofs}
	
	\begin{proof}[Proof of Lemma \ref{lem:uC2alpha}]
		Let $f = f(x,t) = H\del{\epsilon(t),Q^*(t),\dpd{u}{x}(x,t)}$.
		From Lemma \ref{lem:C2alpha estimates} we have
		\begin{equation} \label{eq:uC2alpha1}
			\enVert{u}_{\s{C}^{2+\alpha,1+\alpha/2}(\overline{\s{D}} \times \intcc{0,T})} \leq C(\sigma,\alpha)\del{\sbr{f}_{\alpha,\alpha/2} + r^{\frac{\alpha}{2}}\enVert{f}_0 + C_\alpha
				+ r^{1 + \frac{\alpha}{2}}c_1}.
		\end{equation}
		We now estimate $f$ in $\s{C}^{\alpha,\alpha/2}$.
		First, because $H$ is decreasing in all variables, we (again) deduce
		\begin{equation}
			0 \leq f(x,t) \leq H(0,0,0).
		\end{equation}
		Because $H$ is locally Lipschitz by Lemma \ref{lem:opt quant}, and because $\epsilon,Q^*,$ and $\pd{u}{x}$ are bounded with estimates given in Assumption \ref{as:epsilon(t)}, Lemma \ref{lem: a priori HJ} and Lemma \ref{lem:Q(t) Holder}, we have a constant $C = C(\bar \rho,\epsilon(0),\sigma,M,\alpha)$ such that
		\begin{equation*}
			\enVert{f}_{\s{C}^{\alpha,\alpha/2}} \leq C\del{1 + \enVert{Q^*}_{\s{C}^{\alpha/2}} + \enVert{\dpd{u}{x}}_{\s{C}^{\alpha,\alpha/2}}},
		\end{equation*}
		where $\enVert{\epsilon}_{\s{C}^{\alpha/2}}$ is also estimated using Assumption \ref{as:epsilon(t)}.
		Using Lemma \ref{lem:Q(t) Holder} and interpolation on H\"older spaces, we see that for an arbitrary $\delta > 0$, there exists  $C_\delta = C(\delta,\bar \rho,\epsilon(0),\sigma,M,\enVert{m_0}_{\sr{M}_{\alpha/2}},\alpha)$ such that
		\begin{equation*}
			\enVert{f}_{\s{C}^{\alpha,\alpha/2}} \leq \delta\enVert{u}_{\s{C}^{2+\alpha,1+\alpha/2}} + C_\delta.
		\end{equation*}
		Taking $\delta > 0$ small enough, \eqref{eq:uC2alpha1} becomes
		\begin{equation} \label{eq:uC2alpha2}
			\enVert{u}_{\s{C}^{2+\alpha,1+\alpha/2}(\overline{\s{D}} \times \intcc{0,T})} \leq C(\sigma,\alpha)\del{C_\delta + r^{\frac{\alpha}{2}}H(0,0,0) + C_\alpha
				+ r^{1 + \frac{\alpha}{2}}c_1},
		\end{equation}
		which proves \eqref{eq:uC2alpha}.
	\end{proof}

	Before getting to the proof of Theorem \ref{thm:exist mfg T}, we establish the following lemma:
	\begin{lemma} \label{lem:HJ existence}
		Let $Q^* \in \s{C}^{\alpha/2}\del{[0,T]; \intco{0,\infty}}$ be given, and let Assumptions \ref{as:epsilon(t)} and \ref{as:uT} hold.
		Then there exists a unique solution $u$ to the PDE
		\begin{equation} \label{eq:hj equation}
			\pd{u}{t} + \frac{\sigma^2}{2} \pd[2]{u}{x} + H\del{\epsilon(t),Q^*(t),\dpd{u}{x}} - r u = 0, \quad u(0,t) = 0, \quad u(x,T) = u_T(x).
		\end{equation}
		This solution satisfies $\dpd{u}{x} \geq 0$ and the a priori estimate
		\begin{equation} \label{eq:hj eq estimate}
			\enVert{u}_{\s{C}^{2+\alpha,1+\alpha/2}} \leq C\del{\sigma,r,\alpha,H(0,0,0)}\del{1 + \enVert{Q^*}_{\s{C}^{\alpha/2}} + \enVert{u_T}_{\s{C}^{2+\alpha}}}.
		\end{equation}
	\end{lemma}

	\begin{proof}
		As above we set $p_\infty := \lim_{q \to \infty} P(q) < 0$.
		Fix a $\s{C}^\infty$ function $\psi:\bb{R} \to (p_\infty,\infty)$ such that $p(a) = a$ for all $a \geq 0$.
		Define $X = \s{C}^{2,1}\del{\overline{\s{D}} \times [0,T]}$.
		Let $v \in X$ and $\lambda \in [0,1]$.
		By Lemma \ref{lem:C2alpha estimates}, Assumption \ref{as:uT}, and the local Lipschitz property of $H$, we get a unique solution $u$ to the equation
		\begin{equation} \label{eq:hj fixed point equation}
			\pd{u}{t} + \frac{\sigma^2}{2} \pd[2]{u}{x} + \lambda H\del{\epsilon(t),Q^*(t),\psi\del{\dpd{v}{x}}} - r u = 0, \quad u(0,t) = 0, \quad u(x,T) = \lambda u_T(x),
		\end{equation}
		and $u$ satisfies
		\begin{equation*}
			\enVert{u}_{\s{C}^{2+\alpha,1+\alpha/2}} \leq C\del{\enVert{Q^*}_\infty,\enVert{\dpd{v}{x}}_\infty}\del{1 + \enVert{Q^*}_{\s{C}^{\alpha/2}} + \enVert{\dpd{v}{x}}_{\s{C}^{\alpha,\alpha/2}}}.
		\end{equation*}
		This defines a map $\s{T}:X \times [0,1] \to X$.
		We claim that $\s{T}$ is continuous and compact.
		Suppose $\{v_n,\lambda_n\}$ is a bounded sequence in $X \times [0,1]$ and let $u_n = \s{T}(v_n,\lambda_n)$.
		Then $\{u_n\}$ is bounded in $\s{C}^{2+\alpha,1+\alpha/2}$, which is compactly embedded in $X$, so it has a subsequence that converges to some $u$ in $X$.
		To conclude that $\s{T}$ is both continuous and compact, it is enough to show that whenever $(v_n,\lambda_n) \to (v,\lambda)$ in $X \times [0,1]$, then $u = \s{T}(v,\lambda)$.
		But this can be deduced from plugging $v_n,\lambda_n$ into \eqref{eq:hj fixed point equation} in place of $v,\lambda$, then passing to the limit using the local Lipschitz property of $H$.
		
		Notice that $\s{T}(v,0) = 0$.
		To apply the Leray-Schauder fixed point theorem, it remains to find an a priori bound on solutions to the fixed point equation $u = \s{T}(u,\lambda)$.
		Note that for any such fixed point, $w = \dpd{u}{x}$ satisfies, in a weak sense,
		\begin{equation}
			\pd{w}{t} + \frac{\sigma^2}{2} \pd[2]{w}{x} + \lambda \dpd{H}{a}\del{\epsilon(t),Q^*(t),\psi\del{\dpd{u}{x}}}\psi'\del{\dpd{u}{x}}\dpd{w}{x} - r w = 0,  \quad w(x,T) = \lambda u_T'(x).
		\end{equation}
		Since $u_T' \geq 0$, by the maximum principle we deduce $\dpd{u}{x} = w \geq 0$.
		It follows that $u$ satisfies
		\begin{equation}
			\pd{u}{t} + \frac{\sigma^2}{2} \pd[2]{u}{x} + \lambda H\del{\epsilon(t),Q^*(t),\dpd{u}{x}} - r u = 0, \quad u(0,t) = 0, \quad u(x,T) = \lambda u_T(x),
		\end{equation}
		Lemma \ref{lem: a priori HJ} establishes an a priori bound on $u$; combined with Lemma \ref{lem:C2alpha estimates} and using interpolation, we deduce that \eqref{eq:hj eq estimate} holds for any $u$ satisfying $u = \s{T}(u,\lambda)$.
		By the Leray-Schauder fixed point theorem \cite[Theorem 11.6]{gilbarg2015elliptic}, there exists $u \in X$ such that $u = \s{T}(u,1)$, which means $u$ is a solution to \eqref{eq:hj equation}.
		Uniqueness follows from the maximum principle by standard arguments.
	\end{proof}

	\begin{proof}[Proof of Theorem \ref{thm:exist mfg T}]
		Set $X$ to be the set of all $(v,Q) \in \s{C}^{2,1}\del{\overline{\s{D}} \times [0,T]} \times \s{C}^{0}\del{[0,T]}$ such that $\dpd{v}{x} \geq 0$ and $Q \geq 0$, and define $\s{T}:X \times [0,1] \to X$ as follows.
		Let $(v,Q) \in X, \lambda \in [0,1]$.
		From Lemma \ref{lem:opt quant} we know that the function $\dpd{H}{a}\del{\epsilon(t),Q(t),\dpd{v}{x}}$ is bounded and continuous with
		\begin{equation*}
			\enVert{\dpd{H}{a}\del{\epsilon(t),Q(t),\dpd{v}{x}}}_\infty \leq C\del{\epsilon(0),\enVert{Q}_\infty,\enVert{\dpd{v}{x}}_\infty}
		\end{equation*}
		By Lemma \ref{lem:fp}, there exists a unique solution $m$ satisfying
		\begin{equation} \label{eq:fp fixed pt}
			\dpd{m}{t} - \dfrac{\sigma^2}{2}\dpd[2]{m}{x} + \dpd{}{x}\del{\lambda\dpd{H}{a}\del{\epsilon(t),Q(t),\dpd{v}{x}}m} = 0, \ m|_{x=0} = 0, \
			m|_{t=0} = m_0,
		\end{equation}
		and moreover we have H\"older estimates \eqref{eq:m holder in time} and \eqref{eq:eta holder}.
		Now by Lemma \ref{lem:market clearing Cournot} we can define $Q^*(t)$ by
		\begin{equation} \label{eq:Q^* fixed pt}
			Q^*(t) = -\int_0^\infty \lambda\dpd{H}{a}\del{\epsilon(t),Q^*(t),\dpd{v}{x}}\dif m(t),
		\end{equation}
		and combining \eqref{eq:Q upper bound}, \eqref{eq:Q1-Q2}, \eqref{eq:m holder in time} and \eqref{eq:eta holder}, we have
		\begin{equation} \label{eq:Q^*smoothing}
			\enVert{Q^*}_{\s{C}^{\alpha/2}\del{[0,T]}} \leq C\del{\epsilon(0),\enVert{Q}_\infty,\enVert{\dpd{v}{x}}_\infty,\enVert{m_0}_{\sr{M}_{\alpha/2}},\sigma,\alpha}\del{\enVert{v}_{\s{C}^{2,1}} + 1}.
		\end{equation}
		Setting $f(x,t) = \lambda H\del{\epsilon(t),Q^*(t),\dpd{v}{x}}$, we have, as in the proof of Lemma \ref{lem:uC2alpha},
		\begin{equation*}
			\enVert{f}_{\s{C}^{\alpha,\alpha/2}} \leq C\del{\enVert{Q^*}_\infty,\enVert{\dpd{v}{x}}_\infty}\del{1 + \enVert{Q^*}_{\s{C}^{\alpha/2}} + \enVert{\dpd{v}{x}}_{\s{C}^{\alpha,\alpha/2}}}.
		\end{equation*}
		Thus, by Lemma \ref{lem:HJ existence} there exists a unique solution $u$ of
		\begin{equation} \label{eq:HJ fixed pt}
			\dpd{u}{t} + \dfrac{\sigma^2}{2}\dpd[2]{u}{x} + \lambda H\del{\epsilon(t),Q^*(t),\dpd{u}{x}} - ru = 0,\
			u|_{x=0} = 0, \ u|_{t=T} = \lambda u_T
		\end{equation}
		satisfying \eqref{eq:hj eq estimate}, which in this case can be written
		\begin{equation}\label{eq:u smoothing}
			\enVert{u}_{\s{C}^{2+\alpha,1+\alpha/2}} \leq C\del{\epsilon(0),\enVert{Q}_\infty,\enVert{v}_{\s{C}^{2,1}},\enVert{m_0}_{\sr{M}_{\alpha/2}},\sigma,\alpha}.
		\end{equation}
		Then we set $\s{T}(v,Q,\lambda) = (u,Q^*) \in X$.
		We need to show that $\s{T}$ is continuous and compact.
		Suppose $\{(v_n,Q_n,\lambda_n)\}$ is a sequence in $X \times [0,1]$, and let $(u_n,Q^*_n) = \s{T}(v_n,Q_n,\lambda_n)$.
		Note that by \eqref{eq:Q^*smoothing} and \eqref{eq:u smoothing}, $(u_n,Q^*_n)$ must have a subsequence converging to $(u,Q^*) \in X$, because $\s{C}^{2+\alpha,1+\alpha/2} \times \s{C}^{\alpha/2}$ is compactly embedded in $\s{C}^{2,1} \times \s{C}^0$.
		We now show that if $(v_n,Q_n,\lambda_n) \to (v,Q,\lambda)$, then $(u,Q^*) = \s{T}(v,Q,\lambda)$.
		First let $m_n$ be the solution of \eqref{eq:fp fixed pt} corresponding to $(v_n,Q_n,\lambda_n)$.
		By Lemma \ref{lem:fp} we have that $m_n$ is uniformly H\"older in the ${\bf d}_1$ metric, hence by passing to a subsequence it converges to some $m$ in $\s{C}\del{[0,T];\s{M}_{1,+}}$.
		Since $\dpd{H}{a}$ is locally Lipschitz, we have that $\lambda_n\dpd{H}{a}\del{\epsilon(t),Q_n(t),\dpd{v_n}{x}} \to \lambda\dpd{H}{a}\del{\epsilon(t),Q(t),\dpd{v}{x}}$ uniformly.
		Combining these facts we deduce that $m$ is really the solution to \eqref{eq:fp fixed pt} and $Q^*$ is the solution of \eqref{eq:Q^* fixed pt}.
		Finally, we deduce that $u$ is really the solution to \eqref{eq:HJ fixed pt} by taking the corresponding equation for $u_n$ and passing to the limit.
		We have thus proved that $\s{T}$ is continuous and compact.
		
		It remains to show there exists a constant $C$ such that whenever $\s{T}(u,Q^*,\lambda) = (u,Q^*)$, then
		\begin{equation*}
			\enVert{(u,Q^*)}_X \leq C.
		\end{equation*}
		But this is a consequence of Lemmas \ref{lem:Q(t) Holder} and \ref{lem:uC2alpha}, since $\lambda H$ and $\lambda u_T$ satisfy all the same estimates as $H$ and $u_T$.
		Now we can apply the Leray-Schauder fixed point theorem, which says that there exists $(u,Q^*)$ such that $\s{T}(u,Q^*,1) = (u,Q^*)$.
		Letting $m$ now be defined by solving \eqref{eq:mfg T}(ii), we deduce that $(u,m)$ solves the system \eqref{eq:mfg T}.
		The regularity of this solution follows by once more appealing to Lemmas \ref{lem:fp} and \ref{lem:uC2alpha}.
	\end{proof}
\end{proof}

\section{Proof of the integral estimate used in Section \ref{sec:sensitivity}}

\label{ap:int estimate}

The following proof is more or less the same as that of \cite[Lemma 2.1]{graber2022parameter}.
We include it for completeness.
\begin{proof}[Proof of Lemma \ref{lem:int estimate}]
	Set $h(t) = Bf(t) + g(t)$, so that \eqref{eq:ineq} reads simply
	\begin{equation} \label{eq:ineq1}
		f(t_1) \leq Af(t_0) + \int_{t_0}^{t_1}(t_1-s)^{-1/2}h(s)\dif s \quad  \forall 0 \leq t_0 \leq t_1 \leq t_0 + \delta
	\end{equation}
	For arbitrary $t > 0$ let $n = \floor{\frac{t}{\delta}}$.
	Use \eqref{eq:ineq1} $n+1$ times to get
	\begin{equation} \label{eq:ineq2}
		f(t) \leq A^{n+1}f(0) + \sum_{j=0}^n A^j\int_{\del{t-(j+1)\delta}_+}^{t-j\delta}(t-j\delta-s)^{-1/2}h(s)\dif s,
	\end{equation}
	where $s_+ := \max\{s,0\}$.
	Note that
	\begin{equation*}
		t - (j+1)\delta < s \leq t-j\delta \Rightarrow
		j = \floor{\frac{t-s}{\delta}}
	\end{equation*}
	So we define $\phi(s) = \del{s - \floor{\frac{s}{\delta}}\delta}^{-1/2}$.		
	Then \eqref{eq:ineq2} implies
	\begin{equation} \label{eq:ineq3}
		f(t) \leq A^{\frac{t}{\delta}+1}f(0) + \sum_{j=0}^n  \int_{\del{t-(j+1)\delta}_+}^{t-j\delta}
		A^{\frac{t-s}{\delta}}\phi(t-s)h(s)\dif s
		= A^{\frac{t}{\delta}+1}f(0) + \int_0^t A^{\frac{t-s}{\delta}}\phi(t-s)h(s)\dif s.
	\end{equation}
	Let $\lambda > \frac{1}{\delta}\ln(A)$ and set $\kappa = \lambda - \frac{1}{\delta}\ln(A) > 0$.
	Multiply \eqref{eq:ineq3} by $e^{-\lambda t}$, then integrate from 0 to $T$ to get
	\begin{equation} \label{eq:ineq4}
		\begin{split}
			\int_0^T e^{-\lambda t}f(t)\dif t &\leq \frac{A}{\kappa}f(0) 
			+ \int_0^T\int_0^t e^{-\kappa(t-s)}\phi(t-s)e^{-\lambda s}h(s)\dif s \dif t\\
			&= \frac{A}{\kappa}f(0) 
			+ \int_0^T\int_0^{T-s} e^{-\kappa t}\phi(t)e^{-\lambda s}h(s)\dif t \dif s.
		\end{split}
	\end{equation}
	We now observe that
	\begin{equation} \label{eq:int phi}
		\begin{split}
			\int_0^\infty e^{-\kappa t}\phi(t)\dif t
			&=
			\sum_{n=0}^\infty \int_{n\delta}^{(n+1)\delta} e^{-\kappa t}(t-n\delta)^{-1/2}\dif t\\
			&= \sum_{n=0}^\infty e^{-n\kappa \delta}\int_{0}^{\delta} e^{-\kappa t}t^{-1/2}\dif t\\
			&\leq \frac{1}{1 - e^{-\kappa \delta}}\int_0^\delta t^{-1/2}\dif t\\
			&= \frac{2\delta^{1/2}}{1 - e^{-\kappa \delta}}
			= \frac{2\delta^{1/2}}{1 - Ae^{-\lambda \delta}}
		\end{split}
	\end{equation}
	Applying \eqref{eq:int phi} to \eqref{eq:ineq4}, we get
	\begin{equation} \label{eq:ineq5}
		\int_0^T e^{-\lambda t}f(t)\dif t \leq \frac{A}{\kappa}f(0) 
		+
		\frac{2\delta^{1/2}}{1 - Ae^{-\lambda \delta}}\int_0^T e^{-\lambda s}\del{B f(s) + g(s)} \dif s,
	\end{equation}
	which implies \eqref{eq:estimate}.
\end{proof}

\subsection{Proofs of Results from Section \ref{sec:interior estimates}}

\subsubsection{Proofs of Results from Section \ref{sec:interior heat}}

\label{ap:proofs interior heat}

\begin{proof}[Proof of Proposition \ref{pr:potentials}]
	\firststep
	For a fixed $x > 0$ set $z = d(x)/2$.
	We have chosen $z$ so that $x-y \geq z$ for all $y \in [0,z]$.
	Integrate by parts $n$ times to get
	\begin{multline}
		\dpd[n]{u}{x}(x,t)\\ = \int_{0}^z \dpd[n]{S}{x}(x-y,t)u_0(y)\dif y
		+ \sum_{j=1}^n \dpd[n-j]{S}{x}(x-z,t)u_0^{(j-1)}(z)
		+ \int_{z}^\infty S(x-y,t)u_0^{(n)}(y)\dif y.
	\end{multline}
	Now multiply by $z^n$:
	\begin{multline}
		z^n\dpd[n]{u}{x}(x,t)\\ = \frac{1}{z}\int_{0}^z z^{n+1}\dpd[n]{S}{x}(x-y,t)u_0(y)\dif y
		+ \sum_{j=0}^{n-1} z^{n-j}\dpd[n-j-1]{S}{x}(x-z,t)z^{j}u_0^{(j)}(z)
		+ \int_{z}^\infty S(x-y,t)z^n u_0^{(n)}(y)\dif y.
	\end{multline}
	By Corollary \ref{lem:m_n} and the fact that $S(x-\cdot,t)$ is a density, we get
	\begin{equation} 
		\abs{z^n\dpd[n]{u}{x}(x,t)}
		\leq m_n \enVert{u_0}_0 + \sum_{j=0}^{n-1} m_{n-j-1}\enVert{d^j u_0^{(j)}}_0
		+ \enVert{d^n u_0^{(n)}}_0.
	\end{equation}
	Taking the supremum over all $x$, we get
	\begin{equation} \label{eq:potential1}
		\enVert{d^n\dpd[n]{u}{x}(\cdot,t)}_0
		\leq 2^n\del{m_n \enVert{u_0}_0 + \sum_{j=0}^{n-1} m_{n-j-1}\enVert{d^j u_0^{(j)}}_0
			+ \enVert{d^n u_0^{(n)}}_0}.
	\end{equation}
	
	\nextstep
	We proceed similarly to estimate $v$, but first we define
	\begin{equation*}
		F(y,s) := \int_0^{y} f(\xi,s)\dif \xi.
	\end{equation*}
	By integration by parts we have
	\begin{equation*}
		v(x,t) = \int_0^t \int_0^\infty \dpd{S}{x}(x-y,t-s)F(y,s)\dif y\dif s.
	\end{equation*}
	Calculating as before, we get
	\begin{multline} \label{eq:v expanded}
		z^n\dpd[n]{v}{x}(x,t) = \int_0^t\frac{1}{z}\int_{0}^z z^{n+1}\dpd[n+1]{S}{x}(x-y,t-s)F(y,s)\dif y \dif s\\
		+ \int_0^t\sum_{j=0}^{n-1} z^{n-j}\dpd[n-j]{S}{x}(x-z,t-s)z^j\dpd[j]{F}{x}(z,s)\dif s
		+ \int_0^t\int_{z}^\infty \dpd{S}{x}(x-y,t)z^n\dpd[n]{F}{x}(y,s)\dif y\dif s.
	\end{multline}
	\begin{comment}
	Note that
	\begin{equation} \label{eq:int Sx estimate}
	\int_0^\infty \del{\abs{S_x(x-y,t)} +\abs{S_x(x+y,t)}} \dif x
	= 2S(0,t) = 2(2\sigma^2 \pi t)^{-1/2} =: m_{0,\sigma} t^{-1/2}.
	\end{equation}
	To prove \eqref{eq:int Sx estimate}, observe that $\abs{S_x(x,t)} = -S_x(\abs{x},t) = S_x(-\abs{x},t)$, which we use to get
	\begin{equation*}
	\int_0^\infty \abs{S_x(x-y,t)}\dif x
	= \int_0^y S_x(x-y,t)\dif x - \int_y^\infty S_x(x-y,t)\dif x
	= 2S(0,t) - S(-y,t)
	\end{equation*}
	and
	\begin{equation*}
	\int_0^\infty \abs{S_x(x+y,t)}\dif x
	=
	-\int_0^\infty S_x(x+y,t)\dif x
	= S(y,t).
	\end{equation*}
	Add these together to get \eqref{eq:int Sx estimate}.
	\end{comment}
	Now applying Corollary \ref{lem:m_n} in \eqref{eq:v expanded}, we get
	\begin{multline}
		\abs{z^n \dpd[n]{v}{x}(x,t)}
		\leq \del{m_{n+1,\sigma} + m_{n,1}}\int_0^t (t-s)^{-1/2}\sup_{0 \leq y \leq 1}\abs{F(y,s)}\dif s\\
		+ \sum_{j=1}^{n}\int_0^t m_{n-j,\sigma}(t-s)^{-1/2}\enVert{d^j \dpd[j-1]{f}{x}(\cdot,s)}_0.
	\end{multline}
	Thus,
	\begin{multline} \label{eq:potential2}
		\enVert{d^n \dpd[n]{v}{x}(\cdot,t)}_0
		\leq 2^n \del{m_{n+1,\sigma} + m_{n,1}}\int_0^t (t-s)^{-1/2}\sup_{0 \leq y \leq 1}\abs{\int_0^y f(\xi,s)\dif \xi}\dif s\\
		+ 2^n\sum_{j=1}^{n}\int_0^t m_{n-j,\sigma}(t-s)^{-1/2}\enVert{d^j \dpd[j-1]{f}{x}(\cdot,s)}_0.
	\end{multline}
	
	\nextstep
	Finally,
	\begin{equation}
		d(x)^n\dpd[n]{w}{x}(x,t) = -2\int_{0}^t d(x)^n\dpd[n+1]{S}{x}(x,t-s)\psi(s)\dif s.
	\end{equation}
	By induction we can establish a formula
	\begin{equation}
		\label{eq:dnSdx N=1}
		\dpd[n+1]{S}{x}(x,t) = S(x,t)\sum_{j=0}^{\floor{\frac{n+1}{2}}} (\sigma^2 t)^{j-n-1} c_{n+1,j} x^{n+1-2j},
	\end{equation}
	where $c_{n,j}$ are  coefficients defined recursively with respect to $n$.
	Multiply by $x^n$ to get, using \eqref{eq:heat kernel},
	\begin{equation}
		\label{eq:dnSdx N=1 1}
		x^n\dpd[n+1]{S}{x}(x,t) = (2 \pi)^{-1/2}(\sigma^2 t)^{-3/2}x\exp\cbr{-\frac{\abs{x}^2}{2\sigma^2 t}}\sum_{j=0}^{\floor{\frac{n+1}{2}}} c_{n+1,j} \del{\frac{x^2}{\sigma^2 t}}^{n-j}
	\end{equation}
	and thus
	\begin{equation}
		\int_0^\infty \abs{x^n\dpd[n+1]{S}{x}(x,t)}\dif t
		\leq \sum_{j=0}^{\floor{\frac{n+1}{2}}}\abs{c_{n+1,j}} \int_0^\infty (2 \pi)^{-1/2}(\sigma^2 t)^{-3/2}x\del{\frac{x^2}{\sigma^2 t}}^{n-j}\exp\cbr{-\frac{\abs{x}^2}{2\sigma^2 t}}\dif t.
	\end{equation}
	Use the substitution $t = \frac{x^2}{\sigma^2 s}$ to get
	\begin{multline}
		\int_0^\infty (2 \pi)^{-1/2}(\sigma^2 t)^{-3/2}x\del{\frac{x^2}{\sigma^2 t}}^{n-j}\exp\cbr{-\frac{\abs{x}^2}{2\sigma^2 t}}\dif t\\
		= \sigma^{-2}\int_0^\infty (2 \pi)^{-1/2}s^{n-j-1/2}\exp\cbr{-\frac{s}{2}} \dif s < \infty.
	\end{multline}
	We deduce that for some constant $\iota_n$, not depending on $x$,
	\begin{equation}
		2\int_0^\infty \abs{x^n\dpd[n+1]{S}{x}(x,t)}\dif t
		\leq \iota_n,
	\end{equation}
	and thus
	\begin{equation} \label{eq:potential3}
		\enVert{d^n \dpd[n]{w}{x}(\cdot,t)}_0 \leq \iota_n \sup_{0 \leq s \leq t}\abs{\psi(s)}.
	\end{equation}
	The estimates \eqref{eq:potential1}, \eqref{eq:potential2}, and \eqref{eq:potential3} result in \eqref{eq:potential estimates}.
\end{proof}

\begin{proof}[Proof of Theorem \ref{thm:dirichlet estimate}]
	Define
	\begin{equation} \label{eq:potentials dir}
		\begin{split}
			u_1(x,t) &= \int_{0}^\infty S(x-y,t)u_0(y)\dif y,\\
			u_2(x,t) &= \int_{0}^t\int_{0}^\infty S(x-y,t-s)f(y,s)\dif y\dif s,\\
			u_3(x,t) &= -2\int_{0}^t \dpd{S}{x}(x,t-s)\del{\psi(s) - u_1(0,s) - u_2(0,s)}\dif s.
		\end{split}
	\end{equation}
	Then by classical arguments (cf.~\cite[Section IV.1]{ladyzhenskaia1968linear}) we see that $u = u_1 + u_2 + u_3$ is a solution to \eqref{eq:dirichlet}.
	By the maximum principle, this solution is unique.
	
	By Proposition \ref{pr:potentials}, we have
	\begin{equation} \label{eq:potential estimates dir}
		\begin{split}
			\enVert{u_1(\cdot,t)}_n
			&\leq M_n \enVert{u_0}_n,\\
			\enVert{u_2(\cdot,t)}_n
			&\leq M_n\int_0^t (t-s)^{-1/2}\enVert{f(\cdot,s)}_{n-1,1}^*\dif s,\\
			\enVert{w(\cdot,t)}_n &\leq M_n \sup_{0 \leq s \leq t}\abs{\psi(s) - u_1(0,s) - u_2(0,s)}.
		\end{split}
	\end{equation}
	It also follows from Proposition \ref{pr:potentials} that
	\begin{equation} \label{eq:potential estimates dir2}
		\begin{split}
			\sup_{0 \leq s \leq t}\abs{u_1(0,s)} &\leq  M_n\enVert{u_0}_n,\\
			\sup_{0 \leq s \leq t} \abs{u_2(0,s)} &\leq \sup_{0 \leq s \leq t}M_n\int_0^s (s-s')^{-1/2}\enVert{f(\cdot,s')}_{n-1,1}^*\dif s'
			= 2M_nt^{1/2}\sup_{0 \leq s \leq t}\enVert{f(\cdot,s)}_{n-1,1}^*.
		\end{split}
	\end{equation}
	Combining \eqref{eq:potential estimates dir} and \eqref{eq:potential estimates dir2}, then modifying the constant $M_n$, we deduce \eqref{eq:dirichlet estimate}.
\end{proof}

\subsubsection{Proofs of Results from Section \ref{sec:interior mfg}}

\label{ap:proofs interior mfg}

Let $(u,m)$ be the solution to the finite or infinite time-horizon problem, i.e.~to System \eqref{eq:mfg T} or \eqref{eq:mfg infty}.
For a finite time-horizon we assume $u(x,T) = u_T(x)$ satisfies Assumption \ref{as:uT}.
In addition, we will impose that $\enVert{u_T}_{\s{C}^n} \leq \tilde C_n$ for each $n = 1,2,\ldots$.
(For $n = 1,2$, this is not a new assumption. For larger $n$, it is always possible to impose this restriction at the same time as Assumption \ref{as:uT}.)
We again take Assumption \ref{as:epsilon(t)}, and we denote $\epsilon = \epsilon(0)$.

If $H$ is $n+1$ times differentiable, then, 
under Assumption \ref{as:H must be smooth}, by Corollary \ref{cor:smoothness} we have
\begin{equation} \label{eq:C_ell}
	C_\ell := \max_{0 \leq \tilde \epsilon \leq \epsilon, 0 \leq Q \leq \bar Q, 0 \leq a \leq M} 
	\abs{\dpd[\ell+1]{H}{a}(\tilde{\epsilon},Q,a)} < \infty
	\quad
	\forall \ell \leq n,
\end{equation}
where $\bar Q$ is given by \eqref{eq:M}, $M$ is given in Lemma \ref{lem: a priori HJ}, and $c_2$ is the constant from Assumption \ref{as:uT} and can be made arbitrarily small.
In particular, by Corollary \ref{cor:q* a priori bound}, we have that $C_0$ can be made arbitrarily close to $\bar Q$.
By the a priori bounds proved in Section \ref{sec:fwdbckwd} (see Theorem \ref{thm:exist mfg infty}), we have the following point-wise bound:
\begin{equation*}
	\abs{\dpd[\ell+1]{H}{a}\del{\epsilon,Q^*(t),\dpd{u}{x}}} \leq C_\ell.
\end{equation*}

\begin{proposition} \label{pr:dudxn estimates}
	Let $(u,m)$ be the solution to the mean field games system on a finite or infinite time horizon $T$, i.e.~either of System \eqref{eq:mfg T} or \eqref{eq:mfg infty}.
	Suppose \eqref{eq:r assm n} holds.
	Then for any $n$ such that $H$ is $n+1$ times differentiable, we have
	\begin{equation} \label{eq:uxn}
		\sup_{t \in [0,T]} \enVert{\dpd{u}{x}(\cdot,t)}_n \leq B_n(r)
	\end{equation}
	where $B_n(r)$ is a decreasing function of $r$ that depends on the constants $C_\ell$ for $\ell = 0,1,\ldots, n$.
\end{proposition}
\begin{proof}
	Assume first that $(u,m)$ solves the finite horizon problem.
	We proceed by induction.
	In the first step we prove the base case $n = 1$, and in the second step we prove the inductive step.
	In the final step we extend the result to the infinite-horizon case.
	Note that, by taking $c_2$ small enough in \eqref{eq:C_ell}, the condition \eqref{eq:r assm n} implies
	\begin{equation*}
		r > \max\cbr{(2C_0M_n)^2,1} \ln(2M_n).
	\end{equation*}
	
	\firststep
	Define
	\begin{equation}
		w(x,t) = e^{rt}\dpd{u}{x}(x,T-t), \quad
		f(x,t) = e^{rt}\dpd{}{x}\del{H\del{\epsilon(T-t),Q^*(T-t),\dpd{u}{x}(x,T-t)}}.
	\end{equation}
	Then $w$ satisfies
	\begin{equation}
		\dpd{w}{t} = \dfrac{\sigma^2}{2}\dpd[2]{w}{x} + f(x,t).
	\end{equation}
	\begin{comment}
	By 
	\end{comment}
	
	We first calculate
	\begin{equation} \label{eq:f int}
		\begin{split}
			&\abs{\int_0^x f(y,t)\dif y}\\
			&= e^{rt}\abs{H\del{\epsilon(T-t),Q^*(T-t),\dpd{u}{x}(x,T-t)}
				- H\del{\epsilon(T-t),Q^*(T-t),\dpd{u}{x}(0,T-t)}}\\
			&\leq 2e^{rt}H(0,0,0),
		\end{split}
	\end{equation}
	using the fact that $H$ is decreasing in all its variables (Lemma \ref{lem:opt quant}).
	Next, since
	\begin{equation*}
		f(x,t) = \dpd{H}{a}\del{\epsilon(T-t),Q^*(T-t),\dpd{u}{x}(x,T-t)}\dpd{w}{x}(x,t),
	\end{equation*}
	we have
	\begin{equation} \label{eq:f0}
		\abs{d(x)f(x,t)}
		\leq C_0\enVert{w(\cdot,t)}_1.
	\end{equation}
	By \eqref{eq:f int} and \eqref{eq:f0}, we deduce
	\begin{equation} \label{eq:f01}
		\enVert{f(\cdot,t)}_{0,1}^* \leq C_0 \enVert{w(\cdot,t)}_1 + 2e^{rt}H(0,0,0).
	\end{equation}
	We also know that $\abs{w(0,t)} \leq Me^{rt}$.
	Now we apply Theorem \ref{thm:dirichlet estimate} to get
	\begin{equation}
		\enVert{w(\cdot,t)}_1 \leq M_1\del{\enVert{w(\cdot,t_0)}_1 + (t-t_0)^{1/2} C_0\sup_{t_0 \leq s \leq t} \enVert{w(\cdot,s)}_1 + A_1 e^{rt}}
	\end{equation}
	for all $0 \leq t_0 \leq t \leq t_0 + 1$, where
	\begin{equation*}
		A_1 := 2H(0,0,0) + M,
	\end{equation*}
	which can be made arbitrarily close to $2H(0,0,0) + M$.
	Set $\delta = \min\cbr{(2C_0 M_1)^{-2},1}$. 
	Then for any $0 \leq t_0 \leq t \leq t_0 + \delta$, we deduce
	\begin{equation} \label{eq:recurrence}
		\sup_{t_0 \leq s \leq t} \enVert{w(\cdot,s)}_1 \leq 2M_1\del{\enVert{w(\cdot,t_0)}_1 + A_1 e^{rt}}.
	\end{equation}
	By using \eqref{eq:recurrence} repeatedly, we deduce
	\begin{equation}
		\begin{split}
			\enVert{w(\cdot,t)}_1 &\leq (2M_1)^{\floor{\frac{t}{\delta}} + 1}\enVert{w(\cdot,0)}_1 + \sum_{j=0}^{\floor{\frac{t}{\delta}}} (2M_1)^{j+1} A_1 e^{r(t-j\delta)}\\
			&= (2M_1)^{\floor{\frac{t}{\delta}} + 1}\enVert{u_T'}_1 
			+ 2M_1A_1e^{rt}\frac{1 - (2M_1e^{-r\delta})^{\floor{\frac{t}{\delta}} + 1}}{1 - 2M_1e^{-r\delta}}.
		\end{split}
	\end{equation}
	We use the assumption
	\begin{equation} \label{eq:r assm1}
		r > \frac{\ln(2M_1)}{\delta} = \max\cbr{(2C_0M_1)^2,1} \ln(2M_1)
	\end{equation}
	and divide by $e^{rt}$ to deduce
	\begin{equation}
		\enVert{\dpd{u}{x}(\cdot,T-t)}_1 \leq 2M_1 \enVert{u_T}_1 + \frac{2M_1 A_1}{1 - 2M_1 e^{-r(2C_0 M_1)^{-2}}}.
	\end{equation}
	and since $\enVert{u_T'}_1 \leq \enVert{u_T}_{\s{C}^2} \leq \tilde C_2$ we deduce
	\begin{equation} \label{eq:ux1}
		\sup_{t \in [0,T]}\enVert{\dpd{u}{x}(\cdot,t)}_1 \leq 2M_1 \tilde C_2 + \frac{2M_1 A_1}{1 - 2M_1 e^{-r(2C_0 M_1)^{-2}}} =: B_1(r),
	\end{equation}
	which is the base case.
	
	\nextstep
	Suppose for now that \eqref{eq:uxn} holds for $n-1$; we will prove it holds for $n$.
	By using the chain and product rules, we have
	\begin{multline*}
		\dpd[m-1]{f}{x}(x,t) = e^{rt}\dpd[m]{}{x}\del{H\del{\epsilon(T-t),Q^*(T-t),\dpd{u}{x}(x,T-t)}}\\
		= e^{rt}\sum_{\ell = 0}^{m-1} \sum_{1 \leq k_\ell < k_{\ell - 1} < \cdots < k_1 < k_0 = m}
		\prod_{j=0}^{\ell - 1} {k_j - 1 \choose k_{j+1}}\dpd[k_j - k_{j+1}+1]{u}{x}(x,T-t)\dpd[k_\ell+1]{u}{x}(x,T-t)\\ \times \dpd[\ell+1]{H}{a}\del{\epsilon(T-t),Q^*(T-t),\dpd{u}{x}(x,T-t)} \quad \forall m = 1,\ldots,n,
	\end{multline*}
	where we interpret an empty product as equal to 1.
	Then using Equation \eqref{eq:C_ell} we have
	\begin{multline*}
		\abs{d^m(x)\dpd[m-1]{f}{x}(x,t)} \leq\\
		e^{rt}\sum_{\ell = 0}^{m-1} \sum_{1 \leq k_\ell < \cdots < k_0 = m}
		C_\ell\prod_{j=0}^{\ell - 1} {k_j - 1 \choose k_{j+1}}\abs{d(x)^{k_j - k_{j+1}}\dpd[k_j - k_{j+1}+1]{u}{x}(x,T-t)d(x)^{k_\ell}\dpd[k_\ell+1]{u}{x}(x,T-t)}
		\\
		\leq C_0  \enVert{w(\cdot,t)}_m
		+ e^{rt}\sum_{\ell = 1}^{m-1} \sum_{1 \leq k_\ell < \cdots < k_0 = m}
		C_\ell\prod_{j=0}^{\ell - 1} {k_j - 1 \choose k_{j+1}}\enVert{\dpd{u}{x}(\cdot,T-t)}_{k_j - k_{j+1}}\enVert{\dpd{u}{x}(\cdot,T-t)}_{k_\ell}.
	\end{multline*}
	We deduce that there exists some constant $A_n(r)$, depending only on $C_\ell$ and $B_\ell(r)$ for $\ell \leq n-1$ as well as the constant appearing in estimate \eqref{eq:f int}, such that
	\begin{equation}
		\enVert{f(\cdot,t)}_{n-1,1}^* \leq C_0  \enVert{w(\cdot,t)}_n + (A_n(r) - M) e^{rt}.
	\end{equation}
	Since $B_\ell(r)$ is decreasing with respect to $r$ for $\ell \leq n-1$, the same holds for $A_n(r)$.
	We apply Theorem \ref{thm:dirichlet estimate} again to get
	\begin{equation}
		\enVert{w(\cdot,t)}_2 \leq M_n\del{\enVert{w(\cdot,t_0)}_n + (t-t_0)^{1/2} C_0 \sup_{t_0 \leq s \leq t} \enVert{w(\cdot,s)}_n + A_n(r) e^{rt}}
	\end{equation}
	for all $0 \leq t_0 \leq t \leq t_0 + 1$.
	We will now use the assumption \eqref{eq:r assm n},
	and the exactly same argument as before yields
	\begin{equation} \label{eq:ux2}
		\sup_{t \in [0,T]}\enVert{\dpd{u}{x}(\cdot,t)}_n \leq 2M_n \tilde C_{n+1} + \frac{2M_n A_n(r)}{1 - 2M_n e^{-r(2C_0  M_n)^{-2}}} =: B_n(r).
	\end{equation}
	Since $A_n(r)$ is decreasing with respect to $r$, so is $B_n(r)$.

	For the infinite horizon case, if $(u^T,m^T)$ denotes the solution to the finite time-horizon problem, then its limit as $T \to \infty$ is the solution $(u,m)$ to System \eqref{eq:mfg infty}.
	We deduce that $(u,m)$ satisfies \eqref{eq:uxn}, with $[0,T]$ replaced by $\intco{0,\infty}$.
\end{proof}

As a corollary, we derive \eqref{eq:dHda_n} and \eqref{eq:dHdQ_n}.
To prove \eqref{eq:dHda_n}, observe that
\begin{multline*}
	\dpd[n]{}{x}\del{\dpd{H}{a}\del{\epsilon,Q^*(t),\dpd{u}{x}(x,t)}}\\
	= \sum_{\ell = 0}^{n-1} \sum_{1 \leq k_\ell < \cdots < k_1 < k_0 = n}
	\prod_{j=0}^{\ell - 1} {k_j - 1 \choose k_{j+1}}
	\dpd[k_j - k_{j+1} + 1]{u}{x}(x,t)
	\dpd[k_\ell + 1]{u}{x}(x,t)
	\dpd[\ell + 2]{H}{a}\del{\epsilon,Q^*(t),\dpd{u}{x}(x,t)},
\end{multline*}
so that
\begin{equation*}
	\begin{split}
		&\abs{d(x)^n\dpd[n]{}{x}\del{\dpd{H}{a}\del{\epsilon,Q^*(t),\dpd{u}{x}(x,t)}}}\\
		&\quad \leq \sum_{\ell = 0}^{n-1} \sum_{1 \leq k_\ell < \cdots < k_1 < k_0 = n}
		C_{\ell+1} \prod_{j=0}^{\ell - 1} {k_j - 1 \choose k_{j+1}}
		\abs{d(x)^{k_j - k_{j+1}}\dpd[k_j - k_{j+1} + 1]{u}{x}(x,t)
			d(x)^{k_\ell}\dpd[k_\ell + 1]{u}{x}(x,t)}\\
		&\quad \leq \sum_{\ell = 0}^{n-1} \sum_{1 \leq k_\ell < \cdots < k_1 < k_0 = n}
		C_{\ell+1} \prod_{j=0}^{\ell - 1} {k_j - 1 \choose k_{j+1}}
		\enVert{\dpd{u}{x}(\cdot,t)}_{k_j - k_{j+1}}
		\enVert{\dpd{u}{x}(\cdot,t)}_{k_\ell}.
	\end{split}
\end{equation*}
Thus \eqref{eq:dHda_n} follows from \eqref{eq:uxn}.
The proof of \eqref{eq:dHdQ_n} is similar: use the formulas \eqref{eq:dHdQ} and \eqref{eq:dHda}, taking successive derivatives and applying Equation \eqref{eq:dHda_n}.

\end{document}